\documentclass{article}
\usepackage[utf8]{inputenc}
\usepackage[a4paper]{geometry}
\usepackage{graphicx}
\usepackage[textwidth=8em,textsize=small]{todonotes}
\usepackage{amsmath}
\usepackage{natbib}
\usepackage{amsthm}
\usepackage{amsmath,hyperref}
\usepackage{multirow}
\usepackage{bm}
\usepackage{natbib}
\usepackage{amsfonts}
\usepackage{amsmath}
\usepackage{amssymb}
\usepackage[plain,noend]{algorithm2e}
\usepackage{graphicx,subfigure}
\usepackage{xr}
\usepackage{color}
\usepackage{rotating}
\usepackage{comment}
\usepackage{multicol}
\usepackage{float}
\usepackage{algorithmic}
\newtheorem{theorem}{Theorem}

\newtheorem{assum}{Assumption}
\newtheorem{lemma}{Lemma}

\newtheorem{corollary}{Corollary}
\newtheorem{rmk}{Remark}
\newtheorem{prop}{Proposition}

 \newcommand{\M}{\mathcal{M}}

\newcommand{\diag}{{\text{diag}}}
\newcommand{\bmu}{{\boldsymbol\mu}}

\newcommand{\ba}{\mathbf{a}}
\newcommand{\br}{\mathbf{r}}
\newcommand{\bx}{\mathbf{x}}

\newcommand{\bw}{\mathbf{w}}
\newcommand{\bZ}{\mathbf{Z}}

\newcommand{\bxi}{\boldsymbol{\xi}}
\newcommand{\bzeta}{\boldsymbol{\zeta}}

\def\bX{\mathbf{X}}

\def\bW{\mathbf{W}}
\def\bP{\mathbf{P}}
\def\bB{\mathbf{B}}

\def\bA{\mathbf{A}}
\def\bQ{\mathbf{Q}}

\def\mR{\mathbb{R}}
\def\bu{\mathbf{u}}
\def\bU{\mathbf{U}}
\def\bI{\mathbf{I}}
\def\bG{\mathbf{G}}
\def\bM{\mathbf{M}}
\def\bg{\mathbf{g}}
\def\bR{\mathbf{R}}
\def\bq{\mathbf{q}}
\def\bv{\mathbf{v}}
\def\bV{\mathbf{V}}
\def\be{\mathbf{e}}

\def\bO{\mathbf{O}}

\def\bH{\mathbf{H}}
\def\bh{\mathbf{h}}

\def\bS{\mathbf{S}}
\def\bT{\mathbf{T}}

\def\bD{\mathbf{D}}
\def\bN{\mathbf{N}}

\def\bmu{\boldsymbol{\mu}}

\def\bSigma{\boldsymbol{\Sigma}}

\def\mC{\mathbb{C}}

\def\CC{\mathcal{C}}

\def\VV{\mathcal{V}}
\def\trace{\mbox{tr}}
\def\diag{\mbox{diag}}

\def\Exp{\mbox{E}}

\DeclareMathOperator*{\argmin}{arg\,min}

\definecolor{zx}{RGB}{0,225,0}

\title{%On the
{Asymptotic limits of spiked eigenvalues and eigenvectors of signal-plus-noise matrices with weak signals and heteroskedastic noise}}
\date{}
\begin{document}
\maketitle
\vspace{-4em}
{\centering {\bf Xiaoyu Liu$^1$, Yiming Liu$^1$, Guangming Pan$^{2}$, Lingyue Zhang$^3$, 
Zhixiang Zhang$^{4}$ }\\} 
{\centering 
{\it $^1$School of Economics, Jinan University \\
$^2$ School of Physical and Mathematical Sciences, Nanyang Technological University \\
$^3$School of Mathematical Science, Capital Normal University \\
$^4$Department of Mathematics, Faculty of Science and Technology, University of Macau 
 } \\}

\vspace{1em}
\begin{abstract}
   This paper is to study a signal-plus-noise model in high dimensional settings when the dimension and the sample size are comparable. Specifically, we assume that the noise has a general covariance matrix that allows for heteroskedasticity, and that the deterministic signal has the same magnitude as the noise and can have a rank that tends to infinity. %The key technical contribution is to characterize
We develop the asymptotic limits of the
left and right spiked singular vectors of the signal-plus-noise data matrix and the limits of the spiked eigenvalues of the corresponding Gram matrix. As an application, we propose a new criterion to estimate the number of clusters in clustering problems.

    %Under  mild conditions with respect to the leading eigenvalue of the underlying covariance matrix and the noises, we investigate the asymptotic properties of the spiked eigenvalues and the corresponding eigenvalues of the sample covariance matrix. Based on the discovered results, 
    
    %some related applications  are also explored. Specifically, for a general mixture model, a new criterion to estimate the number of clusters is proposed; the properties of spectral clustering is also investigated. In addition, some classification and dimension reduction problems are also considered.  
      ~ \\~\\
    \noindent\textbf{keywords}: signal-plus-noise matrices; spiked eigenvalues and eigenvectors; deterministic equivalents; spectral clustering.
\end{abstract}

\section{Introduction}\label{sec:intro}
Consider a $p\times n$ signal-plus-noise model with the form% of
\begin{eqnarray}\label{datamat}
	\bX_n=\bA_n+\bSigma_n^{1/2}\bW_n, %\in\mR^{p\times n},
\end{eqnarray}
where $\bA_n$ is a deterministic signal matrix, $\bW_n$ consists of i.i.d.~random variables, and $\bSigma_n$ accounts for the covariance structure in the noise.  Such a model is popular in the study of machine learning \citep{yang2016rate}, matrix denoising \citep{nadakuditi2014optshrink} or signal processing \citep{vallet2012improved}. 
%When  $\bSigma_n$ is an identity matrix, there has been a huge amount of work on eigenvalues and eigenvectors for such signal-plus-noise type matrices. %To name a few, 
%\cite{loubaton2011almost} derived the almost sure limits of eigenvalues. \cite{benaych2012singular} used different approaches to study the phase transition of its singular values and vectors assuming a bi-unitarily invariant distribution on $\bW_n$. \cite{ding2020high} obtained the limits and convergent rates of the leading eigenvalues and eigenvectors. \cite{bao2021singular} showed the distributions of the principal singular vectors and singular subspaces.
When $\bSigma_n$ is an identity matrix, there has been a huge amount of work on eigenvalues and eigenvectors for such signal-plus-noise type matrices. For instance, \cite{loubaton2011almost} derived the almost sure limits of eigenvalues. \cite{benaych2012singular} used different approaches to study the phase transition of singular values and vectors assuming a bi-unitarily invariant distribution on $\bW_n$. Further contributions include \cite{ding2020high} who obtained the limits and convergence rates of the leading eigenvalues and eigenvectors, and \cite{bao2021singular} who showed the distributions of the principal singular vectors and singular subspaces. When $\bSigma_n$ is set to be a diagonal matrix, \cite{hachem2013bilinear} investigated 
the limiting behaviour of the
random bilinear form of the sample covariance matrix under a separable model, which includes the case of  $\bSigma_n$ being diagonal in \eqref{datamat}. 
  When the signal-to-noise ratio tends to infinity, i.e., the ratio of the spectral norm of the signal part to the noise part tends to infinity, \cite{cape2019signal} also considered the asymptotic properties of spiked eigenvectors under Model \eqref{datamat}. By imposing Gaussianity on $\bW_n$, \cite{han2021eigen} provided a
 eigen-selected spectral clustering method with theoretical justifications.

However, the assumptions that $\bSigma_n$ is an identity or diagonal matrix, and the signal-to-noise ratio tends to infinity, seem to be restricted and hard to verify in practice. In this paper,
%we aim to investigate the asymptotic properties of the eigenvalues and eigenvectors of the sample covariance matrix under Model \eqref{datamat} in the regime of $p/n\rightarrow c>0$, 
we aim to investigate the asymptotic properties of the eigenvalues of $\bX_n \bX_n^*$, as well as both the left and right spiked singular vectors of $\bX_n$ under the regime where $p/n\rightarrow c>0$,
with mild regularity conditions on %towards 
$\bSigma_n$ and $\bA_n$, and mild moment assumptions on $\bW_n$. %Compared with the existing literature, 
To the best of our knowledge, we first systematically study the properties of the spiked eigenvalues and eigenvectors of  Model \eqref{datamat} under such mild conditions.  Specifically, we consider \begin{equation}\label{scm2}
 \bS_{n}:=\bX_n\bX_n^*=(\bA_n+\bSigma_n^{1/2}\bW_n) (\bA_n+\bSigma_n^{1/2}\bW_n)^*
 \end{equation}
and \begin{eqnarray}\label{step2}
		\tilde\bS_{n}:=\bX_n^* \bX_n=(\bA_n+\bSigma_n^{1/2}\bW_n)^* (\bA_n+\bSigma_n^{1/2}\bW_n).
\end{eqnarray} 
In order to obtain the asymptotic properties of spiked eigenvectors of $\bS_n$ and $\tilde\bS_{n}$, we analyze the quadratic forms involving the resolvents $\bQ_n(z)$ and $\tilde \bQ_n(z)$ of matrices $ \bS_{n}$ and  $\tilde\bS_{n}$ defined as 
 \begin{eqnarray}\label{res22}
	  \bQ_n(z)=(\bS_{ n}-z\bI)^{-1}
\end{eqnarray}
and
\begin{eqnarray}\label{res}
	 \tilde \bQ_n(z)=(\tilde \bS_{ n}-z\bI)^{-1},
\end{eqnarray}
respectively, where $z$ lies in the complex upper half plane and $\bI$ refer to an identity matrix with comparable sizes. The study on the quadratic forms involving $\bQ_n$ and $\tilde \bQ_n$ can be traced back to \cite{hachem2013bilinear}, where the signal part can have arbitrary rank, and the noise component is of the form  $\bSigma_1^{1/2}\bW_n\bSigma_2^{1/2}$ where $\bSigma_1$ and $\bSigma_2$ are diagonal. 

%By the discovered limits of the quadratic forms (e.g., $\bu^*Q_n(z)\bv$) or the normalized trace (e.g., $\trace Q_n(z)/p$) in terms of  \eqref{res22} and \eqref{res}, asymptotic properties of eigenvalues and eigenvectors of \eqref{scm2} and \eqref{step2} can be obtained.
%Moreover, to demonstrate the usefulness of the derived theoretical results, below are some potential related applications:

To demonstrate the use of the theoretical results, we consider applications in spectral clustering.
When each column of $\bA_n$ can be only chosen from a finite number of distinct unknown deterministic vectors, % \eqref{datamat} can be regarded as a mixture model. 
\eqref{datamat} can be regarded as a collection of samples generated from a mixture model. 
Thus, in a vector form, the $i$-th column of Model \eqref{datamat} can be written as \begin{equation}\label{mo2} \mathbf{x}_{i}=\mathbf{a}_{i}+\bSigma_n^{1 / 2} \mathbf{w}_{i},  \end{equation}
 where $\ba_i=\bmu_s/\sqrt{n}$ for some $s\in\{1,\ldots,K\}$ if $i\in\VV_s\subseteq\{1,\ldots,n\}$. The normalized constant $\sqrt{n}$ in $\ba_i$ is to unify the Assumption \ref{ass1} below. Here $\cup_{s=1}^k\VV_s=\{1\ldots,n\}$ and $\VV_s\cap\VV_t=\emptyset$ for any $s\neq t$, and $K$ actually refers to the number of the different distributions (i.e., clusters) in a mixture model. One should also note that the labels are unknown in clustering problems. %For a mixture model, 
 Numerous literatures investigate mixture models. In statistics, \cite{redner1984mixture} considered the clustering problem for the Gaussian mixture model in low dimensional cases, while \cite{cai2019chime} considered the high dimensional cases. Some classical techniques about clustering were also proposed in past decades; see e.g., \citep{macqueen1967some,bradley1999mathematical,rdusseeun1987clustering,maimon2005data,duda1973pattern}.
  In applied %empirical
  economics, mixture models are used to introduce unobserved heterogeneity. An important example of this setup from the econometrics literature is \cite{keane1997career}, which investigated the clustering problem in labor markets. Such models also arise in analyzing some classes of games with multiple Nash equilibria. See for example, \cite{berry2006identification}, \cite{chen2014likelihood} and others. 
  %The existed approaches with theoretical justifications  related to  \eqref{mo2} mostly relies on the underlying distribution of $\bw_i$ in \eqref{mo2}. 

To put our work into a broader background, our model belongs to the deformed random matrix models in the random matrix theory. 
The behaviours of the extreme eigenvalues are closely related to the signal strength and exhibit the so-called BBP phase transition due to the seminal work \cite{baik2005phase}. The extreme eigenvectors also exhibit different asymptotic behaviours when the associated eigenvalues are above or below the critical values.  \cite{benaych2012singular} studied the phase transition of the singular values and vectors of the low-rank perturbations of large rectangular random matrices assuming that the noise part has a bi-unitarily invariant distribution. There are fruitful results on the study of the asymptotic properties of the eigenvalues and eigenvectors of the deformed random matrix models including \citep{paul2007asymptotics,el2007tracy,capitaine2009largest,benaych2011eigenvalues, benaych2011fluctuations,knowles2013isotropic,huang2018mesoscopic,bao2021singular,ding2020high}.  Another line of work with fruitful results focuses on the non-asymptotic properties of the spectrum of deformed matrix models via matrix perturbation theory, such as  \cite{davis1970rotation,ipsen2009refined,yu2015useful,abbe2020entrywise,chen2021spectral}. Our work belongs to the first line of research, where our main objective is to derive precise asymptotic limits of the spiked eigenvalues and eigenvectors of rectangular deformed models. 
This work distinguishes itself from the aforementioned studies by allowing for heteroskedastic noise and considering a deterministic perturbation matrix with a rank that tends to infinity.

%The model studied in this work is also related to the information-plus-noise model in random matrix theory literature \cite{dozier2007empirical}, where the signal matrix is of full rank. Some recent work on this model includes \cite{zhang2024tracy} where the noise matrix has i.i.d. entries and \cite{bai2023exact} which is more general. However, none of the above investigated the spiked eigenvalues and eigenvectors.

The inference and recovery of the signals play crucial roles in many application fields and have been actively studied. 
There is a tremendous literature 
on such topics, among which including \cite{mestre2008improved,nadakuditi2014optshrink,shabalin2013reconstruction,donoho2014minimax} relies on some results on the asymptotic limits of spiked eigenvalues and eigenvectors in random matrix theory. The spiked sample eigenvalues and eigenvectors contain valuable information about the signals. For instance, the number of spiked eigenvalues 
can represent the number of communities in network data analysis or the number of groups in a mixture model as previously discussed. Additionally, the eigenvector captures the group structure of nodes in network data \cite{rohe2011spectral,fan2022asymptotic} or the label information of observations from a mixture model \cite{han2021eigen}. 
The main statistical application in this work is the inference of the number of clusters in clustering, which can be easily extended to estimating the number of spikes in signal-plus-noise models arising in other problems.

%delicate purterbation results on eigenvalues and eigenvectors \cite{yu2015useful,abbe2020entrywise,chen2021spectral}. 
%The number of spiked eigenvalues represent the number of strong signal  

  Our main theoretical contribution is to
  precisely characterize the 
 first-order limits of the eigenvalues and eigenvectors of $\bS_n$ and $\tilde{\bS}_n$. 
  There are two observations that can be obtained based on our main theoretical results that are somewhat surprising, as they exhibit some overlaps with findings in the literature, albeit in different scopes of problems.
  The first is that the limits of the spiked eigenvalues of $\bS_n$ coincide with these of the sample covariance matrices without the signal part, by letting the population covariance matrix be $\bA_n\bA_n^* + \bSigma_n$, see the discussion below Theorem \ref{thm2}.
  The second is that, 
 the spiked right singular vectors of $\bX_n$ have an intrinsic block structure if $\bA_n$ contains a finite number of distinct deterministic factors, even for a moderate signal-to-noise ratio. Our Corollary \ref{cor1} precisely quantifies the deviation of the right singular vector from a vector with entries having a group structure. This finding is highly relevant to the field of spectral clustering, which has been extensively discussed in the literature. It is worth noting that many existing studies assume strong moment conditions on the noise and consider scenarios where the signal-to-noise ratio tends to infinity.

 As applications,  we propose a method to estimate the number of clusters by leveraging the asymptotic limits of sample eigenvalues. We also discuss how the theoretical results intuitively explain why the spiked eigenvectors have clustering power in the context of spectral clustering.

The remaining sections are organized as follows. In Section \ref{secModelResults}, we state the main results for Model \ref{datamat}. Section \ref{app} includes the applications of clustering. In Section \ref{simu}, we provide numerical results related to the applications discussed in Section \ref{app}. The proofs can be found in Sections \ref{secproofs} and \ref{append}.

%The results in this paper also lead to the study of a spiked model that $\bR$ has some spiked singular values. Spiked sample eigenvalues may carry valuable information about the structure of signals which is worth to study.

%From random matrix point of view, our model can be regarded as a rectangular deformed model, which is similar to the deformed Wigner matrix for Hermitian matrices. ?  From statistical point of view, our results will be useful in detection of signals. By treating $\bR$ as mean of a data matrix, our results shed light on the properties of PCA for data with complicated mean structure. In practice of neural network, it is well known that adding noise to the input data  during training can lead to significant improvements in generalization performance  in some circumstances. Our model can be helpful to understand this phenomenon and our result in some sense explain that adding strong noise may bury the signals.
{\it Conventions}: We use $C$ to denote generic large constant independent of $n$, and its value may change from line to line. Denote $a\wedge b=\min\{a,b\}$, and $a \vee b =\max\{a,b\}$. Let $\mathbf1$ and $\bI$ refer to a vector with all entries being one and an identity matrix with a comparable size, respectively. Let $\|\cdot\|$ denote the Euclidean norm of a vector or the spectral norm of a matrix. Denote by $\mathcal{C}^+$ the complex upper half plane. 
We use $O(\cdot)$ and $o(\cdot)$ for the standard big-O and little-o notation. For a sequence of random variables $(X_n)_{n\ge 1}$, we write
$X_n  \to_P X$ to denote that $X_n$ converge in probability to $X$. For a sequence
$(a_n)_{n\ge 1}$ of scalars, we write $a_n = O_P(1)$ if $(a_n)_{n\ge 1}$ is bounded in probability and $a_n = o_P(1)$
if $(a_n)_{n\ge 1}$ converges to zero in probability. 

%For two positive quantities $A_N$ and $B_N$, $A_N \sim B_N$ means that $C^{-1} B_N \leq A_N \leq C B_N$ for some constant $C$ independent of $N$. $A_N = O(B_N)$ means that $A_N \leq C B_N$. For any function $f$ we use $f^{(n)}$ to denote the $n$-th derivative of $f$. The $s$-coordinate of a vector $\bw$ is denoted by $\bw(s)$. We use $\|\cdot\|$ to denote the spectral norm of a matrix or a vector.

\section{The main results}\label{secModelResults}
  In this section, we mainly investigate the limits of the eigenvalues and eigenvectors of $\bS_n$ and $\tilde{\bS}_n$ defined in \eqref{scm2} and \eqref{step2}, respectively. 
We first impose some mild conditions on $\bW_n$ for establishing the asymptotic limits of the eigenvalues and eigenvectors:
 \begin{assum}\label{ass1}
	We assume that $\bW_n =(w_{ij}) $ is an $p\times n$ matrix, whose entries $\{w_{ij}:1\leq i \leq p, 1\leq j \leq n\}$ are independent  real or complex  random variables satisfying
 \begin{equation*} \mathbb{E} w_{ij}=0, \quad   \mathbb{E}|\sqrt{n}w_{ij}|^2 =1 \text{ and }  \mathbb{E}|\sqrt nw_{ij}|^4\le C. 
 \end{equation*}
\end{assum}
We consider the high dimensional setting specified by the following assumption.
 \begin{assum}\label{ass2}
 $p/n\equiv c_n \rightarrow c \in (0,\infty)$.
 \end{assum}
Note that when $\bA_n$ has a rank of $o(n)$, the limiting spectral distribution of $\bX_n \bX_n^*$  is the same as that of the model by setting $\bA_n$ as a zero matrix. This can be directly concluded by the rank inequality, see Theorem A.43 of \cite{bai2010spectral}. However, to investigate the limiting behaviours of the spiked eigenvalues and eigenvectors under Model \eqref{datamat},  more assumptions on $\bA$ and %also 
$\bSigma_n$ are required. 

\begin{assum}\label{ass3new}
Let $\bA_n$ be a $p\times n$ matrix with bounded spectral norm and of rank $O(n^{1/3})${\color{red}}, and $\bSigma_n$ be a symmetric matrix with bounded spectral norm.	Let $\bR_n :=\bA_n\bA_n^* +\bSigma_n$. 
%We assume that there are 
The $\tilde{n}$ distinct eigenvalues of $\bR_n$ are denoted by $\gamma_1>\gamma_2>\cdots>\gamma_{\tilde{n}}$, %in descending order, 
$\tilde{n}\le p$. The multiplicities of these $\tilde{n}$ distinct eigenvalues are $m_1, m_2, \cdots, m_{\tilde{n}}$, respectively. The eigenvector associated with $\gamma_i$, or the eigenspace if the multiplicity is greater than one, is denoted by $\boldsymbol\Xi_i$.
\end{assum}
\begin{rmk}
In this paper, we consider the case where the leading eigenvalue of $\bR_n$ is bounded, and a similar strategy can be adapted to investigate the case of divergent spikes.  
\end{rmk}

 The key technical tool is the deterministic equivalents of
 $\bQ_n$ and $\tilde{\bQ}_n$ in \eqref{res22} and \eqref{res}. 
 We introduce it first as it requires weaker assumptions than the main results on spiked eigenvectors and eigenvalues, and may be of independent interest.
 % To introduce the deterministic equivalents of  
 %$Q_n$ and $\tilde{Q}_n$ in \eqref{res22} and \eqref{res},
 %we need more notations.
 For any $z \in \mathcal{C}^{+}$, let $\tilde{r}_n(z) \in \mathcal{C}^{+}$ be the unique solution to the equation \begin{eqnarray}\label{32ra3}z=-\frac{1}{\tilde{r}_n}+c_n\int\frac{t\text{d}H^{\bR_n}(t)}{1+t\tilde{r}_n},
\end{eqnarray}
where $H^{\bR_n}(t)$ is the empirical spectral distribution of $\bR_n$. Proposition \ref{prop1} below provides the deterministic equivalence of  $\bQ_n$ and $\tilde{\bQ}_n$, and its proof is in Section \ref{secproofs}.

\begin{prop}\label{prop1}
Suppose that Assumptions \ref{ass1} to \ref{ass3new} are satisfied. Let $(\bu_n)_{n\ge 1}, (\bv_n)_{n\ge 1}$ be sequences of deterministic vectors of unit norm. Then for any 
$z\in \CC^+$ with $\Im z$ being bounded %away 
from below by a positive constant,
%zero,
%=n^{-l}, l\in({0,1/10})$,
%and a fixed unit vectors $\bu_n\in\mR^n$ and $\bv_n\in\mR^p$.
we have \begin{eqnarray}\label{key1}
 { \mathbb{E}}|\bu_n^*(\tilde \bQ_n(z)-\tilde{\bD}_n(z))\bu_n|^2=O(n^{-1}),
\end{eqnarray} 
where \begin{equation}\label{Tz}
 \tilde{\bD}_n(z) =\tilde r_n(z)\bI-(\tilde r_n(z))^2\bA_n^*\left[\bI+\tilde r_n(z) \bR_n \right]^{-1}\bA_n
\end{equation} 
and \begin{equation}\label{Qnlim}
	 {\mathbb{E}} |\bv_n^*(\bQ_n(z)-\bD_n(z))\bv_n|^2 =O(n^{-1}),
	\end{equation}
	where \begin{equation}
\bD_n(z) = \left( -z \bI -z \tilde{r}_n \bR_n \right)^{-1}.
	\end{equation}
\end{prop}

%Proposition \ref{prop1} has some similarities to the results in \cite{hachem2013bilinear}
\begin{rmk} The model we studied is similar to that in \cite{hachem2013bilinear}, and the proof of Proposition \ref{prop1} leverages their main result, which is used in the analysis of our Gaussian case. The main difference between our Proposition \ref{prop1} and their results
is that we study the case with a general $\bSigma_n$, while they consider a model with a separable variance profile where the noise part can be written as $ \bSigma_1^{1/2} \bW_n \bSigma_2^{1/2}$ but $\bSigma_1$ and $\bSigma_2$ are both diagonal matrices. A potential direction for future research is to extend our results to the case with a separable covariance structure in the noise component, where both $\bSigma_1$ and $\bSigma_2$ are general matrices instead of being diagonal. This extension would allow for a more comprehensive understanding of the behaviors of the spiked eigenvalues and eigenvectors in a broader range of settings. %by using similar strategies in this work. 
We refer the reader to \cite{ding2021spiked} for related work.

There are two features of Proposition \ref{prop1} that are worth mentioning here. First, the deterministic equivalents of both $\bQ_n(z)$ and $\tilde{\bQ}_n(z)$ involves a quantity 
$\tilde{r}_n$, which is actually the Stieltjes transform of the generalized Marchenko-Pastur law, see \cite{bai2010spectral} for instance. 
This is hidden in \cite{hachem2013bilinear} as their results hold for general $\bA_n$ instead of being of finite rank.
Second,  when $\bX_n$ has columns with the structure specified in \eqref{mo2}, which is of statistical interest, especially in the context of special clustering, $\tilde{\bD}_n(z)$ has a block structure as it is in a form of $c_1\bI+c_2 \bA_n^* \bM \bA_n$ for some constants $c_1, c_2$ and some matrix $\bM$. This can also be inferred from the observation that $ {\mathbb{E}} (\bX_n^* \bX_n) $ can also be written in such a form.

%Third, the form of $R_n$ is the same as the deterministic equivalents for the sample covariance matrices studied in \cite{mestre2008asymptotic}.
\end{rmk}

Based on Proposition \ref{prop1}, we now first focus on the eigenvectors corresponding to the spiked eigenvalues of $\bR_n$. The following assumption is needed.
%imposed on $\gamma_i, i = 1,\cdots, K$ guarantees the existence of spiked sample eigenvalues. 
%Our key assumptions that lead to the existence of spiked sample eigenvectors and sample eigenvalues are imposed on the matrix $\bR_n=\Exp \bS_n$, 
%taking into account the influences of both $\bA$ and $\bSigma_n$ on the sample eigenvalues simultaneously.

%Denote the largest eigenvalue of $\bSigma_n$ by $\sigma_1$. Let $\tau_0 \in [0,1/\sigma_1]$ be the solution to equation \begin{eqnarray}\label{defoftau0}\int \left(  \frac{t \tau_0}{1-t \tau_0}\right)^2 dF^{\bSigma_n}(t) = 1/c_n.\end{eqnarray}

\begin{assum}\label{ass4new} 
The first $L$ distinct eigenvalues satisfy $\min_{1\leq i\leq L}(\gamma_{i}-\gamma_{i+1})>c_0>0$ for some constant $c_0$ independent of $p$ and $n$, and \begin{eqnarray*}\int\frac{t^2 \mathrm{d}H(t)}{(\gamma_k-t)^2}<\frac1c \end{eqnarray*} for $k\le L<\infty$, where $H(t)$ is the limiting spectral distribution of $\bR_n$. The multiplicities of these $L$ distinct eigenvalues satisfy $\sum_{i=1}^L {m_i} = O(n^{1/3})$. 
\end{assum}

Assumption \ref{ass4new} is a variant of the condition given in definition 4.1 of \cite{bai2012sample}, ensuring that the first $\sum_{i=1}^L {m_i}$ largest eigenvalues of $\bS_n$ are spiked eigenvalues. Note that we allow the number of spiked eigenvalues, $\sum_{i=1}^L {m_i}$, to diverge at the rate of $O(n^{1/3})$. This number is not necessarily equal to the rank of $\bA$.

Let $\mathcal{K}_i = \{\sum_{j=1}^{i-1}m_{j}+1, \sum_{j=1}^{i-1}m_{j}+2,\cdots, \sum_{j=1}^{i}m_{j}\}$ be the set of indices associated with the population eigenvalue $\gamma_i$, and we use the convention that $m_0=0$. Let $\hat{\lambda}_1\ge \hat{\lambda}_2 \ge \cdots \ge \hat{\lambda}_p$ be the sample eigenvalues of $\bS_n$, and $\hat{\bv}_1, \cdots, \hat{\bv}_p$ be the corresponding eigenvectors. We also denote the eigenvectors of $\tilde{\bS}_n$ by $\hat{\bu}_1, \cdots, \hat{\bu}_n$, associated with the sample eigenvalues of $\tilde{\bS}_n$ in the descending order.
The following theorem characterizes the asymptotic behaviours of the spiked eigenvectors and is proved in Section \ref{secproofs}.
%Let $\hat\bv_k\in\mR^{p}$ and $\hat{\bu}_k\in\mR^{n}$ be the eigenvector associated with the $k$-th largest (spiked) eigenvalue of $\bS_n$ and $\tilde{\bS}_n$, respectively. The following theorem characterizes the asymptotic behaviours of $\hat\bv_k$ and $\hat{\bu}_k$, and is proved in Section \ref{secproofs}.

%\noindent{\textbf{Theorem 1.}}
\begin{theorem}\label{thmeve1}
Suppose that Assumptions \ref{ass1}-\ref{ass4new} hold. For any $1\leq k \leq L$ and any sequences of deterministic unit vectors  $\{\bv_n\}_{n\ge 1}$ of length $p$ and $\{\bu_n\}_{n\ge 1}$ of length $n$ we have %the following two facts:
 %\begin{equation}\label{l2normvandhatu} \inf_{s \in \{1,-1\}} \|  \bv - s \hat{\bu}_k \| ^2= 2- 2 \left( \theta_k \frac{\bv^* \bA^* \xi_k \xi_k^* \bA \bv}{\gamma_k}\right)^{1/2}+ O_P (\frac{1}{\sqrt{n}}),\end{equation}
%\begin{equation}\| \bv -  \hat{\bu}_k \| ^2 \| \bv +  \hat{\bu}_k \| ^2 \end{equation}
\begin{enumerate}
\item %\begin{equation}\label{quadlim0}	|\bv^*\hat\bv_k\hat\bv_k^*\bv-\bv^* \bP_k \bv|\overset{p}\rightarrow0,\end{equation}
\begin{equation*}
\left|\bv_n^* \left(\sum_{\ell\in \mathcal{K}_k} \hat\bv_\ell\hat\bv_\ell^*\right)\bv_n-\bv_n^* \bP_k \bv_n\right|=O_P\left(\frac{1}{\sqrt{n}}\right),
\end{equation*}
where $\bP_k=\sum_{j=1}^{\tilde{n}} c_k(j)\boldsymbol\Xi_j \boldsymbol\Xi_j^*$, and $\{c_k(j)\}$ are defined by 
\begin{equation*}
c_k(j) = \left\{
\begin{array}{ll}
\displaystyle 1- \frac{1}{m_k}\sum_{i=1, i\neq k}^{\tilde{n}} m_i \left( \frac{\gamma_k}{\gamma_i-\gamma_k}-\frac{\omega_k}{\gamma_i -\omega_k}\right), & j = k\\[10pt]
\displaystyle \frac{\gamma_k}{\gamma_j-\gamma_k}-\frac{\omega_k}{\gamma_j-\omega_k}, & j\neq k\\
\end{array}
\right.
\end{equation*}
 
and  $\omega_1 \geq \omega_2 \geq,\cdots,\geq \omega_p$ are the real solutions to the equation in $\omega$:\begin{equation*}
\frac{1}{n}\sum_{i=1}^{\tilde{n}} m_i \frac{\gamma_i}{\gamma_i-\omega}=1.
\end{equation*}

\item 
\begin{equation}\label{l2normvandhatu}\left|\bu_n^* \left(\sum_{\ell\in \mathcal{K}_k}\hat{\bu}_\ell\hat{\bu}_\ell^* \right)\bu_n -\eta_k \frac{\bu_n^* \bA_n^* \boldsymbol\Xi_k \boldsymbol\Xi_k^* \bA_n \bu_n}{\gamma_k}\right| =  O_P \left(\frac{1}{\sqrt{n}}\right),\end{equation}
where  \begin{equation}\label{etak}\eta_k=\left(1-\frac1n\sum_{i=1,i\neq k}^{\tilde{n}}\frac{m_i \gamma_i^2}{(\gamma_k-\gamma_i)^2}\right).\end{equation}
\end{enumerate}
\end{theorem}

%\begin{assum}\label{ass5} There exists a small constant $\tau>0$ such that for all $n$ large enough, $F^{\bSigma_n}(0,\tau)\leq 1-\tau$ and $\limsup_{n}\sigma_1 \tau_0 <1$.
%$1+\underline{\delta}(\lambda_+)\sigma_1>\tau$.
%	\end{assum}

\begin{rmk}\label{vectormatch}
 It is worth mentioning that the first-order behaviour of the left spiked singular vectors of $\bX_n$ is the same as that of a sample covariance matrix of  $\bR_n^{1/2} \bW_n $, see the main results in \cite{mestre2008asymptotic}, and Table \ref{tabev} below demonstrated by a simulation. However, the behaviour of the right singular vectors is significantly distinct. Specifically, when the entries of $\bW_n$ are Gaussian variables, the matrix composed of the right eigenvectors of $\bR_n^{1/2}\bW_n$
 is asymptotically Haar distributed.  This observation contrasts with the second fact in Theorem \ref{thmeve1}.

 In addition, it is noteworthy that when $\bSigma_n=\bI$, the model reduces to the one studied in \citep{ding2020high,bao2021singular}. 
 In these studies, the results on the left and right singular vectors of $\bX_n$ are observed to be symmetric due to the symmetry of the model structure. However, for a general $\bSigma_n$, we cannot deduce the properties of the right singular vectors of $\bX_n$ solely based on the properties of the left singular vectors, and vice versa. We %will 
 further discuss the relationship between our results and those in \cite{ding2020high} below Theorem \ref{thm2}.
\end{rmk}

%\begin{rmk}
%Note that $\bP_k$ is a $p\times p$ matrix with eigenvalues $c_k(j), j =1, \ldots , p$ and associated eigenvectors $\xi_j$.  %The expressions of $c_k(j)$ are complicated but we can verify that they are positive,  given at the end of the proof of Corollary \ref{cor2} below. 
%\end{rmk}

The asymptotic behaviour of the spiked eigenvalues is also considered, and thus some more notations are also required.
Similar to \cite{bai2012sample}, for the spiked eigenvalue $\gamma$ outside the support of  $H$ and $\gamma\neq0$, we define
\begin{equation}\label{equ1.3}
\varphi(\gamma)=z(-\frac{1}{\gamma})=\gamma\left(1+c\int\frac{t\mathrm{d}H(t)}{\gamma-t}\right),
\end{equation}
where $z$ is regarded as the function defined in \eqref{32ra3} with its domain extended to the real line. As defined in \cite{bai2012sample}, a spiked eigenvalue $\gamma$ is called a distant spike if $\varphi'(\gamma)>0$ which is coincident to Assumption \ref{ass4new}, and a close spike if $\varphi'(\gamma)\leq0$. Note that $\bS_n$ and $\tilde\bS_n$ share the same nonzero eigenvalues, and we denote by $\hat{\lambda}_1\ge\ldots\ge\hat{\lambda}_{p\wedge n}>0.$

It is well known that $\tilde{r}_n$ defined through \eqref{32ra3} corresponds to the Stieltjes transform of a probability measure, denoted as $F^{c_n,H^{\bR_n}}$, see \cite{bai2010spectral} for instance. %The support of $F^{c_n,H^{\bR_n}}$ can be written as $\cup_{j=1}^{n'} [a_{2j-1},a_{2j}]$ where $n'\le \tilde{n}$. 
The classical eigenvalue location $\mu_i$ (\citep{erdHos2012rigidity,knowles2017anisotropic}) is defined  via 
\begin{eqnarray*}
    n \int_{\mu_i}^\infty \mathrm{d} F^{c_n,H^{\bR_n}} = i.
\end{eqnarray*}

The following theorem characterizes the asymptotic limits of eigenvalues of $\bS_n$, and is proved in Section \ref{secproofs}.
\begin{theorem}\label{thm2}
Under Assumptions \ref{ass1} to \ref{ass4new},  for $1\le k\le L$, and $\ell \in 
 \mathcal{K}_k$,  we have,
 \begin{equation}\label{spikaslim}\hat{\lambda}_\ell \rightarrow\varphi(\gamma_k) ~~~~ \text{a.s.}.\end{equation}
Moreover, for a nonspiked eigenvalue $\hat{\lambda}_j$ with $\sum_{i=1}^L m_i<j<p$, 
\begin{equation}\label{equ2.5new}
|\hat{\lambda}_{j}-\mu_j|\to 0 ~~~~ \text{a.s.}.
\end{equation}
\end{theorem}

\begin{rmk}\label{valuematch}
%Applying 
By the main results in \cite{bai2012sample}, one could obtain the limits of the spiked eigenvalues of $\bR_n^{1/2}\bW_n\bW_n^* \bR_n^{1/2}$. Theorem \ref{thm2} indicates that 
the asymptotic limits of the spiked eigenvalues of $\bX_n\bX_n^*$ are the same as those of $\bR_n^{1/2}\bW_n\bW_n^* \bR_n^{1/2}$. See Table \ref{tabev} below for an illustration.

The study on the spiked eigenvalues leverages the exact separation of the eigenvalues developed by \cite{liu2022ieee}, which extends similar phenomena in many other random matrix models pioneered by \cite{bai1998no}. However, the signal matrix $\bA$ assumed in \cite{liu2022ieee} can only have a finite number of distinct columns. In this work, we extend their main result to the case when the rank of $\bA$ can tend to infinity at the rate of $O(n^{1/3})$ without any other structural assumptions.
\end{rmk}

%\noindent \textbf{Related work.}
%\cite{ding2020high} investigated the limits of the spiked eigenvalues and eigenvectors of a signal-plus-noise model where $\bSigma_n = \bI$. \cite{bao2021singular} further obtained the fluctuation of quadratic forms of left and right spiked eigenvectors of a signal-plus-noise model where $\bSigma_n = \bI$. The model we considered includes both of these two models as special cases, and our results show that the source of sample spiked eigenvalues can be either from the spikes in the signal matrix $\bA_n$, or spikes from $\bSigma_n$, which is the covariance matrix of the noise part. 

 Our results show that the source of sample spiked eigenvalues can be either from the spikes in the signal matrix $\bA_n$, or spikes from $\bSigma_n$, which is the covariance matrix of the noise part. The model we considered includes the commonly studied additive and multiplicative perturbation models, which can be represented by $\bA_n+\bW_n$ and $\bSigma_n^{1/2}\bW_n$ in our notation, as special cases. In Section \ref{sec:example}, we provide an example where $\bSigma_n$ contains a spike to illustrate that both $\bA_n$ and $\bSigma_n$ can contribute to the sample ``spikes''. %Moreover, the overlap between the eigenspace associated with the signal and that associated with the noise can significantly influence the behavior of the sample spiked eigenvalues.

We verify that our main results can recover some of the results in \cite{ding2020high}, or earlier related work \cite{benaych2012singular} for Gaussian $\bW_n$.  %As the first $K$ eigenvectors of $\bR_n=\bA_n \bA_n^* + \bI_n$ are the same as the left singular vectors of $\bA_n$, the singular value decomposition of $\bA_n$ can be written as $\bA_n = \sum_{i=1}^r d_i \xi_i \zeta_i^*$ where $\xi_i$ denotes the eigenvector associated with the $i$-th largest eigenvalue of $\bR_n$. 
Consider $\bSigma_n = \bI_p$, and let the singular value docomposition of $\bA_n$ be $\bA_n = \sum_{i=1}^r d_i \bxi_i \bzeta_i^*$ where $d_i$ are distinct and $r$ is fixed.
Then $\gamma_k = d_k^2+1$ for $k = 1,\cdots, r$, and $\gamma_k = 1$ for $(r+1)\le k \le p$.  Theorem \ref{thm2} implies that if $d_k>c^{1/4}$, we almost surely have $$\hat{\lambda}_k\to \varphi(\gamma_k) = (d_k^2+1)(1+d_k^{-2}c).$$ 
By taking $\bu_n = \bzeta_k$ in \eqref{l2normvandhatu},  we find $\bu_n^* \bA_n^* \bxi_k =d_k$,  and  $\eta_k = 1-c_n d_k^{-4}+O(n^{-1})$, thus for $k\le r$,
 $$\bzeta_k^* \hat{\bu}_k\hat{\bu}_k^* \bzeta_k - (d_k^4-c_n)/[d_k^2(1+d_k^2)] = O_P(n^{-1/2}).$$ 
 These limits coincide with $p(d_k)$ and $a_2(d_k)$ defined in (2.6) and (2.9) of \cite{ding2020high}, respectively.

 One may wonder whether the asymptotic distributions of the spiked eigenvalues and eigenvectors of $\bX_n\bX_n^*$ are the same as those of $\bR_n^{1/2}\bW_n\bW_n\bR_n^{1/2}$, given that their first order limits coincide. The study of the latter model was pioneered by \cite{johnstone2001distribution} and with further related work including \citep{bai2012sample,jiang2021generalized,zhang2022asymptotic,bao2022statistical}. Through simulations, we observe different asymptotic variances between the two models, as indicated by Table \ref{tabev}. %Therefore it appears that the answer is no.

   %The asymptotic behaviours of spiked eigenvalues and eigenvectors of $\bSigma_n^{1/2}\bW$ when $\bSigma_n$ has some spikes has been derived in \cite{bao2022statistical}. 

The aforementioned theoretical results are all built on $\bS_n$ or $\tilde{\bS}_n$ that refer to the noncentral covariance matrices. In some situations, the centered versions are also of interest. Specifically, we consider the corresponding covariance matrices $$\mathbf{\mathcal{S}}_n=(\bX_n-\bar\bX_n)(\bX_n-\bar\bX_n)^*,$$
and $${\boldsymbol\tilde{\mathcal{S}}_n}=(\bX_n-\bar\bX_n)^*(\bX_n-\bar\bX_n),$$
	where $\bar\bX_n=\bar\bx_n\mathbf1^*$ and $\bar\bx_n=\sum_
	{k=1}^n\bx_k/n$.  Let $\boldsymbol\Phi=\bI-\mathbf1\mathbf1^*/n$ and
denote the spectral decomposition of $\bar{\bR}_n=\bA_n\boldsymbol\Phi \bA_n^*+\bSigma_n$ by $\bar{\bR}_n=\sum_{k=1}^{\tilde{n}}\bar{\gamma}_k\bar{\boldsymbol\Xi}_k\bar{\boldsymbol\Xi}_k^*$, where $\bar{\gamma}_1>\ldots> \bar{\gamma}_{\tilde{n}}$ are distinct eigenvalues and $\bar{\boldsymbol\Xi}_k$ are the associated eigenvector (or eigenspace). 
%Here $\bar K$ may be equal to $K$ or $K-1$ in some different cases.
Moreover,  define the corresponding resolvent $\mathcal{\bQ}_n(z)$  and  $\tilde{\mathcal{\bQ}}_n(z)$ of matrix $\boldsymbol{\mathcal{S}}_n$ and $\tilde{\bm{\mathcal{S}}}_n$, respectively:
\begin{eqnarray*}
	 \boldsymbol{\mathcal{Q}}_n(z)=(\bold{\mathcal{S}}_n-z\bI)^{-1},\quad \boldsymbol{\tilde{\mathcal{Q}}}_n(z)=(\bold{\tilde{\mathcal{S}}}_n-z\bI)^{-1}. 
\end{eqnarray*}
Based on the given notations, we %extend Proposition \ref{prop1}:
established the corresponding results for the centralized sample covariance matrices. The proofs are relegated to Section S.3 in the supplementary file.

\begin{prop}\label{cenprop}
	Suppose that Assumptions \ref{ass1} and \ref{ass2}  are satisfied, replace $\bR_n$ in Assumption   \ref{ass3new} by $\bar\bR_n$. Then, we have \begin{eqnarray*}
		&\left|\bu_n^*(\boldsymbol{\tilde{\mathcal{Q}}}_n(z)-\boldsymbol{\tilde{\mathcal{D}}}_n(z))\bu_n\right|=O_P(1/\sqrt n),\\
		&\left|\bv_n^*(\boldsymbol{\mathcal{Q}}_n(z)-\boldsymbol{\mathcal{D}}_n(z))\bv_n\right|=O_P(1/\sqrt n)
	\end{eqnarray*} 	
where 
%$$D(z)=\left(-z(\bI+\tilde m(z)\bSigma_n)+\bA_n(\bI+m(z)\boldsymbol\Phi)^{-1}\bA_n^*\right)^{-1}$$
%$$\tilde D(z)=\left(-z(\bI+m(z)\boldsymbol\Phi)+\bA_n^*(\bI+\tilde m(z)\Sigma)^{-1}\bA_n\right)^{-1},$$
 %$m(z)=\frac1n\trace (\Sigma D(z))$ and $\tilde m(z)=\frac1n\trace(\boldsymbol\Phi\tilde D(z))$, for any fixed unit vectors $\bu_n\in\mR^n$ and $\bv_n\in\mR^p$  and some $z\in\CC^+$ with $\Im z=n^{-l}$ and  $l>0$.
\begin{equation*}\begin{aligned}
\boldsymbol{\tilde{\mathcal{D}}}_n(z)&=\tilde{r}_n(z)\boldsymbol\Phi  - \tilde{r}_n(z)^2 \boldsymbol\Phi \bA_n^* (\bI+\tilde{r}_n(z)\bar{\bR}_n)^{-1}\bA_n\boldsymbol\Phi  -z^{-1} n^{-1}\mathbf{1}\mathbf{1}^*,\\
  \boldsymbol{\mathcal{D}}_n(z)&= (-z - z \tilde{r}_n\bar{\bR}_n)^{-1}.
\end{aligned}\end{equation*}
\end{prop}

%\begin{equation}\label{woodburrytQ}
 %   (\bX^*\bX-z\bI)^{-1} = (-z\bI)^{-1}+z^{-1}\bX^*(\bX\bX^*-z\bI)^{-1}\bX
%\end{equation}

%Denote $M=\bX^*(\bX\bX^*-z\bI)^{-1}\bX$. We know that %$\bu_1^*M\bu_2 = z \bu_1^* \tilde{Q}_n \bu_2 -\bu_1^*\bu_2$
%\begin{equation*}
%\begin{aligned}
%    &\bu^*\Phi M \Phi \bu - \bu^* M \bu = -\frac{2\bu^*\mathbf{1}}{n}\mathbf{1}^*M \bu + \frac{(\bu^*\mathbf{1})^2}{n^2}\mathbf{1}^* M \mathbf{1}\\
%    & = -\frac{2\bu^*\mathbf{1}}{n}z\mathbf{1}^*\tilde{Q}_n \bu+\frac{(\mathbf{1}^*\bu)^2}{n} +\frac{(\bu^*\mathbf{1})^2}{n^2}z\mathbf{1}^*\tilde{Q}_n\mathbf{1}
%\end{aligned}
%\end{equation*}
%\begin{equation*}\begin{aligned}
%     & \tilde{r}^2\bu^*\Phi \bA_n^* (\bI+\tilde{r}{\bR}_n)^{-1}\bA_n\Phi \bu-\tilde{r}^2\bu^*\bA_n^* (\bI+\tilde{r}{\bR}_n)^{-1}\bA_n \bu\\
%     &=-2\tilde{r}^2\bu^*\mathbf{1}\mathbf{1}^*\bA_n^* (\bI+\tilde{r}{\bR}_n)^{-1}\bA_n \bu + \tilde{r}^2(\bu^*\mathbf{1})^2 \mathbf{1}^*\bA_n^* (\bI+\tilde{r}{\bR}_n)^{-1}\bA_n\mathbf{1}
%\end{aligned}\end{equation*}

%With the tools prepared in 
Relying on Proposition \ref{cenprop},  we have the following conclusion for the spiked eigenvalues and the corresponding eigenvectors of $\bold{\mathcal{S}}_n$ and $\tilde{\bold{\mathcal{S}}}_n$.
\begin{theorem}\label{centhm3}
Assume that the conditions of Proposition \ref{cenprop} are satisfied with $\bR_n$ in Assumption \ref{ass4new} replaced by $\bar{\bR}_n$. By replacing $\bS_n$, $\tilde\bS_n$, $\bR_n$ and their latent symbols (e.g., $\gamma_k$) with the counterparts of $\bold{\mathcal{S}}_n$, $\tilde{\bold{\mathcal{S}}}_n$ and $\bar{\bR}_n$,
	 the conclusions in Theorems \ref{thmeve1} and \ref{thm2} still hold.
\end{theorem}

\section{Applications }\label{app}
In this section, based on the results in Section \ref{secModelResults}, we aim to develop some potential applications. 
%\subsection{Clustering}
Spectral clustering has been used in practice frequently in data science and the theoretical underpinning of such a method has received extensive interest in recent years; see e.g., \citep{couillet2016kernel, zhou2019analysis,loffler2021optimality}. This section is to have a deep insight into the spectral clustering based on the Model \eqref{mo2}. Moreover, we also propose a new criterion to estimate the number of clusters. 
   Recalling \eqref{mo2}, for any $i\in \VV_s$, there is ${\mathbb{E}}\bx_i=\bmu_s/\sqrt{n}$, where $s=1,\ldots, K$. 
Let $\bN=[\bmu_1,\ldots,\bmu_{K}]/\sqrt{n}\in\mR^{p\times K}$, $\bH=[\bh_1,\ldots,\bh_{K}]\in \mR^{n\times K}$, $\bh_s=(\bh_{s}(1),\ldots,\bh_{s}(n))^*\in \mR^n$, where $\bh_{s}(i)=1$ if $i\in\VV_s$ and  $\bh_{s}(i)=0$ otherwise.
In a matrix form, write \begin{equation*}\label{datamat2}
 \bX_n=[\bx_1,\ldots,\bx_n]=\bN\bH^*+\bSigma_n^{1/2}\bW_n.
 \end{equation*}
  Notice that \begin{equation*}
{\mathbb{E}} (\tilde{\bS}_n) = \bH \bN^* \bN \bH^*+ \frac{\trace \bSigma_n}{n} \bI_n. 
\end{equation*}
The block structure of ${\mathbb{E}}( \tilde{\bS}_n)$ (except the diagonal positions) is similar to that of stochastic block models (SBM). This motivates one to use spectral clustering for high dimensional data with different means across groups.

  To do the clustering, it is of interest to estimate the number of clusters, i.e., estimation of $K$. There exist plenty of approaches to estimate the number of clusters. %To name a few, 
  For instance, \cite{Robert1953} proposed the Elbow method that aims to minimize the within-group sum of squares (WSS); Silhouette index \citep{Rousseeuw1987} is a measure of how similar an object is to its own cluster compared to other clusters, which takes values in $[-1,1]$; and \cite{Robert2001} proposed a gap statistic to estimate the number of clusters.  These methods either lack theoretical guarantees or have some restrictions in computation or settings. Hence, here we propose a theoretical guarantee and easily implemented approach to estimate the number of clusters.  Notice that under Model \eqref{mo2}, the number of the spiked eigenvalues of $\bS_n$ or $\tilde\bS_n$ is the same as the number of clusters if the means in terms of the different clusters are not linearly correlated and the noise is relatively weak. In the remainder of this section, we consider the scenario that the number of population spiked eigenvalues satisfying Assumption \ref{ass4new} equals the number of clusters, i.e., $K=\sum_{i=1}^L m_i$. The estimation of the number of spikes in different models has been discussed in multiple literatures, and mostly are based on the setting of $\bSigma_n=\bI$; see e.g., \cite{Bai2018}.

%{\color{blue}We assume that only the signal parts contribute to the spikes }
  Motivated by the work of \cite{Bai2018} and Theorem \ref{thm2}, we propose two criteria to estimate the number of clusters. Without loss of generality, we assume $0<c<1$. Let 
\begin{equation}\label{crrr}
\begin{split}
& \text{EDA}_{k}=-n(\hat{\lambda}_{1}-\hat{\lambda}_{k+1})+n(p-k-1)\log\widetilde{\theta}_{p,k}+2pk,\\
& \text{EDB}_{k}=-n\log(p)\cdot(\hat{\lambda}_{1}-\hat{\lambda}_{k+1})+n(p-k-1)\log\widetilde{\theta}_{p,k}+(\log n)pk,
\end{split}
\end{equation}
where  $\widetilde{\theta}_{p,k}=\frac{1}{p-k-1}\sum_{i=k+1}^{p-1}\theta_{i}^2$, and $\theta_{k}=\exp\{\hat{\lambda}_{n,k}-\hat{\lambda}_{n,k+1}\},k=1,2,\ldots,p-1$.
\begin{rmk}
	The first two main terms aim to capture the difference between eigenvalues, and the third term is the penalty term for the number of unknown parameters in the model. The values of EDA and EDB are expected to reach a minimum when $k=K$. From (\ref{crrr}), it can be seen that, as $k$ increases, the first and second terms decrease while the third term increases. For more discussion about \eqref{crrr} and the case of $c>1$, one may refer to the supplementary material.
\end{rmk}
 %Then,  
 We estimate the number of clusters by
\begin{eqnarray*}
\hat{K}_{\text{EDA}}&=&\arg\min\limits_{k=1,\ldots,w}\frac{1}{n}\text{EDA}_{k},\label{equ3.2}\\
\hat{K}_{\text{EDB}}&=&\arg\min\limits_{k=1,\ldots,w}\frac{1}{n}\text{EDB}_{k},\label{equ3.2.5}
\end{eqnarray*}
where $w$ is the prespecified number of clusters satisfying  $w=o(p)$. We denote the eigenvalues of $\bR_n$ by $\gamma_1\ge \gamma_2 \ge \gamma_3,\cdots, \ge \gamma_p$ in the following but we still allow the existence of multiple eigenvalues, as in Assumption \ref{ass3new}. Note that under conditions of Theorem \ref{thm2}, it follows that for $k=1,2,\ldots,K-1$, with probability one,
\begin{equation}\label{equ3.3}
\theta_{k}\to \exp\{\varphi(\gamma_{k})-\varphi(\gamma_{k+1})\},~~\theta_{K}\to \exp\{\varphi(\gamma_{K})-\mu_{K+1}\}, 
\end{equation}
where function $\varphi$ and the limit of the largest nonspiked eignvalue $\mu_{K+1}$ are defined in  (\ref{equ1.3})  and \eqref{equ2.5new}, respectively. For simplicity, denote the limit of $\theta_{k}$ by $\zeta_{k}$ for $k=1,\ldots,K$. Define two sequences $\{a_{s}\}_{s=2}^{K}$ and $\{b_{s}\}_{s=2}^{K}$ as follows
\begin{equation}\label{equ3.5.5}
\begin{split}
	&a_{s}=\zeta_{s}^{2}+\log\zeta_{s}-2c-1+a_{s+1} \text{ and } a_{K+1}=0 \text{ for } s=2,\ldots,K,\\
	&b_{s}=\zeta_{s}^{2}+{\log p}\log\zeta_{s}-c\log n-1+b_{s+1}\text{ and } b_{K+1}=0 \text{ for } s=2,\ldots,K.
\end{split}
\end{equation}
We propose two gap conditions for EDA and EDB, respectively, i.e.,
\begin{eqnarray}
\min\limits_{s=2,\ldots,K}a_{s}&>&0,\label{equ3.6}\\
\min\limits_{s=2,\ldots,K}b_{s}&>&0.\label{equ3.6.5}
\end{eqnarray}

\begin{rmk}
The gap condition in \cite{Bai2018}	was proposed for the population covariance matrix with distant spikes larger than one and other eigenvalues equal to one. While the model studied in this paper imposes no restriction to the non-spiked eigenvalues, the gap conditions in (\ref{equ3.6}) and (\ref{equ3.6.5}) are more easily satisfied and have a wider range of applications.
\end{rmk}
Note that Theorem \ref{thm2} and \eqref{equ3.3} are obtained when the leading eigenvalues are bounded. Here we also investigate the cases of the leading eigenvalues tending to infinity. We need the following result, and its proof is in Section S.2.3 of the supplementary material.
\begin{lemma}\label{lem2.4}
In the same setup of Theorem \ref{thm2}, instead of assuming $\gamma_{1}$ being bounded, suppose that $\gamma_{K}\rightarrow\infty$, as $n\rightarrow\infty$. Then, for any $k=1,\ldots,K$, we have
\begin{equation*}
\lim_{n\rightarrow\infty}\hat{\lambda}_k/\gamma_{k}=1~~~~ \text{a.s.}.
\end{equation*}
\end{lemma}

We say that $\hat{K}$, the estimator of $K$, is strongly consistent if $\hat{K}=K$ almost surely.
Based on Theorem \ref{thm2} and Lemma \ref{lem2.4}, we derive the consistency of $\hat{K}_{\text{EDA}}$ as follows. The proof is in Section 2.3 of the supplementary material.
\begin{theorem}\label{thm3.1}
Under conditions of Theorem \ref{thm2}, %We have that
if the gap condition (\ref{equ3.6}) does not hold, then $\hat{K}_{\text{EDA}}$ is not consistent; if the gap condition holds, then $\hat{K}_{\text{EDA}}$ is strongly consistent.
In particular, if $\gamma_{K}$ tends to infinity, then $\hat{K}_{\text{EDA}}$ is strongly consistent.
\end{theorem}

In \cite{Bai2018}, BIC is consistent when $\gamma_{K}\rightarrow\infty$ at a rate faster than $\log n$, which makes BIC less capable of detecting signals. This is because BIC has a more strict penalty coefficient $\log n$ compared to the penalty coefficient 2 in AIC. For the EDB construction of selecting the number of clusters, we add the coefficient $\log p$ to the first term so that the spikes do not need to be very large and only the corresponding gap condition for EDB is required. By the analogous proof strategy of Theorem \ref{thm3.1}, see Section S.2.3 in the supplementary file, we obtain the consistency of EDB as follows.
\begin{theorem}\label{thm3.2}
Under the same setting of Theorem \ref{thm3.1},
%we also have that
if the gap condition (\ref{equ3.6.5}) does not hold, then $\hat{K}_{\text{EDB}}$ is not consistent; if (\ref{equ3.6.5}) holds, then $\hat{K}_{\text{EDB}}$ is strongly consistent.
Moreover,  if $\gamma_{K}$ tends to infinity, then $\hat{K}_{\text{EDB}}$ is strongly consistent.
\end{theorem}

%From Theorems \ref{thm3.1} and \ref{thm3.2}, we can have an estimator of the number of clusters, and hence we can conduct spectral clustering in the next step.
Once the estimator of the number of clusters is available, we can conduct spectral clustering. Specifically, 
   let the eigenvectors corresponding to the first $\hat K$ eigenvalues of $\tilde\bS_{n}$ be $\widehat{\mathbf U}=[\hat\bu_1,\ldots,\hat\bu_{\hat K}]\in\mR^{n\times \hat K}$. % and 
   We then apply the following K-means optimization to the $\widehat\bU$, i.e., \begin{eqnarray}\label{cr3}
	\bf\tilde{U}=\argmin_{\bU\in\M_{n,\hat K}}\|\bU-\widehat \bU\|_F^2,
\end{eqnarray}
where $\M_{n,K}=\{\bU\in\mR^{n\times K}:\bU\text{ has at most $K$ distinct rows}\}$. Then, we return $\hat\VV_1,\ldots,\hat\VV_{\hat K}$ as the indices for each cluster. From \eqref{cr3}, we see that the spectral clustering is conducted from the obtained $\widehat\bU$, and hence we look into the properties of $\widehat\bU$. In the following, we consider the case that the first $K$ eigenvalues of $\bR_n = \bN\bH^*\bH\bN^*+\bSigma$, denoted as $\gamma_1,\cdots,\gamma_K$, are all single eigenvalues satisfying Assumption \ref{ass4new}, and the associated eigenvectors are denoted as $\bxi_1,\cdots,\bxi_K$. The proof of the following result is in Section S.2.3 of the supplementary material.

% an extend the above result to the eigenvectors corresponding to the largest $r$ eigenvalues, where $r\leq K$.
\begin{corollary}\label{cor1}
Under the conditions of Theorem \ref{thmeve1},  % we have that,
 in the set of all deterministic unit vectors $\bu_n$,  $\tilde{\bu} = \bA_n^*\bxi_k / \|\bA_n^*\bxi_k\|$ %{\color{red} 
maximizes the non-random term   $\gamma_k^{-1}\eta_k\bu_n^* \bA_n^* \bxi_k \bxi_k^* \bA_n \bu_n$ %, how do you minimize it? or are you referring to the left hand of 16 which is random rather than non-random?} 
 in \eqref{l2normvandhatu}, and
% \begin{equation}\label{eivecconsistent}\left\|\frac{\bA_n^*\xi_k}{\|\bA_n^*\xi_k\|}-\hat\bu_k\right\|^2=2\left[1-\theta_k\frac{\xi_k^*\bA_n\bA_n^*\xi_k}{n\gamma_k}\right]+O_P(\frac{1}{\sqrt{n}}).\end{equation}
\begin{equation}\label{eivecconsistent} \left\|(\hat{\bu}_k^* \tilde{\bu})\tilde{\bu} - \hat\bu_k\right\|^2=1-\eta_k \left (1- \frac{\bxi_k^* \bSigma_n \bxi_k}{\gamma_k}\right)+O_P\left(\frac{1}{\sqrt{n}}\right).\end{equation}
%where $\eta_k$ is defined in \eqref{etak}.
Moreover, let $\widehat\bU_r$ be the eigenvectors corresponding to the largest $r$ eigenvalues of $\tilde\bS_n$, where $r\leq K$.
For any deterministic $\bV_r$ that contains $r$ column vectors of unit length, we have
%{\color{red} Tr?}
\begin{equation}\label{FnormhatU}
 \inf_{\boldsymbol\Lambda \in \mathbb{R}^{r\times r}} \| \bV_r \boldsymbol\Lambda - \widehat{\bU}_r \|_F ^2= r -  \trace \left(\bV_r^* \widehat{\bU}_r \widehat{\bU}_r^* \bV_r\right) = r- \trace \left(\bV_r^* \bA^* \bP_R \bA \bV_r\right)+ O_P \left(\frac{1}{\sqrt{n}}\right),
\end{equation}
where $$\bP_R = \sum_{k=1}^r \frac{\eta_k}{\gamma_k}\bxi_k \bxi_k^*.$$
%Specifically, denote by $\bV_r^*$ the matrix formed by  the first $r$ eigenvectors of $\bA^* \bP_R \bA $. Then $\bV_r^*$ minimize the right side of \eqref{FnormhatU}.
%\begin{equation}\inf_{s \in \{1,-1\}} \| \bV \bV^* -  \hat{\bU}_r \hat{\bU}_r^* \|_F ^2= 2 r - 2...+ O_P (\frac{1}{\sqrt{n}}\end{equation}
\end{corollary}
\begin{rmk}
	From Corollary \ref{cor1}, we see that if   $\gamma_k$ tends to infinity,  and $\gamma_{i-1}/\gamma_{i}>1+\delta$  for $1\leq i\leq K$ with $\delta$ being a positive constant independent of $n$, we have $\eta_k\rightarrow1$ %and $\frac{\xi_k^*\bA_n\bA_n^*\xi_k}{n\gamma_k}\rightarrow 0$ as $n\rightarrow\infty$.
thus the right side of  \eqref{eivecconsistent} converges to zero in probability. Consequently, $\hat\bu_k$ is an asymptotic consistent estimator of $\frac{\bA_n^*\bxi_k}{\|\bA_n^*\bxi_k\|}$. Note that $\bA_n=\bN\bH^*$, which has $K$ distinct columns and represents $K$ different means. Hence, under mild conditions, there are $K$ different rows in $\widehat\bU$, and one can use it to find the corresponding clusters. When $\gamma_k$ is bounded, $\hat{\bu}_k$ is not a consistent estimator for the block-wise constant vector $\bA_n^*\bxi_k/\|\bA_n^*\bxi_k\|\in\mR^n$. However, in this case, following the proof of Theorem 2.2 in \cite{jin2015fast}, an elementary misclustering error rate by spectral clustering can be also obtained, which is a new observation based on the proposed results. 
\end{rmk}

\section{Simulation}\label{simu}
%This section provides some simulation studies with respect to the applications mentioned in Section \ref{app}. Additional simulation can be also found in the supplementary material.%Specifically,  Subsection \eqref{sim:clu} 
%{\subsection{Clustering}}\label{sim:clu}

In this section, we first evaluate the performance of the proposed criteria in the estimation of the number of clusters discussed in Section \ref{app}. Denote the sets of under-estimated, exactly estimated and over-estimated models by $\mathcal {F}_{-}, \mathcal {F}_{*}$ and $\mathcal {F}_{+}$, respectively, i.e.,
\begin{equation*}
\mathcal {F}_{-}=\{1,\ldots,K-1\},~~\mathcal {F}_{*}=\{K\},~~\mathcal {F}_{+}=\{K+1,\ldots,w\}.
\end{equation*}
The selection percentages corresponding to $\mathcal {F}_{-}, \mathcal {F}_{*}$ and $\mathcal {F}_{+}$ are computed by $1000$ repetitions. Suppose that the entries of $\bW_{n}$ are i.i.d. with the following distributions:\begin{itemize}
\item 	Standard normal distribution: $w_{i,j}\sim \mathcal {N}(0,1)$.
\item Standardized $t$ distribution with 8 degrees of freedom: $w_{i,j}\sim t_{8}/\sqrt{\text{Var}(t_{8})}$.
\item Standardized Bernoulli distribution with probability $1/2$: $w_{i,j}\sim(\text{Bernoulli}(1,1/2)-1/2)/(1/2)$.
\item  Standardized chi-square distribution with 3 degrees of freedom: $w_{i,j}\sim (\chi^{2}(3)-3)/\sqrt{\text{Var}(\chi^{2}(3))}\\=(\chi^{2}(3)-3)/\sqrt{6}$
\end{itemize}
For comparison,  three different methods are also considered: Average Silhouette Index (\cite{Rousseeuw1987}), Gap Statistic (\cite{Robert2001}) and BIC with degrees of freedom (\cite{David2020}), denoted by ASI, GS and BICdf, respectively. This section considers the situations with $0<c<1$, and some more cases with  $0<c<1$ and   $c>1$ are demonstrated in the supplementary material. 
Here we set $c=1/3,3/4$ and the largest number of possible clusters $w=\lfloor 6\cdot n^{0.1}\rfloor$. Different means in terms of different clusters and the covariance matrices are set as follows
:

\textbf{Case 1.} Let $\bmu_{1}=(3,0,0,0,\ldots,0)^{*}\in\mathbb{R}^{p}$, $\bmu_{2}=(0,3,0,0,\ldots,0)^{*}\in\mathbb{R}^{p}$, $\bmu_{3}=(0,0,3,0,\ldots,0)^{*}\in\mathbb{R}^{p}$,
$\bSigma_n=\bI$, where $\bI$ is the identity matrix of size $p$. Then,
$$\bA_{n}=(\underbrace{\bmu_{1},\ldots,\bmu_{1}}_{n_{1}},\underbrace{\bmu_{2},\ldots,\bmu_{2}}_{n_{2}},\underbrace{\bmu_{3},\ldots,\bmu_{3}}_{n_{3}}),$$
where $n_{1}=n_{2}=0.3n,~n_{3}=0.4n$. Therefore, the true number of clusters is $K=3$.

Note that the spikes in the above case are bounded. We also consider a case of %infinite
spikes with $\gamma_{K}\rightarrow\infty$ at a rate faster than $\log n$ and $\gamma_{1}=O(p)$.

\textbf{Case 2.} Let $\bmu_{1}=(2a,a,-a,a,1,\ldots,1)^{*}\in\mathbb{R}^{p}$, $\bmu_{2}=(a,a,2a,-3a,1,\ldots,1)^{*}\in\mathbb{R}^{p}$, \\$\bmu_{3}=(a,-2a,-a,a,1,\ldots,1)^{*}\in\mathbb{R}^{p}$,
$\bmu_{4}=(-2a,a,a,a,1,\ldots,1)^{*}\in\mathbb{R}^{p}$, and the sample size of cluster corresponding to each center be $n_{1}=n_{3}=0.3n,~n_{2}=n_{4}=0.2n$, such that the true number of clusters $K=4$. Suppose $\bSigma_n=(\sigma_{i,j})_{p\times p}$, where $a=\sqrt{p/10}$, $\sigma_{i,j}=0.2^{|i-j|}$. 

Tables \ref{table2} to \ref{table4} report the percentages of under-estimated, exactly estimated and over-estimated under 1000 replications. From the reported results, we see  the criteria based on EDA and EDB work better and better as $n,p$ become larger. When $c=1/3$, the probabilities of the under-estimated number of clusters are equal to $0$ and increase when $c$ is getting closer to $1$. From (\ref{equ3.5.5}), it is shown that the larger $c$ is, the harder the gap conditions are to be satisfied.
%Thus, some observations may not meet the those gap conditions which hold theoretically in Case 1 such that the under-estimated probabilities are not identically equal to zero. 
EDB generally outperforms EDA except the case of $c=3/4$, when $p,n$ are large. It can be seen that when $c=3/4$, as $n$ increases, the probability of $\mathcal {F}_{-}$ estimated by EDB becomes larger, and is uniformly greater than that by EDA. This is due to the fact that the coefficient in the penalty term of EDB criterion is $\log n$ which is different from the coefficient $2$ in EDA, so that the gap condition of EDB is more stronger than of EDA, that is, (\ref{equ3.6.5}) is more difficult to be satisfied than (\ref{equ3.6}). The criteria based on EDA and EDB show the highest accuracy under Bernoulli distribution, followed by normal, $t_{8}$ and $\chi^{2}(3)$ with relatively heavy right tail which may be destructive to the results.

%\addtocounter{table}{-3}
%%%%%%% Table 2 %%%%%%%%%%%
\begin{table}[H]\centering
\caption{Selection percentages of EDA, EDB, ASI, GS and BICdf in Case 1. Entries in the $\mathcal {F}_{*}$ rows indicate that EDA and EDB exhibit higher accuracy in estimating the number of clusters compared to other criteria.}\label{table2}
\begin{tabular}{c|c c   c c c c c | c c c c c}\hline
    &     &  & EDA & EDB & ASI & GS & BICdf & EDA & EDB & ASI & GS & BICdf\\\hline
$c$ & $n$ &  & \multicolumn{5}{c|}{$\mathcal {N}(0,1)$} & \multicolumn{5}{c}{$t_{8}$}\\\hline
\multirow{6}[6]*{$\frac{1}{3}$} & \multirow{3}[3]*{$180$} & $\mathcal {F}_{-}$  &
                                $0$ & $0$ & $6.1$ & $64.8$ & $34.7$ & $0$ & $0$ & $5.3$ & $80.2$ & $34.8$\\
                                &                         & $\mathcal {F}_{*}$  &
                                $80.4$ & $95.3$ & $93.9$ & $35.2$ & $59.8$ & $78.4$ & $91.1$ & $94.2$ & $19.8$ & $59.3$\\
                                &                         & $\mathcal {F}_{+}$  &
                                $19.6$ & $4.7$ & $0$ & $0$ & $5.5$ & $21.6$ & $8.9$ & $0.5$ & $0$ & $5.9$\\\cline{2-13}
                                & \multirow{3}[3]*{$450$} & $\mathcal {F}_{-}$  &
                                $0$ & $0$ & $3.4$ & $41.9$ & $98.5$ & $0$ & $0$ & $4.3$ & $72.5$ & $98.3$\\
                                &                         & $\mathcal {F}_{*}$  &
                                $99.4$ & $100$ & $96.6$ & $58.1$ & $1.5$ & $98.9$ & $99.8$ & $95.6$ & $27.5$ & $1.7$\\
                                &                         & $\mathcal {F}_{+}$  &
                                $0.6$ & $0$ & $0$ & $0$ & $0$ & $1.1$ & $0.2$ & $0.1$ & $0$ & $0$\\\cline{2-13}
\multirow{6}[6]*{$\frac{1}{2}$} & \multirow{3}[3]*{$120$} & $\mathcal {F}_{-}$  &
                                $0$ & $0$ & $11.2$ & $97.1$ & $98$ & $0$ & $0$ & $12.1$ & $98.8$ & $98.5$\\
                                &                         & $\mathcal {F}_{*}$  &
                                $83.4$ & $93.7$ & $88.1$ & $2.9$ & $2$ & $82$ & $91.3$ & $86.1$ & $1.2$ & $1.5$\\
                                &                         & $\mathcal {F}_{+}$  &
                                $16.6$ & $6.3$ & $0.7$ & $0$ & $0$ & $18$ & $8.7$ & $1.8$ & $0$ & $0$\\\cline{2-13}
                                & \multirow{3}[3]*{$300$} & $\mathcal {F}_{-}$  &
                                $0$ & $0$ & $5.7$ & $97.1$ & $100$ & $0$ & $0$ & $7.9$ & $98.6$ & $100$\\
                                &                         & $\mathcal {F}_{*}$  &
                                $99.3$ & $100$ & $94.3$ & $2.9$ & $0$ & $99$ & $99.8$ & $91.8$ & $1.4$ & $0$\\
                                &                         & $\mathcal {F}_{+}$  &
                                $0.7$ & $0$ & $0$ & $0$ & $0$ & $1$ & $0.2$ & $0.3$ & $0$ & $0$\\\hline
\multirow{6}[6]*{$\frac{3}{4}$} & \multirow{3}[3]*{$80$}  & $\mathcal {F}_{-}$  &
                                $1.1$ & $2.8$ & $22.1$ & $100$ & $100$ & $1.3$ & $2.7$ & $21.5$ & $100$ & $100$\\
                                &                         & $\mathcal {F}_{*}$  &
                                $86.9$ & $91.1$ & $76.1$ & $0$ & $0$ & $82.4$ & $86.6$ & $73.6$ & $0$ & $0$\\
                                &                         & $\mathcal {F}_{+}$  &
                                $12$ & $6.1$ & $1.8$ & $0$ & $0$ & $16.3$ & $10.7$ & $4.9$ & $0$ & $0$\\\cline{2-13}
                                & \multirow{3}[3]*{$200$} & $\mathcal {F}_{-}$  &
                                $0.1$ & $0.1$ & $13.2$ & $100$ & $100$ & $0$ & $0.4$ & $13.3$ & $100$ & $100$\\
                                &                         & $\mathcal {F}_{*}$  &
                                $99.1$ & $99.8$ & $86.4$ & $0$ & $0$ & $99.4$ & $99.6$ & $85.5$ & $0$ & $0$\\
                                &                         & $\mathcal {F}_{+}$  &
                                $0.8$ & $0.1$ & $0.4$ & $0$ & $0$ & $0.6$ & $0$ & $1.2$ & $0$ & $0$\\\hline

$c$ & $n$ &  & \multicolumn{5}{c|}{Bernoulli} & \multicolumn{5}{c}{$\chi^{2}(3)$}\\\hline
\multirow{6}[6]*{$\frac{1}{3}$} & \multirow{3}[3]*{$180$} & $\mathcal {F}_{-}$  &
$0$ & $0$ & $2.8$ & $23.8$ & $33.8$ & $0$ & $0$ & $7.4$ & $88.4$ & $33.8$\\
                                &                         & $\mathcal {F}_{*}$  &
                                $84.6$ & $94.3$ & $97.1$ & $76.2$ & $62.4$ & $71.9$ & $87.7$ & $86.7$ & $11.6$ & $60.8$\\
                                &                         & $\mathcal {F}_{+}$  &
                                $15.4$ & $5.7$ & $0.1$ & $0$ & $3.8$ & $28.1$ & $12.3$ & $5.9$ & $0$ & $5.4$\\\cline{2-13}
                                & \multirow{3}[3]*{$450$} & $\mathcal {F}_{-}$  & $0$ & $0$ & $1.1$ & $8.5$ & $98.5$ & $0$ & $0$ & $4.1$ & $85.6$ & $98.9$\\
                                &                         & $\mathcal {F}_{*}$  &
                                $99.5$ & $100$ & $98.9$ & $91.5$ & $1.5$ & $98.4$ & $99.9$ & $95.2$ & $14.4$ & $1.1$\\
                                &                         & $\mathcal {F}_{+}$  &
                                $0.5$ & $0$ & $0$ & $0$ & $0$ & $1.6$ & $0.1$ & $0.7$ & $0$ & $0$\\\hline
\multirow{6}[6]*{$\frac{1}{2}$} & \multirow{3}[3]*{$120$} & $\mathcal {F}_{-}$  &
$0$ & $0$ & $11$ & $91$ & $99.6$ & $0.1$ & $0$ & $13.3$ & $99.6$ & $98.3$\\
                                &                         & $\mathcal {F}_{*}$  &
                                $88.3$ & $95.4$ & $89$ & $9$ & $0.4$ & $74.5$ & $84.4$ & $78.7$ & $0.4$ & $1.1$\\
                                &                         & $\mathcal {F}_{+}$  &
                                $11.7$ & $4.6$ & $0$ & $0$ & $0$ & $25.4$ & $15.6$ & $8$ & $0$ & $0$\\\cline{2-13}
                                & \multirow{3}[3]*{$300$} & $\mathcal {F}_{-}$  & $0$ & $0$ & $3$ & $85.4$ & $100$ & $0$ & $0$ & $9.2$ & $99.8$ & $98.3$\\
                                &                         & $\mathcal {F}_{*}$  &
                                $99.8$ & $100$ & $97$ & $14.6$ & $0$ & $99.5$ & $100$ & $89.8$ & $0.2$ & $1.7$\\
                                &                         & $\mathcal {F}_{+}$  & $0.2$ & $0$ & $0$ & $0$ & $0$ & $0.5$ & $0$ & $1$ & $0$ & $0$\\\hline
\multirow{6}[6]*{$\frac{3}{4}$} & \multirow{3}[3]*{$80$} & $\mathcal {F}_{-}$   &
$0.4$ & $1$ & $23.5$ & $100$ & $99.9$ & $5.4$ & $7$ & $23.4$ & $100$ & $100$\\
                                &                         & $\mathcal {F}_{*}$  &
                                $90.2$ & $94.4$ & $76.4$ & $0$ & $0.1$ & $75$ & $79.6$ & $58.1$ & $0$ & $0$\\
                                &                         & $\mathcal {F}_{+}$  &
                                $9.4$ & $4.6$ & $0.1$ & $0$ & $0$ & $19.6$ & $13.4$ & $18.5$ & $0$ & $0$\\\cline{2-13}
                                & \multirow{3}[3]*{$200$} & $\mathcal {F}_{-}$  &
                                $0$ & $0.2$ & $7.7$ & $100$ & $100$ & $0.3$ & $1.1$ & $20.1$ & $100$ & $100$\\
                                &                         & $\mathcal {F}_{*}$  &
                                $99.7$ & $99.8$ & $92.2$ & $0$ & $0$ & $99.1$ & $98.9$ & $76.4$ & $0$ & $0$\\
                                &                         & $\mathcal {F}_{+}$  & $0.3$ & $0$ & $0.1$ & $0$ & $0$ & $0.6$ & $0$ & $3.5$ & $0$ & $0$\\\hline
\end{tabular}
\end{table}

%%%%%%% Table 4 %%%%%%%%%%%
\begin{table}[H]
\centering
\caption{Selection percentages of EDA, EDB, ASI, GS and BICdf in Case 2. Entries in the $\mathcal {F}_{*}$ rows indicate that EDA and EDB exhibit higher accuracy in estimating the number of clusters compared to other criteria.}\label{table4}
\begin{tabular}{c|c c   c c c c c | c c c c c}\hline    &     &  & EDA & EDB & ASI & GS & BICdf & EDA & EDB & ASI & GS & BICdf\\\hline$c$ & $n$ &  & \multicolumn{5}{c|}{$\mathcal {N}(0,1)$} & \multicolumn{5}{c}{$t_{8}$}\\\hline\multirow{6}[6]*{$\frac{1}{3}$} & \multirow{3}[3]*{$180$} & $\mathcal {F}_{-}$  & $0$ & $0$ & $70.4$ & $25.1$ & $0.2$ & $0$ & $0$ & $70$ & $28.1$ & $0$\\                                &                         & $\mathcal {F}_{*}$  &                                $60.4$ & $82.4$ & $28.9$ & $65.5$ & $67.1$ & $57.1$ & $78.3$ & $28.7$ & $62.1$ & $66.5$\\                                &                         & $\mathcal {F}_{+}$  &                                $39.6$ & $17.6$ & $0.7$ & $9.4$ & $32.7$ & $42.9$ & $21.7$ & $1.3$ & $9.8$ & $33.5$\\\cline{2-13}                                & \multirow{3}[3]*{$450$} & $\mathcal {F}_{-}$  &                                $0$ & $0$ & $67.2$ & $28.4$ & $0$ & $0$ & $0$ & $69.2$ & $31.4$ & $0.1$\\                                &                         & $\mathcal {F}_{*}$  &                                $94.8$ & $99.6$ & $29.8$ & $61.2$ & $62.7$ & $93.5$ & $99$ & $27.8$ & $58$ & $63$\\                                &                         & $\mathcal {F}_{+}$  &                                $5.2$ & $0.4$ & $3$ & $10.4$ & $37.3$ & $6.5$ & $1$ & $3$ & $10.6$ & $36.9$\\\cline{2-13}\multirow{6}[6]*{$\frac{1}{2}$} & \multirow{3}[3]*{$120$} & $\mathcal {F}_{-}$  &$0$ & $0$ & $71.2$ & $28.3$ & $1.8$ & $0$ & $0$ & $69.2$ & $28.8$ & $2.5$\\                                &                         & $\mathcal {F}_{*}$  &                                $70.4$ & $83.1$ & $27.3$ & $63.4$ & $70.8$ & $62.9$ & $83.4$ & $29.8$ & $63.4$ & $66.2$\\                                &                         & $\mathcal {F}_{+}$  &                                $29.6$ & $16.9$ & $1.5$ & $8.3$ & $27.4$ & $37.1$ & $16.6$ & $1$ & $7.8$ & $31.3$\\\cline{2-13}                                & \multirow{3}[3]*{$300$} & $\mathcal {F}_{-}$  & $0$ & $0$ & $66.7$ & $28.3$ & $0$ & $0$ & $0$ & $65.5$ & $26.7$ & $0$\\                                &                         & $\mathcal {F}_{*}$  &                                $96.7$ & $99.7$ & $30.5$ & $59.5$ & $59.6$ & $96.4$ & $99.8$ & $31.1$ & $62.2$ & $61.7$\\                                &                         & $\mathcal {F}_{+}$  &                                $3.3$ & $0.3$ & $2.8$ & $12.2$ & $40.4$ & $3.6$ & $0.2$ & $3.4$ & $11.1$ & $38.3$\\\hline\multirow{6}[6]*{$\frac{3}{4}$} & \multirow{3}[3]*{$80$}  & $\mathcal {F}_{-}$  &$0$ & $0$ & $67.1$ & $30.6$ & $13.9$ & $0$ & $0$ & $69.5$ & $36.2$ & $13$\\                                &                         & $\mathcal {F}_{*}$  &                                $75$ & $84.4$ & $31.4$ & $61.2$ & $68.4$ & $68.8$ & $78.4$ & $29.5$ & $56.8$ & $68.9$\\                                &                         & $\mathcal {F}_{+}$  &                                $25$ & $15.6$ & $1.5$ & $8.2$ & $17.7$ & $31.2$ & $21.6$ & $1$ & $7$ & $18.1$\\\cline{2-13}                                & \multirow{3}[3]*{$200$} & $\mathcal {F}_{-}$  &                                $0$ & $0$ & $68.7$ & $31.1$ & $0.2$ & $0$ & $0$ & $67.7$ & $30.1$ & $0$\\                                &                         & $\mathcal {F}_{*}$  &                                $97$ & $99.7$ & $28.1$ & $58.8$ & $63$ & $96$ & $99.3$ & $29.5$ & $59.8$ & $60.1$\\                                &                         & $\mathcal {F}_{+}$  &                                $3$ & $0.3$ & $3.2$ & $10.1$ & $36.8$ & $4$ & $0.7$ & $2.8$ & $10.1$ & $39.9$\\\hline$c$ & $n$ &  & \multicolumn{5}{c|}{Bernoulli} & \multicolumn{5}{c}{$\chi^{2}(3)$}\\\hline\multirow{6}[6]*{$\frac{1}{3}$} & \multirow{3}[3]*{$180$} & $\mathcal {F}_{-}$  &$0$ & $0$ & $67.5$ & $24.3$ & $0$ & $0$ & $0$ & $66.4$ & $24.2$ & $0.1$\\                                &                         & $\mathcal {F}_{*}$  &                                $60.7$ & $83.6$ & $32.2$ & $66.9$ & $66.3$ & $49.5$ & $75.9$ & $31.7$ & $67.5$ & $66.6$\\                                &                         & $\mathcal {F}_{+}$  &                                $39.3$ & $16.4$ & $0.3$ & $8.8$ & $33.7$ & $50.5$ & $24.1$ & $1.9$ & $8.1$ & $33.3$\\\cline{2-13}                                & \multirow{3}[3]*{$450$} & $\mathcal {F}_{-}$  & $0$ & $0$ & $70.7$ & $29.3$ & $0$ & $0$ & $0$ & $67.1$ & $30.1$ & $0$\\                                &                         & $\mathcal {F}_{*}$  &                                $95.1$ & $99.4$ & $26.6$ & $58.3$ & $62.4$ & $91.9$ & $99.4$ & $29.3$ & $56.6$ & $60.9$\\                                &                         & $\mathcal {F}_{+}$  &                                $4.9$ & $0.6$ & $2.7$ & $12.4$ & $37.6$ & $8.1$ & $0.6$ & $3.6$ & $13.3$ & $39.1$\\\hline\multirow{6}[6]*{$\frac{1}{2}$} & \multirow{3}[3]*{$120$} & $\mathcal {F}_{-}$  &$0$ & $0$ & $72.6$ & $29.2$ & $1.9$ & $0$ & $0$ & $68.1$ & $30.2$ & $1.9$\\                                &                         & $\mathcal {F}_{*}$  &                                $72.6$ & $85$ & $25.6$ & $63.6$ & $66.7$ & $62.2$ & $74.7$ & $29.3$ & $61.5$ & $66.6$\\                                &                         & $\mathcal {F}_{+}$  &                                $27.4$ & $15$ & $1.8$ & $7.2$ & $31.4$ & $37.8$ & $25.3$ & $2.6$ & $8.3$ & $31.5$\\\cline{2-13}                                & \multirow{3}[3]*{$300$} & $\mathcal {F}_{-}$  &                                $0$ & $0$ & $69.3$ & $28.9$ & $0.1$ & $0$ & $0$ & $66.9$ & $29.3$ & $0$\\                                &                         & $\mathcal {F}_{*}$  &                                $96.9$ & $99.7$ & $28$ & $60.7$ & $59.7$ & $94.5$ & $98.9$ & $28.9$ & $59.8$ & $61.6$\\                                &                         & $\mathcal {F}_{+}$  &                                $3.1$ & $0.3$ & $2.7$ & $10.4$ & $40.2$ & $5.5$ & $1.1$ & $4.2$ & $10.9$ & $38.4$\\\hline\multirow{6}[6]*{$\frac{3}{4}$} & \multirow{3}[3]*{$80$} & $\mathcal {F}_{-}$   &$0$ & $0$ & $70.9$ & $31.7$ & $14.4$ & $0$ & $0$ & $65$ & $35.5$ & $15.1$\\                                &                         & $\mathcal {F}_{*}$  &                                $77.6$ & $89.1$ & $28.1$ & $61.5$ & $68.3$ & $65.9$ & $75.1$ & $31.1$ & $57.7$ & $66$\\                                &                         & $\mathcal {F}_{+}$  &                                $22.4$ & $10.9$ & $1$ & $6.8$ & $17.3$ & $34.1$ & $24.9$ & $3.9$ & $6.8$ & $18.9$\\\cline{2-13}                                & \multirow{3}[3]*{$200$} & $\mathcal {F}_{-}$  &                                $0$ & $0$ & $63.8$ & $29.4$ & $0.2$ & $0$ & $0$ & $66.8$ & $32.4$ & $0.3$\\                                &                         & $\mathcal {F}_{*}$  &                                $97.7$ & $99.7$ & $32.4$ & $60.1$ & $61.2$ & $96.1$ & $99.2$ & $30.2$ & $57.2$ & $59.1$\\                                &                         & $\mathcal {F}_{+}$  &                                $2.3$ & $0.3$ & $3.8$ & $10.5$ & $38.6$ & $3.9$ & $0.8$ & $3$ & $10.4$ & $40.6$\\
\hline\end{tabular}\end{table}
%\subsection{Relation to the  model $\bR_n^{1/2}\bW_n$}

%\section{Discussion}\label{diss}

%In addition to the aforementioned applications, the proposed main results is also useful in the classifications and some other statistical inferences. The details of the classifications are discussed in the Supplementary material. For some other statistical inferences, consider a high dimension mean test hypothesis:
%$H_0:\boldsymbol{\mu}=\boldsymbol{\mu}_0\quad \text{ versus }\quad H_1: \boldsymbol{\mu}\neq\boldsymbol{\mu}_0.$
%To test $H_0$,   the  asymptotic properties  of Hotelling's $T^2$ statistic is  investigated by \cite{pan2011central} under Assumption \ref{ass2}. However, if $p\gg n$, the existed theory becomes invalid. To tackle this issue, many approaches, e.g., random projection, are proposed under different conditions. By Theorem \ref{thmeve1}, we can construct a novel projection matrix to reduce the dimension of data, and use the reduced observations to do statistical inference. For more details, one can refer to \cite{yang2022center}. 

%Moreover, setting $\ba_i=\bB\by_i$ with $\bB\in\mR^{n\times K}$ and $\by_i\in\mR^{K}$ in \eqref{datamat} is also a common model in practice; see for e.g., factor models. Under such models, many statistical problems are also of interest. For instance, the estimation of the factor numbers and the factor loadings  (see \cite{fan2011high}). 

At the end of this section, we use a simple simulation to demonstrate the matching properties of the left spiked eigenvectors and spiked eigenvalues between a signal-plus-noise matrix and a sample covariance matrix, which have been discussed in Remarks \ref{vectormatch} and \ref{valuematch}. Let 
$\bA_n =\sum_{i=1}^2d_i \bg_i \br_i^*$ where $d_1=3,d_2=2$, $\bg_1=2^{-1/2}(1,1,0,\cdots, 0)^*$, $\bg_2=2^{-1/2}(-1,1,0,\cdots, 0)^*$, $\br_1$ and $\br_2$ are two right eigenvectors of a $p\times n$ Gaussian matrix.
%$\bA_n = U \Lambda V^\top \in \mathbb{R}^{p\times n}$, where $U$ has two column vectors $2^{-1/2}(1,1,0,\cdots, 0)$ and  $2^{-1/2}(-1,1,0,\cdots, 0)$, $\Lambda = \diag(3,2)$ and $V$ consists the first two right eigenvectors of a $p\times n$ Gaussian matrix, and  
Take $\bSigma_n = (0.4^{|i-j|}) + \diag(0,0,6,0,\cdots 0)\in \mathbb{R}^{p\times p}$. Let Model 1 be $\bA_n+\bSigma_n^{1/2} \bW_n$ where $\bW_n$ consists of independent $\mathcal{N}(0,1/n)$, and Model 2 be $\bR_n^{1/2}\bW_n \bW_n^* \bR_n^{1/2}$ where $\bR_n = \bA_n \bA_n^* +\bSigma_n$, and $\bW_n$ the same as Model 1. There are three spiked eigenvalues satisfying Assumption \ref{ass4new}. Table \ref{tabev} reports the three largest eigenvalues and eigenvectors of $\bX_n\bX_n^*$ with $p=100,n=200$ averaged from 500 replications generated by Model 1 and 2, respectively.  

\setcounter{table}{2}
\begin{table}[hb]
\centering
\caption{The first three eigenvalues and eigenvectors of $\bX_n\bX_n^*$ where $\bX_n$ are generated by Model 1 and 2, averaging from 500 replications each, with $v = (1,0,\cdots, 0)$
(values in parentheses indicate the standard deviations).}
\renewcommand{\arraystretch}{1.2}
\begin{tabular}{ccccccc}
\hline
    &  $\lambda_1$ &   $\lambda_2$ &   $\lambda_3$ &  $(v^* \hat{v}_1)^2$ &  $(v^* \hat{v}_2)^2$ &  $(v^* \hat{v}_3)^2$\\
 \multirow{2}{*}{$\bX_n= \bA_n + \bSigma_n^{1/2}\bW_n$}  & 11.122 & 7.574   & 5.238  &0.447 & 0.040 & 0.377   \\
   & (0.550) & (0.633)   & (0.261)  & (0.052) & (0.050) & (0.053) \\
\multirow{2}{*}{$\bX_n= \bR_n^{1/2}\bW_n(\bR_n =\bA_n\bA_n^* + \bSigma_n)$} &11.114 & 7.583 & 5.212 &0.444 & 0.047 & 0.371 \\
 & (1.026) & (0.640)   & (0.432)  & (0.089) & (0.065) & (0.080) \\
  \hline
\end{tabular}\label{tabev}
\end{table}

We observe that the first-order limits are almost the same for the two types of models. Moreover, the fluctuation behaviour is possibly different which can be inferred from the different standard deviations in Table \ref{tabev}.

\section{Proofs of main results}\label{secproofs}

\subsection{Proof of Theorems \ref{thmeve1} and \ref{thm2}}
In this section, we prove the main results in Section \ref{secModelResults} .
Proposition \ref{prop1} plays an important role in the proof of Theorem \ref{thmeve1}. To prove Proposition \ref{prop1},  the following Proposition \ref{prop01} is required, whose proof is provided in the supplementary material \cite{liu2024asymptotic}. 
%and use the conjugate transpose $``^*"$ to replace the common transpose $``^*"$, which are the same in real cases. 
%The supplementary material also includes the proof of other theoretical results in Section  \ref{app}.
%Theorem \ref{thm2}, Theorem \ref{thm3.1} and Corollary \ref{cor1} are also provided. 
\begin{prop}\label{prop01}
	Under the conditions of Proposition \ref{prop1}, for any sequence of $n\times 1$ deterministic unit norm vectors $(\bu_n)_{n\geq 1}$ and $z\in \CC^+$ with $\Im z$ being bounded %away 
from below by a positive constant, we have \begin{eqnarray}\label{tQtTdif}
{\mathbb{E}}|\bu_n^*(\tilde \bQ_n(z)-\tilde \bT_n(z))\bu_n|^2=O(n^{-1}),
\end{eqnarray} 
	where 
 \begin{equation*}\begin{aligned}
\tilde \bT_n(z)&=\left(-z(1+\delta_n(z))\bI+\bA_n^*(\bI+\tilde\delta(z)\bSigma_n)^{-1}\bA_n\right)^{-1}\\
\bT_n(z)&=\left(-z(\bI+\tilde\delta_n(z)\bSigma_n)+\frac1{1+\delta(z)}\bA_n\bA_n^*\right)^{-1},\\
 \end{aligned}\end{equation*}
 $\delta_n(z)=\frac1n\trace (\bSigma_n \bT_n(z))$ and $\tilde\delta_n(z)=\frac1n\trace(\tilde \bT_n(z))$. 
\end{prop}
The proof strategy of this result is to deal with the Gaussian case first and then use the interpolation methods to deal with $\bW$ under Assumption \ref{ass1}.  The Gaussian case is a simple consequence of the main result in \cite{hachem2013bilinear}; see the related discussion around equation (S.22) in the supplementary file. 

%To give the theoretical justifications, 
%we first introduce a necessary lemma.
We also need the following result to prove Propositon \ref{prop1}.
\begin{lemma}\label{wbm}(Woodbury matrix identity) Suppose that $A$ is $n\times n$, $D$ is ${k\times k}$, $U$ is ${n\times k}$, 
and $V$ is ${k\times n}$. If $A$ and $D$ are invertible, we have%there is
	\begin{equation*}
(A+U D V)^{-1}=A^{-1}-A^{-1} U\left(D^{-1}+V A^{-1} U\right)^{-1} V A^{-1}.
\end{equation*}
\end{lemma}

Now we start to prove Proposition \ref{prop1}.\\

\noindent \textbf{Proof of Proposition \ref{prop1}.}
We omit the subscripts ``n'' for ease of notation.
First, we prove \eqref{key1}. 
Assuming the validity of Proposition 3, it suffices to prove that for any $z \in \mathcal{C}^+$
\begin{eqnarray}\label{tTtRdif}
		\left|\bu^*\left[\tilde \bT(z) - \tilde \bD(z)\right]\bu\right|=O(n^{-1/2}).
	\end{eqnarray}

Using the Woodbury matrix identity in Lemma \ref{wbm}, we have \begin{equation}\label{tildeTz1}\tilde \bT(z)=-\frac1{z(1+\delta(z))}\bI-\left(-\frac1{z(1+\delta(z))}\right)^2\bA^*\left[\bI+\tilde \delta(z)\bSigma+\frac{-1}{z(1+\delta(z))}\bA\bA^* \right]^{-1}\bA.\end{equation}
Define
 $$\tilde{\boldsymbol{\Delta}}(z)=\tilde\delta(z)\bI-(\tilde\delta(z))^2\bA^*\left[\bI+\tilde \delta(z)\left(\bSigma+\bA\bA^*\right) \right]^{-1}\bA. $$
We have \begin{eqnarray}\label{ineq_p1}
 &&\left|\bu^*\left(\tilde \bT(z)-\tilde {\boldsymbol{\Delta}}(z)\right)\bu\right|\le \left|-\frac1{z(1+\delta(z))}-\tilde \delta(z)\right||\bu^*\bu|\nonumber\\
 &&+\left|\left(-\frac1{z(1+\delta(z))}\right)^2-\tilde \delta^2(z) \right|\left|\bu^*\left( \bA^*\left[\bI+\tilde \delta(z)\left(\bSigma+\bA\bA^*\right) \right]^{-1}\bA\right)\bu \right|\nonumber\\
 &&+\left|\left(-\frac1{z(1+\delta(z))}\right)^2\right|\Bigg|\bu^*\Bigg(\bA^*\left[\bI+\tilde \delta(z)\left(\bSigma+\bA\bA^*\right) \right]^{-1}\bA\nonumber\\
 &&~-\bA^*\left[\bI+\tilde \delta(z)\bSigma+\frac{-1}{z(1+\delta(z))}\bA\bA^* \right]^{-1}\bA\Bigg)\bu\Bigg|.
 \end{eqnarray}
We first consider  the convergence rate of  \begin{equation}\label{appro}
 -\frac1{z(1+\delta(z))}-\tilde \delta(z).
 \end{equation}
By \eqref{tildeTz1} we have \begin{eqnarray}\label{343wo}
 	-\frac1{z(1+\delta(z))}-\tilde \delta(z)=\frac{1}{n}\left(\frac1{z(1+\delta(z))}\right)^2\trace\bA^*\bT(z)\bA.
 \end{eqnarray}
Proposition 2.2 in \cite{hachem2007deterministic} yields $\|\bT(z)\|\le \frac{1}{\Im z}$, and one can see in \cite{hachem2013bilinear} as well. Also, by Lemma 2.3 of \cite{silverstein1995empirical}, there is $\|(\bI+\tilde\delta(z)\bSigma)^{-1}\|\le\max(\frac{4}{\Im z},2)$. Combining the fact that $\trace \bA\bA^*=O(n^{1/3})$, we have \begin{equation}\label{tdelta_delta}
 	|-\frac1{z(1+\delta(z))}-\tilde \delta(z)|=O\left(\frac1{n^{2/3}(\Im z)^3}\right).
 \end{equation}
  Thus, a direct calculation shows that  \begin{eqnarray}\label{tdel1}
 \left|u^*\left(\tilde \bT(z)-\tilde{\boldsymbol{\Delta}}(z)\right)v\right|\le O\left(\frac1{n^{2/3}(\Im z)^7}\right).
 \end{eqnarray}

 Then, to conclude the bound in \eqref{tTtRdif}, we require a bound on the difference between $\tilde{\delta}$ and $\tilde{r}$, where $\tilde{r} = \tilde{r}_n$ is defined through \eqref{32ra3}.
% Next, let \ $$ \tilde R(z)=\tilde r(z)\bI-(\tilde r(z))^2\bA^*\left[\bI+\tilde r(z)\left(\bSigma+\bA\bA^*\right) \right]^{-1}\bA,$$where $\tilde{r}(z) $ in $\mathcal{C}^{+}$ solves the equation \begin{eqnarray*}z=-\frac{1}{\tilde{r}(z)}+c_n\int\frac{t{\color{blue}\mathrm{d}}H^{\bR_n}(t)}{1+t\tilde{r}(z)},\end{eqnarray*}
 %and $H^{\bR_n}(t)$ is the empirical spectral distribution of $\bR_n=\bSigma+\bA\bA^*$. 
 If we denote the right hand side of (\ref{343wo}) by $\omega$, then (\ref{343wo}) can be rewritten as $z=\frac{1}{\tilde{\delta}}-z\delta+\omega_1,$ where $\omega_1=-\frac{1}{\tilde{\delta}}-\frac{1}{\tilde{\delta}+\omega}$.
 We also let  $$\bT'(z)=\left(-z(\bI+\tilde\delta(z)\bSigma)-z\tilde \delta(z)\bA\bA^*\right)^{-1}.$$
 By the definition of $\delta$, this equation can be further written as \begin{eqnarray}\label{34rb2}
 \begin{aligned}
 z &=-\frac{1}{\tilde{\delta}}-\frac{z}{n}\trace\bSigma \bT +\omega_1\\&=
-\frac{1}{\tilde{\delta}}-\frac{z}{n}\trace(\bSigma+\bA\bA^*)\bT+\frac{z}{n}\trace\bA\bA^*\bT+\omega_1  \\&=-\frac{1}{\tilde{\delta}}-\frac{z}{n}\trace(\bSigma+\bA\bA^*)\bT'+\frac{z}{n}\trace(\bSigma+\bA\bA^*)(\bT'-\bT)+\frac{z}{n}\trace\bA\bA^*\bT+\omega_1\\&=-\frac{1}{\tilde{\delta}}+c_n\int \frac{t{\mathrm{d}}H^{\bR_n}(t)}{1+t\tilde{\delta}}+\omega_2,
\end{aligned}
 \end{eqnarray}
 where $\omega_2=\omega_1+\frac{z}{n}\trace(\bSigma+\bA\bA^*)(\bT'-\bT)+\frac{z}{n}\trace \bA\bA^*\bT.$
We have that $ |\omega_1| = O\left(n^{-2/3}(\Im z)^{-5}\right)$, $|\frac{z}{n}\trace(\bSigma+\bA
\bA^*)(\bT'-\bT)|=O\left(n^{-2/3}(\Im z)^{-5}\right)$,and $|\frac{z}{n}\trace\bA\bA^*\bT|=O\left(n^{-2/3}(\Im z)^{-1}\right)$. Then it follows that $|\omega_2|=O\left(n^{-2/3}(\Im z)^{-5}\right)$.
 With equations (\ref{32ra3}) and (\ref{34rb2}) at hand,  we have
  \begin{eqnarray*}
 	\tilde{\delta}-\tilde{r}&=&(\tilde{\delta}-\tilde{r})\left(\tilde{\delta}\tilde{r}c_n\int\frac{t^2{\mathrm{d}}H^{\bR_n}(t)}{(1+t\tilde r)(1+t\tilde\delta)}\right)-\tilde\delta\tilde r\omega_2.
 \end{eqnarray*}
Similar to (6.2.26) in \cite{bai2010spectral}, we also have $$\left|\tilde{\delta}\tilde{r}c_n\int\frac{t^2{\mathrm{d}}H^{\bR_n}(t)}{(1+t\tilde r)(1+t\tilde\delta)}\right|\le 1-C(\Im z)^2.$$
 Therefore, we obtain  
\begin{equation*}|\tilde{\delta}-\tilde{r}|=O\left(\frac1{n^{2/3}(\Im z)^7}\right).
\end{equation*}
Using the same arguments as in  \eqref{ineq_p1}, it follows that \begin{eqnarray*}\label{rdel001}
 | \bu^*(\tilde \bD(z)- \tilde{\boldsymbol{\Delta}} (z))\bu |=O\left(\frac1{n^{2/3}(\Im z)^{11}}\right).
 \end{eqnarray*}
Combining this with \eqref{tdel1} yields \eqref{tTtRdif} for any fixed $z \in \mathcal{C}^+$.
This concludes the proof of \eqref{key1}. 

The proof of \eqref{Qnlim} is very similar by using ${\mathbb{E}}|\bv_n^*(\bQ_n(z)-\bT_n(z))\bv_n|^2=O(n^{-1})$ and \eqref{tdelta_delta}.
The former bound can be verified similarly to its conjugate version \eqref{tQtTdif} and even more simply; see (S.19) of the supplementary material for comparison. 
\qed

To prove Theorem \ref{thmeve1}, we also need the following result on the exact separation of the eigenvalues of $\bS_{n}$. This result improves Theorem 1 of \cite{liu2022ieee} by relaxing the condition on $\bA$, see the discussion in Remark \ref{valuematch}, and its proof is in the supplementary file. Recall that $\bR_n=\bA_n\bA_n^*+\bSigma=\sum_{j=1}^{\tilde{n}}\gamma_j \boldsymbol\Xi_j \boldsymbol\Xi_j^*$.

\begin{lemma}\label{sepl1} Assume that Assumptions \ref{ass1},\ref{ass2}, and \ref{ass3new} hold. Let $(a, b)$ be an interval with $a>0$ and lie in an open interval outside the support of $F^{c_n,H^{\bR_n}}$ for all large $n$.  For $0\le k \le \tilde{n}-1$, denote $\ell_k=\sum_{i=1}^k m_{i}$ with the convention $\ell_0=0$. If  $[-\tilde r(a)^{-1},-\tilde r(b)^{-1}]\subset (\gamma_{k+1},\gamma_k)$  where $\tilde r(z)$ are given in \eqref{32ra3}, with the convention 
$\gamma_0=\hat{\lambda}_0=\infty$, we have $$\mathbb{P}(\hat{\lambda}_{\ell_k}>b \text{ and }\hat{\lambda}_{\ell_k+1}<a, \; \text{for all large} \; n) = 1.
%\text{ as $n\rightarrow\infty$},
$$

%where $\ell_k=\sum_{i=1}^k m_{i}$ and $\lambda_j$ is the $j$-th largest eigenvalue of $\bS_{n}$.
\end{lemma}

\noindent \textbf{Proof of Theorem \ref{thmeve1}.}
We first prove \eqref{l2normvandhatu}. %The key step to conclude \eqref{l2normvandhatu} is to find the  asymptotic limit of $\bu^* \hat{\bu}_k \hat{\bu}_k^* \bu$ for any deterministic vector $\bu$. 
Define $$\mR_y(k)=\{z\in\mC:~\hat\sigma_1\le\Re z\le \hat\sigma_2,~|\Im z|\le y\},$$
where $y>0$, $[\hat\sigma_1,\hat\sigma_2]$ encloses the sample eigenvalues $\{\hat{\lambda}_{\ell}, \ell\in \mathcal{K}_k\}$ and excludes all other sample eigenvalues. The existence of $\mR_y(k)$ is guaranteed by Assumption \ref{ass4new} and Lemma \ref{sepl1}. The reasoning is as follows. 
It is known that if $\alpha$ satisfies $\varphi'(\alpha)>0$ and $\alpha \neq 0$ is not in the support of $H$, then $\varphi(\alpha)$ is not in the support of $F^{c,H}$, see \cite{bai2010spectral} for instance.  This reasoning also applies to finding intervals outside the support of $F^{c_n,H^{\bR_n}}$ by using $\varphi_n(\gamma)$, which is the non-asymptotic version of $\varphi(\gamma)$ by replacing $c$ and $H$ with $c_n =p/n$ and $H^{\bR_n}$, respectively. By Assumption \ref{ass4new} we can find two intervals lying on the left and right side of $\gamma_k$, respectively, and that they lie outside the support of $F^{c,H}$ and $F^{c_n,H^{\bR_n}}$ for all large $n$.  According to Lemma \ref{sepl1} we can find an interval enclosing only $\{\hat{\lambda}_{\ell}, \ell\in \mathcal{K}_k\}$ with probability one for sufficiently large $n$. % We refer the reader to the proof of Theorem 2 below for a more detailed explanation. 

By the Cauchy integral formula, we have\begin{equation}\label{cif} \frac{1}{2\pi i}\oint_{\partial \mR_y^-(k)}\bu^*\tilde{\bQ}_n(z)\bu {\mathrm{d}}z=\bu^*\left(\sum_{\ell\in \mathcal{K}_k}\hat{\bu}_\ell\hat{\bu}_\ell^* \right) \bu:=\hat r_k,\end{equation}
where $\bu$ is any $n\times 1$ deterministic unit vector, and $\partial \mR_y^-(k)$ represents negatively oriented boundary of $\mR_y(k)$.
\begin{lemma}\label{eigvec1}
Under Assumptions of Theorem \ref{thmeve1}, we have 
$$\left|\hat r_k-\frac{1}{2\pi i}\oint_{\partial \mR_y^-(k)}\bu^*\tilde \bD_n(z)\bu {\mathrm{d}}z\right| =O_P \left(\frac{1}{\sqrt{n}}\right).$$
where $\tilde \bD(z)$ is defined in \eqref{Tz}.
\end{lemma}
\begin{proof}
The proof is in the same spirit as that of Proposition 1 in \cite{mestre2008asymptotic}. Since our result provides a convergence rate of error, we use a slightly different argument by considering the second moment of the left term.
 Define an event $\Omega:=\{ \hat{\sigma}_1+\delta < \hat{\lambda}_\ell < \hat{\sigma}_2-\delta , \ell\in \mathcal{K}_k \}$,  which holds almost surely for some small positive $\delta>0$ independent of $n$. We have
\begin{equation}\label{secbd234}\begin{aligned}
   & \mathbb{E} \left|\oint_{\partial \mR_y^-(k)} \left(\bu^*(\tilde \bQ_n(z)-\tilde \bD_n(z))\bu \right)I(\Omega) {\mathrm{d}}z \right|^2\\
   & \leq C  \oint_{\partial \mR_y^-(k)}\mathbb{E} \left(|\bu^*(\tilde \bQ_n(z)-\tilde \bD_n(z))\bu|^2 I(\Omega)\right) |{\mathrm{d}}z| = O(n^{-1})
\end{aligned}\end{equation}
where the first step uses H\"older's inequality and the second step follows from  Proposition \ref{prop01} and Lemma \ref{sepl1}.
Then the conclusion follows from Chebyshev's inequality.
\end{proof}
The above lemma reduces the proof to calculating the deterministic integral $$F=\frac{1}{2\pi i}\oint_{\partial \mR_y^-(k)} \bu^*\tilde \bD(z)\bu {\mathrm{d}}z.$$ Let $w(z)=-\frac1{\tilde r(z)}$, where $\tilde r(z)$ is introduced in Proposition \ref{prop1}. We find that $w(z)$ satisfies the following equation $$z=w(z)\left(1-c\int\frac{t{\mathrm{d}}H^{\bR_n}(t)}{t-w(z)}\right),$$
which is parallel to equation (24) in \cite{mestre2008improved}. Thus, $w(z)$ satisfies all the properties listed in Proposition 2 in \cite{mestre2008improved}. Write $F=F_1+F_2$, where
\begin{equation}\label{defF1}
F_1=-\frac{1}{2\pi i}\bu^*\bu\oint_{T^-(k)}\frac1w\left[1-\frac1n\sum^L_{k=1}m_k\left(\frac{\gamma_k}{\gamma_k-w}\right)^2 \right]{\mathrm{d}}w,	
\end{equation}
\begin{equation}\label{defF2}
F_2	=-\frac{1}{2\pi i}\oint_{T^-(k)}\frac1w\bu^*\bA^*\sum_{k=1}^L\frac{\boldsymbol\Xi_k\boldsymbol\Xi_k^*}{w-\gamma_k}\bA\bu\left[1-\frac1n\sum^L_{k=1}m_k\left(\frac{\gamma_k}{\gamma_k-w}\right)^2 \right]{\mathrm{d}}w	
\end{equation}
and where $T^-(k)$ is a simple closed curve that includes $\gamma_k$ and excludes all the other population eigenvalues of $\bR_n$ with negative orientation. By a calculation, $$F_1=Res\left(\frac1w\left[1-\frac1n\sum^L_{k=1}m_k\left(\frac{\gamma_k}{\gamma_k-w}\right)^2 \right],\gamma_k\right)=\frac{m_k}{n}.$$
For $F_2$, we further decompose the integrand as $$F_2=-\frac1{2\pi i}\oint_{T^-(k)}(\chi_{1k}(w)+\chi_{2k}(w)+\chi_{3k}(w)+\chi_{4k}(w)){\mathrm{d}}w,$$
where \begin{eqnarray*}
 \chi_{1k}(w)&=&\frac{\bu^*\bA^*\boldsymbol\Xi_k\boldsymbol\Xi_k^*\bA \bu}{w(w-\gamma_k)}, ~\chi_{2k}(w)=-\frac{m_k\gamma_k^2}{n}\frac{\bu^*\bA^*\boldsymbol\Xi_k\boldsymbol\Xi_k^*\bA \bu}{w(w-\gamma_k)^3},\\
 \chi_{3k}(w)&=&-\frac{\bu^*\bA^*\boldsymbol\Xi_k\boldsymbol\Xi_k^*\bA \bu}{nw(w-\gamma_k)}\sum_{i=1,i\neq k}^{\tilde{n}}m_i\left(\frac{\gamma_i}{\gamma_i-w}\right)^2,\\
  \chi_{4k}(w)&=&-\frac1{nw}\bu^*\bA^*\sum_{i=1,i\neq k}^{\tilde{n}}\frac{\boldsymbol\Xi_i\boldsymbol\Xi_i^*}{w-\gamma_i}\bA \bu\frac{m_k\gamma_k^2}{(\gamma_k-w)^2}.
 \end{eqnarray*}
By calculation, we have
\begin{eqnarray*}
 Res(\chi_{1k}(w),\gamma_k)&=&\frac{\bu^*\bA^*\boldsymbol\Xi_k\boldsymbol\Xi_k^*\bA \bu}{\gamma_k}, ~ Res(\chi_{2k}(w),\gamma_k)=-\frac{m_k\bu^*\bA^*\boldsymbol\Xi_k\boldsymbol\Xi_k^*\bA \bu}{n\gamma_k},\\
 Res(\chi_{3k}(w),\gamma_k)&=&-\frac{\bu^*\bA^*\boldsymbol\Xi_k\boldsymbol\Xi_k^*\bA \bu}{n\gamma_k}\sum_{i=1,i\neq k}^{\tilde{n}}m_i\left(\frac{\gamma_i}{\gamma_i-\gamma_k}\right)^2,\\
  Res(\chi_{4k}(w),\gamma_k)&=&-\frac1n \bu^*\bA^*\sum_{i=1,i\neq k}^{\tilde{n}}\frac{\boldsymbol\Xi_i\boldsymbol\Xi_i^*m_k(\gamma_i-2\gamma_k)}{(\gamma_k-\gamma_i)^2}\bA \bu.
 \end{eqnarray*}
Therefore,  we have $$F=\frac{\bu^*\bA^*\boldsymbol\Xi_k\boldsymbol\Xi_k^*\bA \bu}{\gamma_k}\left(1-\frac1n\sum_{i=1,i\neq \tilde{n}}\frac{m_i\gamma_i^2}{(\gamma_k-\gamma_i)^2}\right)+O\left(\frac{m_k}{n}\right).$$
Letting $\eta_k=\left(1-\frac1n\sum_{i=1,i\neq k}\frac{\gamma_i^2}{(\gamma_k-\gamma_i)^2}\right)$ and by the assumption that $m_k=O(n^{1/3})$, we conclude \eqref{l2normvandhatu}. 

The first assertion can be obtained by an argument similar to the one that leads to \eqref{secbd234} and the calculation of the deterministic term is exactly the same as (22) in \cite{mestre2008asymptotic}, where their $d_k(z)$ lines up with  $(-\gamma_k z\tilde{r}(z) - z)^{-1}$ in our case.
\qed 
\\

\noindent \textbf{Proof of Theorem \ref{thm2}.}
We first consider \eqref{spikaslim}.
Denote the support of $H$ by  $\Gamma_H$. Under Assumption \ref{ass4new}, it is easy to obtain that $\varphi'(\gamma_k)>0$ for $1\leq k\leq L$. By the continuity of $\varphi$, there exists $\delta>0$ such that
\begin{equation}\label{equ5.2}
\varphi'(x)>0,\quad\forall x\in(\gamma_{k}-\delta,\gamma_{k}+\delta)
\end{equation}
and $\gamma_{k+1}<\gamma_{k}-\delta<\gamma_{k}+\delta<\gamma_{k-1}$ (by default,  $\gamma_{0}=\infty$).  Then, we can find
 $0<\varepsilon<\delta$ and $\gamma_{k}-\delta<a<b<\gamma_{k}-\varepsilon<\gamma_{k}<\gamma_{k}+\varepsilon<e<f<\gamma_{k}+\delta$ 
 such that $[a,b]$ and $[e,f]$ are outside $\Gamma_H$. Let $K = \sum_{i=1}^L m_i$. For $\gamma\in[a,b]\cup[e,f]$, write
\begin{equation*}
\varphi_{n}(\gamma)=\gamma+\gamma c_{n}\left(\frac{p-K}{p}\int\frac{t}{\gamma-t}\mathrm{d}H_{n}^{\text{Non}}(t)+\frac{1}{p}\sum_{k=1}^{L}\frac{m_k \gamma_{k}}{\gamma-\gamma_{k}}\right),
\end{equation*}
where $H_{n}^{\text{Non}}(t)=\frac{1}{p-K}\sum_{j=K+1}^{p}I_{[\gamma_{j},\infty)}(t)$ is the ESD of nonspikes. Then,
\begin{equation}\label{equ2.3}
\begin{split}
\varphi_{n}(\gamma)-\varphi(\gamma)=&c\gamma\int\frac{t}{\gamma-t}\mathrm{d}H^{\text{Non}}_{n}(t)-c\gamma\int\frac{t\mathrm{d}H(t)}{\gamma-t}+\frac{c_{n}}{p}\gamma\sum_{k=1}^{L}\frac{m_k\gamma_{k}}{\gamma-\gamma_{k}}\\
&+\left(c_{n}\frac{p-K}{p}-c\right)\gamma\int\frac{t}{\gamma-t}\mathrm{d}H^{\text{Non}}_{n}(t).
\end{split}
\end{equation}
Observe that $$\inf\limits_{L+1\leq j\leq p,\gamma\in[a,b]\cup[e,f]}|\gamma_{j}-\gamma|>0 \text{ and } \inf\limits_{1\leq k\leq L,\gamma\in[a,b]\cup[e,f]}|\gamma_{k}-\gamma|>0,$$ so that the third and the fourth term on the right hand of (\ref{equ2.3}) converge uniformly to zero, as $p\rightarrow\infty$. It is shown that the first term on the right hand of (\ref{equ2.3}) converges pointwise to the second one, in which they are all continuous function w.r.t. $\gamma$. Since $\{c\gamma\int\frac{t}{\gamma-t}\mathrm{d}H^{\text{Non}}_{n}(t)\}$ can be regarded as a monotone sequence of functions,
by Dini's theorem, the convergence is uniform. Thus, $\varphi_{n}$ uniformly converges to $\varphi$ on $[a,b]\cup[e,f]$. The proof for the uniform convergence of $\varphi'_{n}$ equal to
\begin{equation*}
%\begin{split}
\varphi'_{n}(\gamma)=1-c_{n}\left(\frac{p-K}{p}\int\frac{t^{2}\mathrm{d}H^{\text{Non}}_{n}(t)}{(\gamma-t)^2}+\frac{1}{p}\sum_{k=1}^{L}\frac{m_k\gamma_{k}^2}{(\gamma-\gamma_k)^2}\right).
%&\left(\varphi^{c_{n},H_{n}}(\lambda)\right)'=1+c_{n}\left(\frac{p-K}{p}\int\frac{(\lambda-t)t-\lambda t}{(\lambda-t)^2}\mathrm{d}H_{n}^{\text{Non}}(t)+\frac{1}{p}\sum_{r=1}^{l}\frac{n_{r}\alpha_{r}(\lambda-\alpha_{r})-n_{r}\alpha_{r}\lambda}{(\lambda-\alpha_{r})^2}\right)\\
%\end{split}
\end{equation*}
is analogous and left out here. Hence, from Lemma 6.1 of \cite{bai2010spectral}, combining (\ref{equ5.2}) and the uniform convergence of $\varphi_{n},\varphi'_{n}$ on $[a,b]\cup[e,f]$, it follows that both $[\varphi(a),\varphi(b)]$ and $[\varphi(e),\varphi(f)]$ are out of the support of $F^{c_n,H^{\bR_n}}$. Then, using Lemma \ref{sepl1},
\begin{equation*}
\begin{split}
& \mathbb{P}(\hat{\lambda}_{\ell_{k+1}}\leq\varphi(a)<\varphi(b)\leq \hat{\lambda}_{\ell_k},~\text{for all large}~n)=1,\\
& \mathbb{P}(\hat{\lambda}_{\ell_{k-1}+1}\leq\varphi(e)<\varphi(f)\leq \hat{\lambda}_{\ell_{k-1}},~\text{for all large}~n)=1.
\end{split}
\end{equation*}
Hence, with probability  one,
\begin{equation*}
\varphi(b)\leq\liminf\limits_{n}\hat{\lambda}_{\ell_k},\quad \limsup\limits_{n}\hat{\lambda}_{\ell_{k-1}+1}\leq\varphi(e),
\end{equation*}
Finally, letting $b\uparrow\gamma_{k}$ and $e\downarrow\gamma_{k}$, we have
\begin{equation}\label{equ2.1}
\varphi(\gamma_{k})\leq\liminf\limits_{n}\hat{\lambda}_{\ell_k},\quad \limsup\limits_{n}\hat{\lambda}_{\ell_{k-1}+1}\leq\varphi(\gamma_{k})\quad \text{with probability  one.}
\end{equation}
From (\ref{equ2.1}), we conclude that for any $\ell \in \mathcal{K}_k$, with probability  one, 
\begin{equation*}
\lim\limits_{n\rightarrow\infty}\hat{\lambda}_{\ell}=\varphi(\gamma_{k}),~k=1,2,\ldots,L.
\end{equation*}

Next we turn to the second assertion \eqref{equ2.5new}. Define $\varphi_n$ by replacing $c$ and $H$ defined in $\varphi$, see \eqref{equ1.3}, with $c_n = p/n$ and $H^{\bR_n}$. The support of $F^{c_n,H^{\bR_n}}$ can be characterized as
 $\cup_{i=1}^{n_0} [a_{2i},a_{2i-1}]$
where $a_{i}$ are critical points of $\varphi_n$, see Lemma 2.6 of \cite{knowles2017anisotropic}. For those spiked eigenvalues, it can be verified that the corresponding support satisfies $|a_{2i-1}-a_{2i}| \to 0$; and for nonspiked eigenvalues, there may exist just one bulk component, or several bulk components that are disjoint with each other. In the latter case, we assume the distance between adjacent bulks is bounded below by a constant, otherwise, we merge the adjacent bulks into one bulk.

The classical number of eigenvalues in the $i$-th bulk component is $n_i := \int_{a_{2i}}^{a_{2i-1}}dF^{c_n,H^{\bR_n}}$. According to Lemma A.1 of \cite{knowles2017anisotropic}, $n_i = \sum_{j=1}^{\tilde{n}} m_j I\left(-[\tilde{r}(a_{2i})]^{-1}< \gamma_j < -[\tilde{r}(a_{2i-1})]^{-1}\right)$.
 There exists a small constant $\epsilon_0$ such that 
 $\varphi_n' ([a_{2i+1}+\epsilon_0,a_{2i}-\epsilon_0])>0$.   By Lemma \ref{sepl1}, if $\gamma_{\ell}>-[\tilde{r}(a_{2i}-\epsilon_0)]^{-1}$ and $\gamma_{\ell+1}<-[\tilde{r}(a_{2i+1}+\epsilon_0)]^{-1}$, we have $\lambda_{\ell}>a_{2i}-\epsilon_0$ and $\lambda_{\ell+1}<a_{2i+1}+\epsilon_0$. Therefore the number of sample eigenvalues in $i$-th bulk is the same as the classical number of eigenvalues $n_i$. Then by the weak convergence of $F^{\bS_n}$ to $F^{c,H}$, which is the limiting measure of $F^{c_n,H^{\bR_n}}$, and the continuity of the density of $F^{c,H}$ inside each bulk component, we conclude \eqref{equ2.5new}.
\qed

\subsection{An example where $\Sigma_n$ contains a spike}\label{sec:example}
In this section, we consider an example where $\bSigma_n$ contains a spike to illustrate that the sample spiked eigenvalues can be sourced from either $\bA_n$ or $\bSigma_n$. Moreover, the overlap between the eigenspace associated with $\bA_n\bA_n^*$ and $\bSigma_n$ can significantly influence the behavior of the sample spiked eigenvalues.

 We consider $\bX = d \bg \br^* + \bSigma^{1/2}\bW$, where $\bg, \br$ are two non-random unit norm vectors and $\bSigma = \operatorname{diag}(\ell+1,1, \cdots, 1)$. This leads to $\bR_n = d^2 \bg \bg^* + \bSigma$ according to the definition of $\bR_n$ in our Assumption \ref{ass3new}. By Wely's inequality, it can be checked that $\bR_n$ has an eigenvalue one with multiplicity  $p-2$. To obtain the other two eigenvalues of $\bR_n$, denoting $\sigma_i$ as the $i$-th diagonal element of $\bSigma$ and $g_i$ as the $i$-th coordinate of $\bg$, we find 
$$
\begin{aligned}
\operatorname{det}\left(\bSigma+d^2 \bg \bg^*-\gamma \bI_p\right) & =\operatorname{det}(\bSigma-\gamma \bI_p) \operatorname{det}\left(\bI_p+d^2(\bSigma-\gamma \bI_p)^{-1} \bg \bg^* \right) \\ %=\prod_{i=1}^p\left(\sigma_i-\gamma\right)\left(1+d^2 \bg^* (\bSigma-\gamma \bI_p)^{-1} \bg \right)\\
& =\prod_{i=1}^p\left(\sigma_i-\gamma\right)\left(1+d^2\sum_{i=1}^p \frac{g_i^2}{\left(\sigma_i-\gamma\right)}\right)
\end{aligned}
$$ 
where in the first equation we assume the invertibility of $\bSigma- \gamma I_p$, and in the second step we use the identity $\operatorname{det}(I+AB)=\operatorname{det}(I+BA)$.
Therefore, $\bR_n$ has an eigenvalue one with multiplicity  $p-2$, and the other two eigenvalues are the solution to the equation
\begin{equation*}
    1 + d^2 \frac{g_1^2}{\ell+1-\gamma}+ d^2 \frac{1-g_1^2}{1-\gamma}=0,
\end{equation*}
%This reduces to $\gamma^2 - (\ell+2 + d^2) \gamma + (1+d^2)(\ell+1) - d^2 g_1^2 \ell = 0$.
given by
\begin{equation*}
    \gamma_{1} = 2^{-1}\left(\ell+2+d^2 +  \sqrt{(d^2-\ell)^2+4 d^2 g_1^2 \ell} \right), \quad \gamma_{2} = 2^{-1}\left(\ell+2+d^2 - \sqrt{(d^2-\ell)^2+4 d^2 g_1^2 \ell} \right).
\end{equation*}
We find that the function $\varphi(\gamma)$ defined in (2.10) in the manuscript is $\varphi(\gamma):=\gamma\left(1+c/(\gamma-1)\right)$.
Therefore, we expect a BBP transition result as follows:  for $i=1$ or 2, if $\gamma_i>1+\sqrt{c}$, then $\hat{\lambda}_i \to \varphi(\gamma_i)$ almost surely; otherwise $\hat{\lambda}_i \to (1+\sqrt{c})^2$ almost surely.

We calculate the solutions $\gamma_i, i=1,2$ in some special cases below to illustrate the influence on the overlap between $\bg$ and the eigenvector associated with $\ell+1$ of $\bSigma$. Assume $d^2>\ell>0$.   When $g_1 =o(n^{-\alpha})$ for some $\alpha>0$,
we find $\gamma_1 = d^2+1+o(n^{-2\alpha})$ and $\gamma_2 = \ell+1+o(n^{-2\alpha})$. Then we have the following with probability one:
\begin{itemize}
    \item if $d^2>\ell> \sqrt{c}$, $\hat{\lambda}_1 \to \varphi(d^2+1)$ and $\hat{\lambda}_2 \to \varphi(\ell+1)$;
    \item if $d^2>\sqrt{c}\ge \ell$, $\hat{\lambda}_1 \to \varphi(d^2+1)$ and 
    $\hat{\lambda}_2 \to (1+\sqrt{c})^2$. The second largest eigenvalue is not an outlier eigenvalue;
    \item if $\sqrt{c}\ge d^2>\ell $, there are no outlier eigenvalues. 
\end{itemize}
We can carry similar analysis when $g_1 = 1$ (or $1-o(n^{-\alpha})$). This implies $\bg$ is aligned with $\be_1$, leading to one population spiked eigenvalue equal $d^2+\ell+1$. Hence the first sample eigenvalue will tend to $\varphi(d^2+\ell+1)$ almost surely if $d^2+\ell>\sqrt{c},$ and is not sample spiked eigenvalue if $d^2+\ell\le \sqrt{c}.$

\begin{figure}
  \begin{minipage}{0.5\textwidth}
    \centering
    \includegraphics[width=0.8\linewidth]{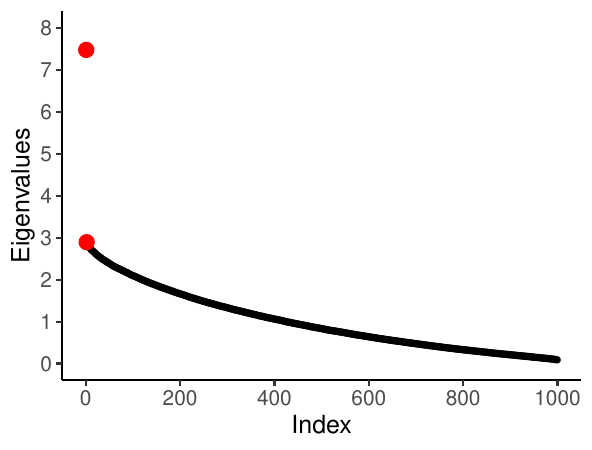}
    %\caption{First Figure}
    \label{fig:figure1}
  \end{minipage}%
  \begin{minipage}{0.5\textwidth}
    \centering
    \includegraphics[width=0.8\linewidth]{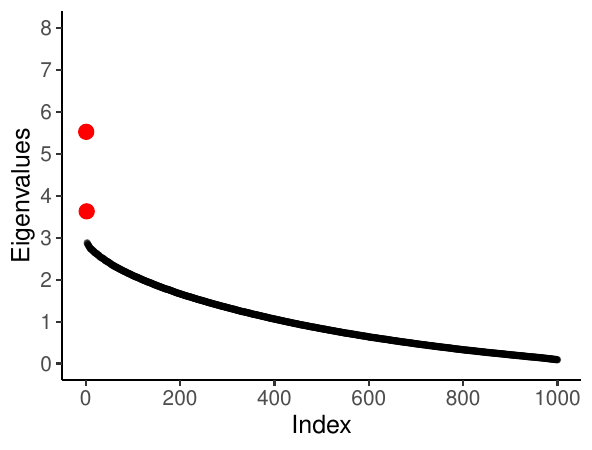}
   % \caption{Second Figure}
    \label{fig:figure2}
  \end{minipage}
  \caption{Scree plots of eigenvalues of $\bX\bX^*$. Here $\bX_{p\times n} = 2\bg \br^* + \bSigma^{1/2}\bW$ where $p=1000,n=2000$,
$\bSigma=\operatorname{diag}\{3,1, \cdots, 1\} 
$, $\bW$ consists of i.i.d. $\mathcal{N}(0,1/n)$ entries, $\br$ is a unit vector generated by normalizing a vector with $Unif[0,1]$ entries. The left panel corresponds to $\bg = (1,0,\cdots, 0)$ and the right panel corresponds to $\bg=(0,1,0,\cdots,0)$. }
  \label{fig:spi}
\end{figure}

Figure \ref{fig:spi} demonstrates the above results by plotting the eigenvalues of the considered model with $p=1000, n =2000$, $d=2$ and $\ell=2$. The left panel shows the result when $\bg = (1,0,\cdots,0)^*$ and $\bg = (0,1,0,\cdots,0)^*$, respectively. The largest two eigenvalues are red. According to the above results, the limits of the two largest eigenvalues in the left panel case are 7.583, 2.914, and 5.625, 3.750 in the right panel case.

\section{Appendix}\label{append}
This section provides the criteria of the estimation of the number of clusters when $c>1$, some additional simulation studies, and the remaining proof of the theoretical results.

\subsection{The estimation of the number of clusters when $c>1$.}
In this section, we consider the case when $p,n\rightarrow\infty$ such that $p/n\rightarrow c\in(1,\infty)$.
%and assume that Condition \ref{con3.5} with $n$ instead of $p$ hold.
Then the smallest $(p-n)$ eigenvalues of $\bm{S}_{n}$ are zero, that is,
\begin{equation*}
\hat{\lambda}_1\geq \hat{\lambda}_2\geq\ldots\geq \hat{\lambda}_K>\hat{\lambda}_{K+1}\geq\ldots\geq \hat{\lambda}_{n-1}\geq \hat{\lambda}_{n}\geq \hat{\lambda}_{n+1}=\ldots=\hat{\lambda}_{p}=0.
\end{equation*}
The modified criteria $\acute{\text{EDA}}_{k}$ and $\acute{\text{EDB}}_{k}$ for selecting the true number of clusters under $c>1$ are obtained by replacing the second term in $\text{EDA}_{k}$ and $\text{EDB}_{k}$ with $(n-k-1)\log\widetilde{\theta}_{n,k}$:
\begin{eqnarray*}
\acute{\text{EDA}}_{k}&=&-n(\hat{\lambda}_1-\hat{\lambda}_{k+1})+n(n-k-1)\log\widetilde{\theta}_{n,k}+2pk,\\
\acute{\text{EDB}}_{k}&=&-n\log(p)\cdot(\hat{\lambda}_1-\hat{\lambda}_{k+1})+n(n-k-1)\log\widetilde{\theta}_{n,k}+(\log n)pk
\end{eqnarray*}
where $\theta_{k}=\exp\{\hat{\lambda}_k-\hat{\lambda}_{k+1}\},k=1,2,\ldots,n-1,~\widetilde{\theta}_{n,k}=\frac{1}{n-k-1}\sum_{i=k+1}^{n-1}\theta_{i}^2$, called pseudo-EDA and pseudo-EDB, respectively. Analogous to the case where $0<c<1$, the modified pseudo-EDA and pseudo-EDB select the number of clusters by
\begin{eqnarray*}
\hat{K}_{\acute{\text{EDA}}}&=&\arg\min\limits_{k=1,\ldots,w}\frac{1}{n}\acute{\text{EDA}}_{k},\\
\hat{K}_{\acute{\text{EDB}}}&=&\arg\min\limits_{k=1,\ldots,w}\frac{1}{n}\acute{\text{EDB}}_{k}.
\end{eqnarray*}

The corresponding gap conditions for pseudo-EDA and pseudo-EDB stay the same as in (\ref{equ3.6}) and (\ref{equ3.6.5}), respectively. The following theorems show that $\hat{K}_{\acute{\text{EDA}}}$ and  $\hat{K}_{\acute{\text{EDB}}}$ possess a similar property as $\hat{K}_{\text{EDA}}$ and $\hat{K}_{\text{EDB}}$.

% Theorem \ref{thm3.1} and Theorem \ref{thm3.2} discuss the case when $0<c<1$. We also investigate the consistency of criteria of \eqref{equ3.2} and \eqref{equ3.2.5} when $c>1$.
% Replacing  $\theta_{k}$ and $\widetilde{\theta}_{p,k}$ in \eqref{cr1} by $\theta_{k}=\exp\{l_{n,k}-l_{n,k+1}\}$ for $k=1,2,\ldots,n-1$ and $\widetilde{\theta}_{n,k}=\frac{1}{n-k-1}\sum_{i=k+1}^{n-1}\theta_{i}^2$, and denoting $\hat{K}_{\acute{\text{EDA}}}$ and $\hat{K}_{\acute{\text{EDB}}}$ by the estimator of the replaced \eqref{equ3.2} and \eqref{equ3.2.5}, one can also have following theorem:

\begin{theorem}\label{thm3.3}
Under conditions of Theorems \ref{thm3.1}  and \ref{thm3.2}, we have the following consistency results of the estimation criteria $\hat{K}_{\acute{\text{EDA}}}$ and $\hat{K}_{\acute{\text{EDB}}}$.\\
(i) Suppose that $\gamma_{1}$ is bounded. If the gap conditions (\ref{equ3.6}), (\ref{equ3.6.5}) do not hold, then $\hat{K}_{\acute{\text{EDA}}}$ and $\hat{K}_{\acute{\text{EDB}}}$ are not consistent. If the gap conditions (\ref{equ3.6}) and (\ref{equ3.6.5}) hold, then $\hat{K}_{\acute{\text{EDA}}}$ and $\hat{K}_{\acute{\text{EDB}}}$ are strongly consistent.\\
(ii) Suppose that $\gamma_{K}$ tends to infinity. Then, $\hat{K}_{\acute{\text{EDA}}}$ and $\hat{K}_{\acute{\text{EDB}}}$ are strongly consistent.
\end{theorem}

\begin{rmk}
To illustrate EDA and EDB, one can refer to the example below, in which the true number of clusters is two.

\noindent
\textbf{Example.}
Let $p=60$, $n=100$ ,  $\bSigma$ be $p\times p$ identity matrix and $c=3/5$ and $\varpi_1=\ldots=\varpi_p=1$. Suppose the means of two clusters are $\bmu_1=(2,0,0,\ldots,0)^{\top}$, $\bmu_2=(0,2,0,\ldots,0)^{\top}$ with equal number of observations in each cluster, that is, $n_1=n_2=50$. From Theorem \ref{thm2}, the limits of first four eigenvalues of $\bm{S}_n$ can be obtained as follows
\begin{equation}\label{exeq1}
\hat{\lambda}_1,~\hat{\lambda}_{2}\rightarrow\varphi(3)=3.9~~\text{i.p.},\quad\quad \hat{\lambda}_3,~\hat{\lambda}_{4}\rightarrow(1+\sqrt{3/5})^2~~\text{a.s.}.
\end{equation}
Then,
\begin{equation}\label{exeq2}
\begin{split}
& \theta_1=\exp\{\hat{\lambda}_1-\hat{\lambda}_{2}\}\rightarrow1,\quad\theta_2=\exp\{\hat{\lambda}_2-\hat{\lambda}_{3}\}\rightarrow\exp\{3.9-(1+\sqrt{3/5})^2\},\\
& \theta_3,\theta_4,\ldots,\theta_{p-1}\rightarrow1,\quad\tilde{\theta}_{p,2}=\frac{1}{p-3}\sum_{i=3}^{p-1}\theta_{i}^2\sim1,\quad\tilde{\theta}_{p,3}\sim1,\\
&\tilde{\theta}_{p,1}=\frac{1}{p-2}\sum_{i=2}^{p-1}\theta_{i}^2=\frac{\theta_2^2}{p-2}+\frac{1}{p-2}\sum_{i=3}^{p-1}\theta_i^2\sim\frac{\exp\{2[3.9-(1+\sqrt{3/5})^2]\}}{p-2}+\frac{p-3}{p-2}.
\end{split}
\end{equation}
Using (\ref{exeq1}) and (\ref{exeq2}), we have
\begin{equation}\label{exeq3}
\begin{split}
&\frac{1}{n}(\text{EDA}_1-\text{EDA}_2)\sim(p-2)\log\tilde{\theta}_{p,1}+(l_{n,1}-l_{n,3})-2\frac{p}{n}\approx2.94>0,\\
&\frac{1}{n}(\text{EDB}_1-\text{EDB}_2)\sim(p-2)\log\tilde{\theta}_{p,1}+\log(p)(l_{n,1}-l_{n,3})-(\log n)\frac{p}{n}\approx3.7>0,
\end{split}
\end{equation}
which means EDA and EDB can not lead to underestimation of the number of clusters, and the following expressions imply that they do not also lead to overestimation
\begin{equation}\label{exeq4}
\begin{split}
&\frac{1}{n}(\text{EDA}_3-\text{EDA}_2)\sim2\frac{p}{n}=1.2>0,\\
&\frac{1}{n}(\text{EDB}_3-\text{EDB}_2)\sim(\log n)\frac{p}{n}\approx2.76>0.
\end{split}
\end{equation}
From (\ref{exeq3}) and (\ref{exeq4}), it follows that both EDA and EDB are able to estimate the number of clusters accurately.

\end{rmk}
\subsection{Additional simulations}
\subsubsection{The cases of $c<1$}
In this part, we also include two more cases when $c<1$ as a supplement. Specifically, we consider the following two cases:

 \textbf{Case 3.} Let $\bmu_{1}=(5,0,-4,0,0,\ldots,0)^{\top}\in\mathbb{R}^{p}$, $\bmu_{2}=(0,4,0,-6,0,\ldots,0)^{\top}\in\mathbb{R}^{p}$, \\$\bmu_{3}=(0,-5,-5,0,0,\ldots,0)^{\top}\in\mathbb{R}^{p}$,
$\bmu_{4}=(-6,0,0,6,0,\ldots,0)^{\top}\in\mathbb{R}^{p}$, and $\bSigma=(\sigma_{i,j})_{p\times p}$, where $\sigma_{i,j}=0.2^{|i-j|}$. Define
\begin{equation*} \bA_{n}=\big(\underbrace{\bmu_{1},\ldots,\bmu_{1}}_{n_{1}},\underbrace{\bmu_{2},\ldots,\bmu_{2}}_{n_{2}},\underbrace{\bmu_{3},\ldots,\bmu_{3}}_{n_{3}},\underbrace{\bmu_{4},\ldots,\bmu_{4}}_{n_{4}}\big),
\end{equation*}
where $n_{1}=n_{3}=0.3n,~n_{2}=n_{4}=0.2n$. Therefore, the true number of clusters is $K=4$.

\textbf{Case 4.} The same setting as in the above Case 2 with $\bSigma=(\sigma_{i,j})_{p\times p}$ instead of $\bI$, where $\sigma_{i,j}=0.2^{|i-j|}$.

Tables \ref{table1} and \ref{table3} illustrate similar performance as those in Case 1 and 2, which suggest the robustness of the proposed method under different scenarios.

%%%%%%% Table 1 %%%%%%%%%%%
\begin{table}[h]\centering
\caption{Selection percentages of EDA, EDB, ASI, GS and BICdf in Case 3. Entries in the $\mathcal {F}_{*}$ rows indicate that EDA and EDB exhibit higher accuracy in estimating the number of clusters compared to other criteria.}\label{table1}
\begin{tabular}{c|c c   c c c c c | c c c c c}\hline
    &     &  & EDA & EDB & ASI & GS & BICdf & EDA & EDB & ASI & GS & BICdf\\\hline
$c$ & $n$ &  & \multicolumn{5}{c|}{$\mathcal {N}(0,1)$} & \multicolumn{5}{c}{$t_{8}$}\\\hline
\multirow{6}[6]*{$\frac{1}{3}$} & \multirow{3}[3]*{$180$} &
                                $\mathcal {F}_{-}$  & $0$ & $0$ & $69.1$ & $32.7$ & $1.1$  & $0$ & $0$ & $68.4$ & $31.4$ & $0.9$\\
                                &                         &
                                $\mathcal {F}_{*}$  & $59.8$ & $83.4$ & $30.7$ & $60.3$ & $67.7$ & $57.4$ & $78.5$ & $31$ & $61.4$ & $67.1$\\
                                &                         &
                                $\mathcal {F}_{+}$  & $40.2$ & $16.6$ & $0.2$ & $7$ & $31.2$ & $42.6$ & $21.5$ & $0.6$ & $7.2$ & $32$\\\cline{2-13}
                                & \multirow{3}[3]*{$450$} &
                                $\mathcal {F}_{-}$  & $0$ & $0$ & $75$ & $24.7$ & $1.6$ & $0$ & $0$ & $73.6$ & $26.1$ & $1.3$\\
                                &                         &
                                $\mathcal {F}_{*}$  & $93.1$ & $98.9$ & $24.8$ & $66.4$ & $71.9$ & $94.1$ & $99.2$ & $26.4$ & $64.6$ & $70.5$\\
                                &                         &
                                $\mathcal {F}_{+}$  & $6.9$  & $1.1$ & $0.2$ & $8.9$ & $26.5$ & $5.9$ & $0.8$ & $0$ & $9.3$ & $28.2$\\\hline
\multirow{6}[6]*{$\frac{1}{2}$} & \multirow{3}[3]*{$120$} &
                                $\mathcal {F}_{-}$  & $0.2$ & $0.9$ & $70.9$ & $32.5$ & $9.5$ & $0.8$ & $1.3$ & $70.7$ & $32.7$ & $8$\\
                                &                         &
                                $\mathcal {F}_{*}$  & $68.8$ & $83.6$ & $28.5$ & $61.1$ & $67.2$ & $67.4$ & $81.7$ & $28.6$ & $61.6$ & $66.1$\\
                                &                         &
                                $\mathcal {F}_{+}$  & $31$ & $15.5$ & $0.6$ & $6.4$ & $23.3$ & $31.8$ & $17$ & $0.7$ & $5.7$ & $25.9$\\\cline{2-13}
                                & \multirow{3}[3]*{$300$} &
                                $\mathcal {F}_{-}$  & $0$ & $0.3$ & $74.8$ & $25.7$ & $18.2$ & $0$ & $0.4$ & $72.2$ & $24.8$ & $15$\\
                                &                         &
                                $\mathcal {F}_{*}$  & $96$ & $99.1$ & $24.8$ & $67.4$ & $68.3$ & $97.4$ & $99.3$ & $27.6$ & $66.2$ & $69.7$\\
                                &                         &
                                $\mathcal {F}_{+}$  & $4$ & $0.6$ & $0.4$ & $6.9$ & $13.5$ & $2.6$ & $0.3$ & $0.2$ & $9$ & $15.3$\\\hline
\multirow{6}[6]*{$\frac{3}{4}$} & \multirow{3}[3]*{$80$}  &
                                $\mathcal {F}_{-}$  & $5.6$ & $10.9$ & $70.2$ & $34.7$ & $17.9$ & $6.6$ & $11.9$ & $70$ & $34.8$ & $19.3$\\
                                &                         &
                                $\mathcal {F}_{*}$  & $72.5$ & $77.9$ & $28.6$ & $59.8$ & $66.6$ & $68.7$ & $75$ & $28.1$ & $59.3$ & $65.6$\\
                                &                         &
                                $\mathcal {F}_{+}$  & $21.9$ & $11.2$ & $1.2$ & $5.5$ & $15.5$ & $24.7$ & $13.1$ & $1.9$ & $5.9$ & $15.1$\\\cline{2-13}
                                & \multirow{3}[3]*{$200$} &
                                $\mathcal {F}_{-}$  & $6.6$ & $15$ & $76.1$ & $29.6$ & $30.9$ & $9$ & $17.9$ & $75.1$ & $29$ & $26.6$\\
                                &                         &
                                $\mathcal {F}_{*}$  & $91$ & $84.6$ & $23.5$ & $64$ & $62.8$ & $88.2$ & $81.6$ & $24.5$ & $64.9$ & $66$\\
                                &                         &
                                $\mathcal {F}_{+}$  & $2.4$ & $0.4$ & $0.4$ & $6.4$ & $6.3$ & $2.8$ & $0.5$ & $0.4$ & $6.1$ & $7.4$\\\hline
$c$ & $n$ &  & \multicolumn{5}{c|}{Bernoulli} & \multicolumn{5}{c}{$\chi^{2}(3)$}\\\hline
\multirow{6}[6]*{$\frac{1}{3}$} & \multirow{3}[3]*{$180$} &
                                $\mathcal {F}_{-}$  & $0$ & $0$ & $71.9$ & $30.3$ & $1.1$ & $0$ & $0$ & $65.7$ & $29.3$ & $1.1$\\
                                &                         &
                                $\mathcal {F}_{*}$  & $64.2$ & $82.4$ & $27.5$ & $60.9$ & $65.1$ & $52.4$ & $75.2$ & $32.9$ & $62.2$ & $68.2$\\
                                &                         &
                                $\mathcal {F}_{+}$  & $35.8$ & $17.6$ & $0.6$ & $8.8$ & $33.8$ & $47.6$ & $24.8$ & $1.4$ & $8.5$ & $30.7$\\\cline{2-13}
                                & \multirow{3}[3]*{$450$} &
                                $\mathcal {F}_{-}$  & $0$ & $0$ & $75.7$ & $26.6$ & $2.3$ & $0$ & $0$ & $68.3$ & $27.1$ & $1.1$\\
                                &                         &
                                $\mathcal {F}_{*}$  & $96.6$ & $98.8$ & $24.2$ & $65.7$ & $66.9$ & $93.1$ & $99$ & $31.3$ & $64$ & $70$\\
                                &                         &
                                $\mathcal {F}_{+}$  & $3.4$ & $1.2$ & $0.1$ & $7.7$ & $30.8$ & $6.9$ & $1$ & $0.4$ & $8.9$ & $28.9$\\\hline
\multirow{6}[6]*{$\frac{1}{2}$} & \multirow{3}[3]*{$120$} &
                                $\mathcal {F}_{-}$  & $0.2$ & $0.5$ & $71.8$ & $31.3$ & $8.8$ & $0.7$ & $2.4$ & $69.5$ & $32.1$ & $8.1$\\
                                &                         &
                                $\mathcal {F}_{*}$  & $72.3$ & $85.2$ & $27.7$ & $62.2$ & $68.5$ & $61.5$ & $77$ & $29$ & $62$ & $68.5$\\
                                &                         &
                                $\mathcal {F}_{+}$  & $27.5$ & $14.3$ & $0.5$ & $6.5$ & $22.7$ & $37.8$ & $20.6$ & $1.5$ & $5.9$ & $23.4$\\\cline{2-13}
                                & \multirow{3}[3]*{$300$} &
                                $\mathcal {F}_{-}$  & $0.1$ & $0.1$ & $76$ & $26.7$ & $16.3$ & $0.1$ & $1.1$ & $69.5$ & $23.4$ & $13$\\
                                &                         &
                                $\mathcal {F}_{*}$  & $97.5$ & $99.4$ & $23.9$ & $64.3$ & $71.5$ & $96.6$ & $98.2$ & $30.2$ & $67.7$ & $70.1$\\
                                &                         &
                                $\mathcal {F}_{+}$  & $2.4$ & $0.5$ & $0.1$ & $9$ & $12.2$ & $3.3$ & $0.7$ & $0.3$ & $8.9$ & $16.9$\\\hline
\multirow{6}[6]*{$\frac{3}{4}$} & \multirow{3}[3]*{$80$}  &
                                $\mathcal {F}_{-}$  & $4.4$ & $6$ & $68.9$ & $33.8$ & $19.9$ & $7.7$ & $13.9$ & $67.7$ & $33.8$ & $17.7$\\
                                &                         &
                                $\mathcal {F}_{*}$  & $74.8$ & $83.1$ & $30.1$ & $62.4$ & $66.2$ & $64.2$ & $69.3$ & $29.5$ & $60.8$ & $67.7$\\
                                &                         &
                                $\mathcal {F}_{+}$  & $20.8$ & $10.9$ & $1$ & $3.8$ & $13.9$ & $28.1$ & $16.8$ & $2.8$ & $5.4$ & $14.6$\\\cline{2-13}
                                & \multirow{3}[3]*{$200$} &
                                $\mathcal {F}_{-}$  & $5.8$ & $12.7$ & $75.6$ & $28.7$ & $30.8$ & $10.5$ & $19.9$ & $71.4$ & $28.7$ & $28.8$\\
                                &                         &
                                $\mathcal {F}_{*}$  & $92.6$ & $87$ & $24.2$ & $65.4$ & $63.9$ & $86.4$ & $79.9$ & $28$ & $65.7$ & $63.4$\\
                                &                         &
                                $\mathcal {F}_{+}$  & $1.6$ & $0.3$ & $0.2$ & $5.9$ & $5.3$ & $3.1$ & $0.2$ & $0.6$ & $5.6$ & $7.8$\\\hline
\end{tabular}
\end{table}
%%%%%%% Table 3 %%%%%%%%%%%

\begin{table}[!h]\centering
\caption{Selection percentages of EDA, EDB, ASI, GS and BICdf in Case 4. Entries in the $\mathcal {F}_{*}$ rows indicate that EDA and EDB exhibit higher accuracy in estimating the number of clusters compared to other criteria.}\label{table3}
\begin{tabular}{c|c c   c c c c c | c c c c c}\hline
    &     &  & EDA & EDB & ASI & GS & BICdf & EDA & EDB & ASI & GS & BICdf\\\hline
$c$ & $n$ &  & \multicolumn{5}{c|}{$\mathcal {N}(0,1)$} & \multicolumn{5}{c}{$t_{8}$}\\\hline
\multirow{6}[6]*{$\frac{1}{3}$} & \multirow{3}[3]*{$180$} & $\mathcal {F}_{-}$  &
$0$ & $0$ & $5.4$ & $69.4$ & $92.3$ & $0$ & $0$ & $3.8$ & 82.4 & $89.4$\\
                                &                         & $\mathcal {F}_{*}$  &
                                $58.1$ & $81.2$ & $94.4$ & $30.6$ & $7.7$ & $58$ & $78.6$ & $95$ & $17.6$ & $10.6$\\
                                &                         & $\mathcal {F}_{+}$  &
                                $41.9$ & $18.8$ & $0.2$ & $0$ & $0$ & $42$ & $21.4$ & $1.2$ & $0$ & $0$\\\cline{2-13}
                                & \multirow{3}[3]*{$450$} & $\mathcal {F}_{-}$  &
                                $0$ & $0$ & $2.7$ & $60.9$ & $99.2$ & $0$ & $0$ & $3.5$ & $77.6$ & $99$\\
                                &                         & $\mathcal {F}_{*}$  &
                                $93.2$ & $99.1$ & $97.3$ & $39.1$ & $0.8$ & $93.7$ & $99.6$ & $96.3$ & $22.4$ & $1$\\
                                &                         & $\mathcal {F}_{+}$  &
                                $6.8$ & $0.9$ & $0$ & $0$ & $0$ & $6.3$ & $0.4$ & $0.2$ & $0$ & $0$\\\cline{2-13}
\multirow{6}[6]*{$\frac{1}{2}$} & \multirow{3}[3]*{$120$} & $\mathcal {F}_{-}$  &
$0.3$ & $0.4$ & $9.3$ & $98.6$ & $98.8$ & $0.1$ & $0.7$ & $11.9$ &  $99.8$& $99.3$\\
                                &                         & $\mathcal {F}_{*}$  &
                                $69.6$ & $82.2$ & $89.6$ & $1.4$ & $1.2$ & $68.1$ & $80.6$ & $85$ & $0.2$ & $0.7$\\
                                &                         & $\mathcal {F}_{+}$  &
                                $30.1$ & $17.4$ & $1.1$ & $0$ & $0$ & $31.8$ & $18.7$ & $3.1$ & $0$ & $0$\\\cline{2-13}
                                & \multirow{3}[3]*{$300$} & $\mathcal {F}_{-}$  & $0$ & $0$ & $4.6$ & $99.5$ & $100$ & $0$ & $0$ & $8.2$ & $100$ & $100$\\
                                &                         & $\mathcal {F}_{*}$  &
                                $95.9$ & $99.2$ & $95.3$ & $0.5$ & $0$ & $96.7$ & $99.6$ & $91.3$ & $0$ & $0$\\
                                &                         & $\mathcal {F}_{+}$  &
                                $4.1$ & $0.8$ & $0.1$ & $0$ & $0$ & $3.3$ & $0.4$ & $0.5$ & $0$ & $0$\\\hline
\multirow{6}[6]*{$\frac{3}{4}$} & \multirow{3}[3]*{$80$}  & $\mathcal {F}_{-}$  &
$3.5$ & $8.4$ & $17.2$ & $100$ & $100$ & $5.5$ & $10.4$ & $19.9$ & $100$ & $99.9$\\
                                &                         & $\mathcal {F}_{*}$  &
                                $72.7$ & $79.6$ & $79$ & $0$ & $0$ & $71$ & $73.3$ & $72.7$ & $0$ & $0.1$\\
                                &                         & $\mathcal {F}_{+}$  &
                                $23.8$ & $12$ & $3.8$ & $0$ & $0$ & $23.5$ & $16.3$ & $7.4$ & $0$ & $0$\\\cline{2-13}
                                & \multirow{3}[3]*{$200$} & $\mathcal {F}_{-}$  &
                                $1$ & $6.9$ & $11$ & $100$ & $100$ & $2.1$ & $8.4$ & $14.2$ & $100$ & $100$\\
                                &                         & $\mathcal {F}_{*}$  &
                                $96.2$ & $92.8$ & $88.2$ & $0$ & $0$ & $95$ & $91.3$ & $84.3$ & $0$ & $0$\\
                                &                         & $\mathcal {F}_{+}$
                                & $2.8$ & $0.3$ & $0.8$ & $0$ & $0$ & $2.9$ & $0.3$ & $1.5$ & $0$ & $0$\\\hline

$c$ & $n$ &  & \multicolumn{5}{c|}{Bernoulli} & \multicolumn{5}{c}{$\chi^{2}(3)$}\\\hline
\multirow{6}[6]*{$\frac{1}{3}$} & \multirow{3}[3]*{$180$} & $\mathcal {F}_{-}$  &
$0$ & $0$ & $2.1$ & $39.6$ & $90.2$ & $0$ & $0$ & $7$ & $90.2$ & $87.6$\\
                                &                         & $\mathcal {F}_{*}$  &
                                $61.2$ & $82.4$ & $97.8$ & $60.4$ & $9.8$ & $51.8$ & $75.1$ & $87.6$ & $9.8$ & $12.4$\\
                                &                         & $\mathcal {F}_{+}$  &
                                $38.8$ & $17.6$ & $0.1$ & $0$ & $0$ & $48.2$ & $24.9$ & $5.4$ & $0$ & $0$\\\cline{2-13}
                                & \multirow{3}[3]*{$450$} & $\mathcal {F}_{-}$  &
                                $0$ & $0$ & $0.5$ & $36.7$ & $99.3$ & $0$ & $0$ & $4.1$ & $89.3$ & $99$\\
                                &                         & $\mathcal {F}_{*}$  &
                                $93.9$ & $99.5$ & $99.5$ & $63.3$ & $0.7$ & $93.4$ & $99.2$ & $94.9$ & $10.7$ & $1$\\
                                &                         & $\mathcal {F}_{+}$  &
                                $6.1$ & $0.5$ & $0$ & $0$ & $0$ & $6.6$ & $0.8$ & $1$ & $0$ & $0$\\\hline
\multirow{6}[6]*{$\frac{1}{2}$} & \multirow{3}[3]*{$120$} & $\mathcal {F}_{-}$  &
$0.3$ & $0$ & $7.1$ & $95$ & $99.6$ & $0.2$ & $1.3$ & $15$ & $99.8$ & $98.2$\\
                                &                         & $\mathcal {F}_{*}$  &
                                $71.8$ & $87.5$ & $92.7$ & $5$ & $0.4$ & $59.3$ & $74.6$ & $75.5$ & $0.2$ & $1.8$\\
                                &                         & $\mathcal {F}_{+}$  &
                                $27.9$ & $12.5$ & $0.2$ & $0$ & $0$ & $40.5$ & $24.1$ & $9.5$ & $0$ & $0$\\\cline{2-13}
                                & \multirow{3}[3]*{$300$} & $\mathcal {F}_{-}$  &
                                $0$ & $0$ & $1.6$ & $96.9$ & $100$ & $0$ & $0.4$ & $10.4$ & $100$ & $100$\\
                                &                         & $\mathcal {F}_{*}$  &
                                $97.1$ & $99.8$ & $98.4$ & $3.1$ & $0$ & $95$ & $99$ & $87.7$ & $0$ & $0$\\
                                &                         & $\mathcal {F}_{+}$  &
                                $2.9$ & $0.2$ & $0$ & $0$ & $0$ & $5$ & $0.6$ & $1.9$ & $0$ & $0$\\\hline
\multirow{6}[6]*{$\frac{3}{4}$} & \multirow{3}[3]*{$80$} & $\mathcal {F}_{-}$   &
$3.2$ & $5.7$ & $18.3$ & $100$ & $100$ & $10.1$ & $15.1$ & $23.2$ & $100$ & $100$\\
                                &                         & $\mathcal {F}_{*}$  &
                                $77$ & $83.5$ & $79.1$ & $0$ & $0$ & $61.5$ & $69.6$ & $57.1$ & $0$ & $0$\\
                                &                         & $\mathcal {F}_{+}$  &
                                $19.8$ & $10.8$ & $2.6$ & $0$ & $0$ & $28.4$ & $15.3$ & $19.6$ & $0$ & $0$\\\cline{2-13}
                                & \multirow{3}[3]*{$200$} & $\mathcal {F}_{-}$  &
                                $0.6$ & $4.6$ & $5.9$ & $100$ & $100$ & $5.2$ & $14.6$ & $21.5$ & $100$ & $100$\\
                                &                         & $\mathcal {F}_{*}$  &
                                $97.3$ & $95.4$ & $93.9$ & $0$ & $0$ & $91.9$ & $85.2$ & $74.7$ & $0$ & $0$\\
                                &                         & $\mathcal {F}_{+}$  &
                                $2.1$ & $0$ & $0.2$ & $0$ & $0$ & $2.9$ & $0.2$ & 3.8 & $0$ & $0$\\\hline
\end{tabular}
\end{table}

\subsubsection{The cases of $c>1$}
We also consider the consistency properties of pseudo-EDA $\hat{K}_{\acute{\text{EDA}}}$ and pseudo-EDB $\hat{K}_{\acute{\text{EDB}}}$ when $c=3/2$ and $3$ under the following situations:

\textbf{Case 5.} Let $\bmu_{1}=(5,0,0,0,\ldots,0)^{\top}\in\mathbb{R}^{p}$, $\bmu_{2}=(0,6,0,0,\ldots,0)^{\top}\in\mathbb{R}^{p}$, $\bmu_{3}=(-2,0,4,0,\ldots,0)^{\top}\in\mathbb{R}^{p}$, $\bSigma=(\sigma_{i,j})_{p\times p}$, where $\sigma_{i,j}=0.2^{|i-j|}$. Then,
\begin{equation*} \bA_n=\big(\underbrace{\bmu_{1},\ldots,\bmu_{1}}_{n_{1}},\underbrace{\bmu_{2},\ldots,\bmu_{2}}_{n_{2}},\underbrace{\bmu_{3},\ldots,\bmu_{3}}_{n_{3}}\big),
\end{equation*}
where $n_{1}=n_{2}=0.3n,~n_{3}=0.4n$. Therefore, the true number of clusters is $K=3$.

\textbf{Case 6.} Let $\bmu_{1}=(4,0,0,0,\ldots,0)^{\top}\in\mathbb{R}^{p}$, $\bmu_{2}=(0,4,0,0,\ldots,0)^{\top}\in\mathbb{R}^{p}$, $\bmu_{3}=(0,0,4,0,\ldots,0)^{\top}\in\mathbb{R}^{p}$,
$\bSigma=\bI$. Then, $\bA_n$ has the same form as above with $n_{1}=n_{2}=0.3n,~n_{3}=0.4n$. Therefore, the true number of clusters is $K=3$.

\textbf{Case 7.} The same setting as in Case 6 with $\bI$ replaced by $\bSigma=(\sigma_{i,j})_{p\times p}$, where $\sigma_{i,j}=0.2^{|i-j|}$.

\textbf{Case 8.} The same setting as in Case 2 with $a=\sqrt{n/10}$ instead of $a=\sqrt{p/10}$.
Generality, as can be seen from Table \ref{table5}, when $c=3/2$, with $p$ increasing, $\hat{K}_{\acute{\text{EDA}}}$ and $\hat{K}_{\acute{\text{EDB}}}$ perform better, especially $\hat{K}_{\acute{\text{EDB}}}$. As $c$ increases (fixed $p$ and reducing $n$), from (\ref{equ3.5.5}), the gap conditions of EDA and EDB are not easy to satisfy. In particular, the gap condition of EDB is more strict than that of EDA when $n(>20)$  and $c$ are large. Therefore, the performance of pseudo-EDA is better than that of pseudo-EDB at $c=3$. Other tables are similarly.

%%%%%%%%%%%%%%%%%%%%%%%%%%%%%% c>1 %%%%%%%%%%%%%%%%%%%%%%%%%%%%%%%%%%%%%%
%%%%%%% Table 5 %%%%%%%%%%%
\begin{table}[!h]\centering
\caption{Selection percentages of EDA, EDB, ASI, GS and BICdf in Case 5. Entries in the $\mathcal {F}_{*}$ rows indicate that EDA and EDB exhibit higher accuracy in estimating the number of clusters compared to other criteria.}\label{table5}
\begin{tabular}{c|c c   c c c c c | c c c c c}\hline
    &     &  & $\acute{\text{EDA}}$ & $\acute{\text{EDB}}$ & ASI & GS & BICdf & $\acute{\text{EDA}}$ & $\acute{\text{EDB}}$ & ASI & GS & BICdf\\\hline
$c$ & $n$ &  & \multicolumn{5}{c|}{$\mathcal {N}(0,1)$} & \multicolumn{5}{c}{$t_{8}$}\\\hline
\multirow{6}[6]*{$\frac{3}{2}$} & \multirow{3}[3]*{$60$}  & $\mathcal {F}_{-}$  &
$0.5$ & $0.6$ & $5.6$ & $21.6$ & $94.5$ & $0.2$ & $0.3$ & $5.4$ & $28.3$ & $95.1$\\
                                &                         & $\mathcal {F}_{*}$  &
                                $89$ & $93.2$ & $91.9$ & $78.4$ & $5.5$ & $84.4$ & $90.1$ & $89.8$ & $71.7$ & $4.9$\\
                                &                         & $\mathcal {F}_{+}$  &
                                $10.5$ & $6.2$ & $2.5$ & $0$ & $0$ & $15.4$ & $9.6$ & $4.8$ & $0$ & $0$\\\cline{2-13}
                                & \multirow{3}[3]*{$300$} & $\mathcal {F}_{-}$  & $0$ & $0$ & $4.9$ & $6.3$ & $100$ & $0$ & $0$ & $5.2$ & $10.9$ & $100$\\
                                &                         & $\mathcal {F}_{*}$  &
                                $100$ & $100$ & $94.2$ & $93.7$ & $0$ & $100$ & $100$ & $93.2$ & $89.1$ & $0$\\
                                &                         & $\mathcal {F}_{+}$  & $0$ & $0$ & $0.9$ & $0$ & $0$ & $0$ & $0$ & $1.6$ & $0$ & $0$\\\hline
\multirow{6}[6]*{$3$}           & \multirow{3}[3]*{$30$}  & $\mathcal {F}_{-}$  &
$15$ & $15.5$ & $6$ & $100$ & $100$ & $13.9$ & $14.9$ & $7.4$ & $99.9$ & $100$\\
                                &                         & $\mathcal {F}_{*}$  &
                                $73.7$ & $74.4$ & $89.8$ & $0$ & $0$ & $71.8$ & $71.7$ & $86.6$ & $0.1$ & $0$\\
                                &                         & $\mathcal {F}_{+}$  &
                                $11.3$ & $10.1$ & $4.2$ & $0$ & $0$ & $14.3$ & $13.4$ & $6$ & $0$ & $0$\\\cline{2-13}
                                & \multirow{3}[3]*{$150$} & $\mathcal {F}_{-}$  &
                                $8.3$ & $19.2$ & $3.7$ & $98.9$ & $100$ & $8.8$ & $19$ & $7.4$ & $99.6$ & $100$\\
                                &                         & $\mathcal {F}_{*}$  &
                                $91.7$ & $80.8$ & $94.8$ & $1.1$ & $0$ & $91.2$ & $81$ & $90.3$ & $0.4$ & $0$\\
                                &                         & $\mathcal {F}_{+}$  & $0$ & $0$ & $1.5$ & $0$ & $0$ & $0$ & $0$ & $2.3$ & $0$ & $0$\\\hline

$c$ & $n$ &  & \multicolumn{5}{c|}{Bernoulli} & \multicolumn{5}{c}{$\chi^{2}(3)$}\\\hline
\multirow{6}[6]*{$\frac{3}{2}$} & \multirow{3}[3]*{$60$}  & $\mathcal {F}_{-}$  &
$0$ & $0$ & $5.7$ & $12.1$ & $95.2$ & $0.4$ & $1.3$ & $4.7$ & $36.4$ & $94.3$\\
                                &                         & $\mathcal {F}_{*}$  &
                                $89.4$ & $93.6$ & $92.8$ & $87.9$ & $4.8$ & $81.2$ & $85.8$ & $87.7$ & $63.4$ & $5.7$\\
                                &                         & $\mathcal {F}_{+}$  &
                                $10.6$ & $6.4$ & $1.5$ & $0$ & $0$ & $18.4$ & $12.9$ & $7.6$ & $0$ & $0$\\\cline{2-13}
                                & \multirow{3}[3]*{$300$} & $\mathcal {F}_{-}$  & $0$ & $0$ & $4.1$ & $3.6$ & $100$ & $0$ & $0$ & $7.2$ & $12$ & $99.9$\\
                                &                         & $\mathcal {F}_{*}$  &
                                $100$ & $100$ & $95.7$ & $96.4$ & $0$ & $100$ & $100$ & $90.5$ & $88$ & $0.1$\\
                                &                         & $\mathcal {F}_{+}$  & $0$ & $0$ & $0.2$ & $0$ & $0$ & $0$ & $0$ & $2.3$ & $0$ & $0$\\\hline
\multirow{6}[6]*{$3$}           & \multirow{3}[3]*{$30$} & $\mathcal {F}_{-}$   &
$12$ & $12.7$ & $6.3$ & $99.8$ & $100$ & $17.1$ & $17.9$ & $6.7$ & $100$ & $100$\\
                                &                         & $\mathcal {F}_{*}$  &
                                $79.1$ & $79.6$ & $91.1$ & $0.2$ & $0$ & $64.6$ & $65.4$ & $80.5$ & $0$ & $0$\\
                                &                         & $\mathcal {F}_{+}$  & $
                                8.9$ & $7.7$ & $2.6$ & $0$ & $0$ & $18.3$ & $16.7$ & $12.8$ & $0$ & $0$\\\cline{2-13}
                                & \multirow{3}[3]*{$150$} & $\mathcal {F}_{-}$  &
                                $7.6$ & $16$ & $3.6$ & $97.3$ & $100$ & $14.2$ & $25.5$ & $8.6$ & $99.9$ & $100$\\
                                &                         & $\mathcal {F}_{*}$  &
                                $92.4$ & $84$ & $96.3$ & $2.7$ & $0$ & $85.8$ & $74.5$ & $86.2$ & $0.1$ & $0$\\
                                &                         & $\mathcal {F}_{+}$  & $0$ & $0$ & $0.1$ & $0$ & $0$ & $0$ & $0$ & $5.2$ & $0$ & $0$\\\hline
\end{tabular}
\end{table}

%%%%%%% Table 6 %%%%%%%%%%%
\begin{table}[!h]\centering
\caption{Selection percentages of EDA, EDB, ASI, GS and BICdf in Case 6. Entries in the $\mathcal {F}_{*}$ rows indicate that EDA and EDB exhibit higher accuracy in estimating the number of clusters compared to other criteria.}\label{table6}
\begin{tabular}{c|c c   c c c c c | c c c c c}\hline
    &     &  & $\acute{\text{EDA}}$ & $\acute{\text{EDB}}$ & ASI & GS & BICdf & $\acute{\text{EDA}}$ & $\acute{\text{EDB}}$ & ASI & GS & BICdf\\\hline
$c$ & $n$ &  & \multicolumn{5}{c|}{$\mathcal {N}(0,1)$} & \multicolumn{5}{c}{$t_{8}$}\\\hline
\multirow{6}[6]*{$\frac{3}{2}$} & \multirow{3}[3]*{$60$}  & $\mathcal {F}_{-}$  & $0$ & $0$ & $6.8$ & $100$ & $100$ & $0$ & $0$ & $9.6$ & $100$ & $100$\\
                                &                         & $\mathcal {F}_{*}$  & $93.3$ & $96.7$ & $91.4$ & $0$ & $0$ & $92.1$ & $95.2$ & $85.7$ & $0$ & $0$ \\
                                &                         & $\mathcal {F}_{+}$  &
                                $6.7$ & $3.3$ & $1.8$ & $0$ & $0$ & $7.9$ & $4.8$ & $4.7$ & $0$ & $0$\\\cline{2-13}
                                & \multirow{3}[3]*{$300$} & $\mathcal {F}_{-}$  & $0$ & $0$ & $2.9$ & $100$ & $100$ & $0$ & $0$ & $2.6$ & $100$ & $100$\\
                                &                         & $\mathcal {F}_{*}$  & $100$ & $100$ & $96.9$ & $0$ & $0$ & $99.8$ & $100$ & $96.6$ & $0$ & $0$\\
                                &                         & $\mathcal {F}_{+}$  & $0$ & $0$ & $0.2$ & $0$ & $0$ & $0.2$ & $0$ & $0.8$ & $0$ & $0$\\\hline
\multirow{6}[6]*{$3$}           & \multirow{3}[3]*{$30$}  & $\mathcal {F}_{-}$  &
$16.8$ & $18.6$ & $19.9$ & $100$ & $100$ & $21.1$ & $25.4$ & $22.4$ & $100$ & $100$\\
                                &                         & $\mathcal {F}_{*}$  & $78.6$ & $75.4$ & $70.2$ & $0$ & $0$ & $69.7$ & $68.3$ & $60.8$ & $0$ & $0$ \\
                                &                         & $\mathcal {F}_{+}$  & $4.6$ & $6$ & $9.9$ & $0$ & $0$ & $9.2$ & $6.3$ & $16.8$ & $0$ & $0$ \\\cline{2-13}
                                & \multirow{3}[3]*{$150$} & $\mathcal {F}_{-}$  & $1.2$ & $10$ & $9.6$ & $100$ & $100$ & $2.8$ & $13.4$ & $13.2$ & $100$ & $100$ \\
                                &                         & $\mathcal {F}_{*}$  & $98.8$ & $90$ & $89.4$ & $0$ & $0$ & $97.2$ & $86.6$ & $83.4$ & $0$ & $0$\\
                                &                         & $\mathcal {F}_{+}$  & $0$ & $0$ & $1$ & $0$ & $0$ & $0$ & $0$ & $3.4$ & $0$ & $0$\\\hline

$c$ & $n$ &  & \multicolumn{5}{c|}{Bernoulli} & \multicolumn{5}{c}{$\chi^{2}(3)$}\\\hline
\multirow{6}[6]*{$\frac{3}{2}$} & \multirow{3}[3]*{$60$}  & $\mathcal {F}_{-}$  & $0$ & $0$ & $6.1$ & $100$ & $100$ & $0.1$ & $0.3$ & $10.4$ & $100$ & $100$\\
                                &                         & $\mathcal {F}_{*}$  & $94.9$ & $97.7$ & $93.3$ & $0$ & $0$ & $88.1$ & $92.5$ & $74.1$ & $0$ & $0$ \\
                                &                         & $\mathcal {F}_{+}$  &
                                $5.1$ & $2.3$ & $0.6$ & $0$ & $0$ & $11.8$ & $7.2$ & $15.5$ & $0$ & $0$\\\cline{2-13}
                                & \multirow{3}[3]*{$300$} & $\mathcal {F}_{-}$  & $0$ & $0$ & $1.1$ & $99.9$ & $100$ & $0$ & $0$ & $8.6$ & $100$ & $100$\\
                                &                         & $\mathcal {F}_{*}$  & $100$ & $100$ & $98.8$ & $0.1$ & $0$ & $100$ & $100$ & $89.6$ & $0$ & $0$\\
                                &                         & $\mathcal {F}_{+}$  & $0$ & $0$ & $0.1$ & $0$ & $0$ & $0$ & $0$ & $1.8$ & $0$ & $0$\\\hline
\multirow{6}[6]*{$3$}           & \multirow{3}[3]*{$30$} & $\mathcal {F}_{-}$   & $12.2$ & $14.6$ & $19$ & $100$ & $100$ & $30.5$ & $30.7$ & $24.6$ & $100$ & $100$\\
                                &                         & $\mathcal {F}_{*}$  & $83.3$ & $82.4$ & $75.4$ & $0$ & $0$ & $60.1$ & $61.9$ & $46$ & $0$ & $0$\\
                                &                         & $\mathcal {F}_{+}$  & $4.5$ & $3$ & $5.6$ & $0$ & $0$ & $9.4$ & $7.4$ & $29.4$ & $0$ & $0$ \\\cline{2-13}
                                & \multirow{3}[3]*{$150$} & $\mathcal {F}_{-}$  & $1.4$ & $8.9$ & $2.1$ & $100$ & $100$ & $4.9$ & $19.7$ & $19.7$ & $100$ & $100$\\
                                &                         & $\mathcal {F}_{*}$  & $98.6$ & $91.1$ & $97.6$ & $0$ & $0$ & $95.1$ & $80.3$ & $74.3$ & $0$ & $0$ \\
                                &                         & $\mathcal {F}_{+}$  & $0$ & $0$ & $0.3$ & $0$ & $0$ & $0$ & $0$ & $6$ & $0$ & $0$\\\hline
\end{tabular}
\end{table}
\clearpage

%%%%%%% Table 7 %%%%%%%%%%%
\begin{table}[!h]\centering
\caption{Selection percentages of EDA, EDB, ASI, GS and BICdf in Case 7}. Entries in the $\mathcal {F}_{*}$ rows indicate that EDA and EDB exhibit higher accuracy in estimating the number of clusters compared to other criteria.\label{table7}
\begin{tabular}{c|c c   c c c c c | c c c c c}\hline
    &     &  & $\acute{\text{EDA}}$ & $\acute{\text{EDB}}$ & ASI & GS & BICdf & $\acute{\text{EDA}}$ & $\acute{\text{EDB}}$ & ASI & GS & BICdf\\\hline
$c$ & $n$ &  & \multicolumn{5}{c|}{$\mathcal {N}(0,1)$} & \multicolumn{5}{c}{$t_{8}$}\\\hline
\multirow{6}[6]*{$\frac{3}{2}$} & \multirow{3}[3]*{$60$}  & $\mathcal {F}_{-}$  & $0$ & $0$ & $5.7$ & $100$ & $100$ & $0$ & $0.1$ & $7.8$ & $100$ & $100$\\
                                &                         & $\mathcal {F}_{*}$  & $89.5$ & $93.5$ & $91.7$ & $0$ & $0$ & $86$ & $90.8$ & $86.8$ & $0$ & $0$\\
                                &                         & $\mathcal {F}_{+}$  & $10.5$ & $6.5$ & $2.6$ & $0$ & $0$ & $14$ & $9.1$ & $5.4$ & $0$ & $0$ \\\cline{2-13}
                                & \multirow{3}[3]*{$300$} & $\mathcal {F}_{-}$  & $0$ & $0$ & $2$ & $100$ & $100$ & $0$ & $0$ & $2.6$ & $100$ & $100$\\
                                &                         & $\mathcal {F}_{*}$  & $100$ & $100$ & $97.7$ & $0$ & $0$ & $99.9$ & $100$ & $96.3$ & $0$ & $0$\\
                                &                         & $\mathcal {F}_{+}$  & $0$ & $0$ & $0.3$ & $0$ & $0$ & $0.1$ & $0$ & $1.1$ & $0$ & $0$\\\hline
\multirow{6}[6]*{$3$}           & \multirow{3}[3]*{$30$}  & $\mathcal {F}_{-}$  & $24.5$ & $26.4$ & $17.8$ & $100$ & $100$ & $24.4$ & $26.6$ & $19.9$ & $100$ & $100$\\
                                &                         & $\mathcal {F}_{*}$  & $67.3$ & $66.7$ & $65.7$ & $0$ & $0$ & $64.9$ & $64.2$ & $61.2$ & $0$ & $0$ \\
                                &                         & $\mathcal {F}_{+}$  & $8.2$ & $6.9$ & $16.5$ & $0$ & $0$ & $10.7$ & $9.2$ & $18.9$ & $0$ & $0$ \\\cline{2-13}
                                & \multirow{3}[3]*{$150$} & $\mathcal {F}_{-}$  & $9.4$ & $34$ & $9.8$ & $100$ & $100$ & $10.2$ & $38$ & $14.7$ & $100$ & $100$ \\
                                &                         & $\mathcal {F}_{*}$  & $90.6$ & $66$ & $88.3$ & $0$ & $0$ & $89.8$ & $62$ & $81.5$ & $0$ & $0$\\
                                &                         & $\mathcal {F}_{+}$  & $0$ & $0$ & $1.9$ & $0$ & $0$ & $0$ & $0$ & $3.8$ & $0$ & $0$\\\hline

$c$ & $n$ &  & \multicolumn{5}{c|}{Bernoulli} & \multicolumn{5}{c}{$\chi^{2}(3)$}\\\hline
\multirow{6}[6]*{$\frac{3}{2}$} & \multirow{3}[3]*{$60$}  & $\mathcal {F}_{-}$  & $0.1$ & $0.2$ & $5.1$ & $99.7$ & $100$ & $0.3$ & $0.6$ & $9.1$ & $100$ & $100$ \\
                                &                         & $\mathcal {F}_{*}$  & $89.4$ & $93.3$ & $93.9$ & $0.3$ & $0$ & $82.6$ & $87.4$ & $76.8$ & $0$ & $0$\\
                                &                         & $\mathcal {F}_{+}$  & $10.5$ & $6.5$ & $1$ & $0$ & $0$ & $17.1$ & $12$ & $14.1$ & $0$ & $0$ \\\cline{2-13}
                                & \multirow{3}[3]*{$300$} & $\mathcal {F}_{-}$  & $0$ & $0$ & $1.4$ & $100$ & $100$ & $0$ & $0$ & $7.9$ & $100$ & $100$\\
                                &                         & $\mathcal {F}_{*}$  & $100$ & $100$ & $98.5$ & $0$ & $0$ & $100$ & $100$ & $90.1$ & $0$ & $0$\\
                                &                         & $\mathcal {F}_{+}$  & $0$ & $0$ & $0.1$ & $0$ & $0$ & $0$ & $0$ & $2$ & $0$ & $0$\\\hline
\multirow{6}[6]*{$3$}           & \multirow{3}[3]*{$30$} & $\mathcal {F}_{-}$   & $20.7$ & $21.9$ & $16.2$ & $100$ & $100$ & $30.9$ & $33.1$ & $25.8$ & $100$ & $100$\\
                                &                         & $\mathcal {F}_{*}$  & $71.5$ & $71.1$ & $73.7$ & $0$ & $0$ & $57.4$ & $56.2$ & $42.4$ & $0$ & $0$ \\
                                &                         & $\mathcal {F}_{+}$  & $7.8$ & $7$ & $10.1$ & $0$ & $0$ & $11.7$ & $10.7$ & $31.8$ & $0$ & $0$ \\\cline{2-13}
                                & \multirow{3}[3]*{$150$} & $\mathcal {F}_{-}$  & $6.9$ & $29.9$ & $2.4$ & $100$ & $100$ & $14.4$ & $45.7$ & $18.6$ & $100$ & $100$\\
                                &                         & $\mathcal {F}_{*}$  & $93.1$ & $70.1$ & $97$ & $0$ & $0$ & $85.6$ & $54.3$ & $73.5$ & $0$ & $0$\\
                                &                         & $\mathcal {F}_{+}$  & $0$ & $0$ & $0.6$ & $0$ & $0$ & $0$ & $0$ & $7.9$ & $0$ & $0$\\\hline
\end{tabular}
\end{table}

%%%%%%% Table 8 %%%%%%%%%%%
\begin{table}[!h]\centering
\caption{Selection percentages of EDA, EDB, ASI, GS and BICdf in Case 8. Entries in the $\mathcal {F}_{*}$ rows indicate that EDA and EDB exhibit higher accuracy in estimating the number of clusters compared to other criteria.}\label{table8}
\begin{tabular}{c|c c   c c c c c | c c c c c}\hline
    &     &  & $\acute{\text{EDA}}$ & $\acute{\text{EDB}}$ & ASI & GS & BICdf & $\acute{\text{EDA}}$ & $\acute{\text{EDB}}$ & ASI & GS & BICdf\\\hline
$c$ & $n$ &  & \multicolumn{5}{c|}{$\mathcal {N}(0,1)$} & \multicolumn{5}{c}{$t_{8}$}\\\hline
\multirow{6}[6]*{$\frac{3}{2}$} & \multirow{3}[3]*{$60$}  & $\mathcal {F}_{-}$  & $0$ & $0$ & $76.2$ & $45.6$ & $76.1$ & $0$ & $0$ & $75.2$ & $45.2$ & $76.7$\\
                                &                         & $\mathcal {F}_{*}$  & $88.9$ & $94$ & $22$ & $51.1$ & $23.9$ & $84.8$ & $90$ & $23.4$ & $52.3$ & 23.3\\
                                &                         & $\mathcal {F}_{+}$  & $11.1$ & $6$ & $1.1$ & $3.3$ & $0$ & $15.2$ & $10$ & $1.4$ & $2.5$ & $0$ \\\cline{2-13}
                                & \multirow{3}[3]*{$300$} & $\mathcal {F}_{-}$  & $0$ & $0$ & $72.4$ & $30.4$ & $1.8$ & $0$ & $0$ & $68.9$ & $30.9$ & $2.2$\\
                                &                         & $\mathcal {F}_{*}$  & $100$ & $100$ & $23.7$ & $58.9$ & $59.2$ & $99.9$ & $100$ & $27.8$ & $59.6$ & $60.7$\\
                                &                         & $\mathcal {F}_{+}$  & $0$ & $0$ & $3.9$ & $10.7$ & $39$ & $0.1$ & $0$ & $3.3$ & $9.5$ & $37.1$\\\hline
\multirow{6}[6]*{$3$}           & \multirow{3}[3]*{$30$} & $\mathcal {F}_{-}$   & $16.9$ & $19$ & $85.1$ & $100$ & $100$ & $20.1$ & $21.9$ & $82.6$ & $100$ & $100$\\
                                &                         & $\mathcal {F}_{*}$  & $73.3$ & $73.1$ & $13.4$ & $0$ & $0$ & $68.8$ & $67.6$ & $14.4$ & $0$ & $0$\\
                                &                         & $\mathcal {F}_{+}$  & $9.8$ & $7.9$ & $1.5$ & $0$ & $0$ & $11.1$ & $10.5$ & $3$ & $0$ & $0$\\\cline{2-13}
                                & \multirow{3}[3]*{$150$} & $\mathcal {F}_{-}$  & $0$ & $0$ & $86$ & $27.6$ & $84.3$ & $0$ & $0$ & $84.4$ & $28.1$ & $81.8$\\
                                &                         & $\mathcal {F}_{*}$  & $100$ & $100$ & $12.3$ & $63.1$ & $15.7$ & $99.8$ & $100$ & $13.2$ & $64.7$ & $18.2$\\
                                &                         & $\mathcal {F}_{+}$  & $0$ & $0$ & $1.7$ & $9.3$ & $0$ & $0.2$ & $0$ & $2.4$ & $7.2$ & $0$\\\hline

$c$ & $n$ &  & \multicolumn{5}{c|}{Bernoulli} & \multicolumn{5}{c}{$\chi^{2}(3)$}\\\hline
\multirow{6}[6]*{$\frac{3}{2}$} & \multirow{3}[3]*{$60$}  & $\mathcal {F}_{-}$  & $0$ & $0$ & $78.9$ & $38.3$ & $79.2$ & $0$ & $0$ & $72.9$ & $48.6$ & $78$\\
                                &                         & $\mathcal {F}_{*}$  &
                                $89.8$ & $93.2$ & $20.3$ & $57.7$ & $20.8$ & $78.6$ & $86.9$ & $24.5$ & $47.9$ & $21.8$\\
                                &                         & $\mathcal {F}_{+}$  &
                                $10.2$ & $6.8$ & $0.8$ & $4$ & $0$ & $21.4$ & $13.1$ & $2.6$ & $3.5$ & $0.2$\\\cline{2-13}
                                & \multirow{3}[3]*{$300$} & $\mathcal {F}_{-}$  &
                                $0$ & $0$ & $71.3$ & $32.3$ & $1.9$ & $0$ & $0$ & $69.9$ & $32.7$ & $1.8$\\
                                &                         & $\mathcal {F}_{*}$  &
                                $100$ & $100$ & $26.1$ & $58.2$ & $59.7$ & $100$ & $100$ & $26.8$ & $56.1$ & $60.1$\\
                                &                         & $\mathcal {F}_{+}$  &
                                $0$ & $0$ & $2.6$ & $9.5$ & $38.4$ & $0$ & $0$ & $3.3$ & $11.2$ & $38.1$\\\hline
\multirow{6}[6]*{$3$}           & \multirow{3}[3]*{$30$} & $\mathcal {F}_{-}$   &
$15.7$ & $14.5$ & $85.3$ & $100$ & $100$ & $24.2$ & $21.3$ & $78.4$ & $100$ & $100$\\
                                &                         & $\mathcal {F}_{*}$  &
                                $77.4$ & $78.4$ & $14.3$ & $0$ & $0$ & $58.6$ & $64.9$ & $16.9$ & $0$ & $0$\\
                                &                         & $\mathcal {F}_{+}$  &
                                $6.9$ & $7.1$ & $0.4$ & $0$ & $0$ & $17.2$ & $13.8$ & $4.7$ & $0$ & $0$\\\cline{2-13}
                                & \multirow{3}[3]*{$150$} & $\mathcal {F}_{-}$  &
                                $0$ & $0$ & $84.6$ & $25$ & $84$ & $0$ & $0$ & $81.8$ & $31.3$ & $82.1$\\
                                &                         & $\mathcal {F}_{*}$  &
                                $99.9$ & $100$ & $12.4$ & $65$ & $16$ & $100$ & $100$ & $15.1$ & $59.6$ & $17.9$\\
                                &                         & $\mathcal {F}_{+}$  &
                                $0.1$ & $0$ & $3$ & $10$ & $0$ & $0$ & $0$ & $3.1$ & $9.1$ & $0$\\\hline
\end{tabular}
\end{table}

\subsection{Proof of results in Section \ref{app}}
This section includes the proofs of Lemma \ref{lem2.4}, Theorems \ref{thm3.1}, \ref{thm3.2}, and \ref{thm3.3}.

\vspace{1em}
\noindent\textbf{Proof of Lemma \ref{lem2.4}.}
For any matrix $\bm{A}$, denote by $\sigma_{i}(\bm{A})$, $\rho_{i}(\bm{A})$ the $i$-th largest eigenvalue and singular value of $\bm{A}$, respectively. From conditions in Theorem \ref{thm2} and the main result of \citet{Yin1988}, it is shown that, with probability $1$, as $n\rightarrow\infty$, for $k=1,\ldots,K$, there is a constant $C$ such that
\begin{eqnarray}
\left|\hat{\lambda}_{k}-\gamma_{k}\right|&=&\left|\sigma_{k}( \bX_n\bX_n^{*})-\sigma_{k}(  \bA_n \bA_n^{*}+\bSigma)\right|\notag\\
&\leq&\left|\sigma_{k}( \bX_n\bX_n^{*})-\sigma_{k}(  \bA_n \bA_n^{*})\right|+\left|\sigma_{1}(\bSigma)\right|\notag\\
&=&\left|\rho_{k}^{2}( \bX_n)-\rho^{2}_{k}(  \bA_n)\right|+\left|\sigma_{1}(\bSigma)\right|\notag\\
&\le&\left|\rho_{k}( \bX_n)+\rho_{k}(  \bA_n)\right|\left|\rho_{1}( \bSigma^{1/2}\bW_{n})\right|+\left|\sigma_{1}(\bSigma)\right|\notag\\
&\le &2\sqrt{C}(1+\sqrt{c})\left|\rho_{k}( \bX_n)+\rho_{k}(  \bA_n)\right|+C.\label{equ5.0}
\end{eqnarray}
Since $\gamma_{K}\rightarrow\infty$, it follows that
\begin{equation}\label{equ5.1}
\frac{\left|\rho_{k}( \bX_n)+\rho_{k}(  \bA_n)\right|}{\left|\sigma_{k}(  \bA_n \bA_n^{*}+\bSigma)\right|}\leq\frac{2\left|\rho_{k}(  \bA_n)\right|+\rho_{k}( \bSigma^{1/2}\bm{W}_{n})}{\left|\sigma_{k}(  \bA_n \bA_n^{*})\right|}\leq\frac{C}{\left|\rho_{k}(  \bA_n)\right|}\rightarrow0.
\end{equation}
Dividing by $\gamma_{k}$ on the both sides of (\ref{equ5.0}), due to (\ref{equ5.1}), we complete the proof.
\qed
\\

\noindent\textbf{Proof of Theorem \ref{thm3.1}.}
We first consider the case where $k<K$. 
Note that the criteria in (\ref{crrr}) can be also expressed as
\begin{equation}\label{cr1}
\begin{split}
& \text{EDA}_{k}=-n\log(\theta_{1}\cdots\theta_{k})+n(p-k-1)\log\widetilde{\theta}_{p,k}+2pk,\\
& \text{EDB}_{k}=-n\log(p)\cdot\log(\theta_{1}\cdots\theta_{k})+n(p-k-1)\log\widetilde{\theta}_{p,k}+(\log n)pk.
\end{split}
\end{equation}
From (\ref{cr1}), write
%±ÈÖµ·½·¨µÄÍÆµ¼
%\begin{eqnarray}
%&&\frac{1}{n}\left(\text{EDA}_{k}(m)-\text{EDA}_{K}(m)\right)=\frac{1}{n}\sum_{i=k+1}^{K}\left(\text{EDA}_{i-1}(m)-\text{EDA}_{i}(m)\right)\notag\\
%&=&\sum_{i=k+1}^{K}\big[-m\log(\alpha_{1}\cdots\alpha_{i-1})+(p-i)\log\overline{\alpha}_{i-1}(m)+2(i-1)p/n-m\log\mathbb{P}(i-1)\notag\\
%&&~~~~~~~~~+m\log(\alpha_{1}\cdots\alpha_{i})-(p-i-1)\log\overline{\alpha}_{i}(m)-2ip/n+m\log\mathbb{P}(i)\big]\notag\\
%&=&\sum_{i=k+1}^{K}\left[m\log\alpha_{i}+(p-i)\log\frac{\overline{\alpha}_{i-1}(m)}{\overline{\alpha}_{i}(m)}+\log\overline{\alpha}_{i}(m)-2\frac{p}{n}+m\log\frac{\mathbb{P}(i)}{\mathbb{P}(i-1)}\right]\notag\\
%&=&\sum_{i=k+1}^{K}\left\{m\log\alpha_{i}+(p-i)\log\left[1-\frac{1}{p-i}\left(1-\frac{\alpha_{i}^{m}}{\overline{\alpha}_{i}(m)}\right)\right]+\log\overline{\alpha}_{i}(m)-2\frac{p}{n}\right\}\notag\\
%&\sim&\sum_{i=k+1}^{K}\left\{\frac{\beta_{i}^{m}}{\overline{\alpha}_{i}(m)}+\log\beta_{i}^m+\log\overline{\alpha}_{i}(m)-1-2c\right\}.\label{equ3.7}
%\end{eqnarray}
\begin{eqnarray}
&&\frac{1}{n}\left(\text{EDA}_{k}-\text{EDA}_{K}\right)=\frac{1}{n}\sum_{i=k+1}^{K}\left(\text{EDA}_{i-1}-\text{EDA}_{i}\right)\notag\\
&=&\sum_{i=k+1}^{K}\left\{\log\theta_{i}+(p-i)\log\left[1-\frac{1}{p-i}\left(1-\frac{\theta_{i}^{2}}{\widetilde{\theta}_{p,i}}\right)\right]+\log\widetilde{\theta}_{p,i}-2\frac{p}{n}\right\}\notag\\
&\sim&\sum_{i=k+1}^{K}\left\{\frac{\zeta_{i}^2}{\widetilde{\theta}_{p,i}}+\log\zeta_{i}+\log\widetilde{\theta}_{p,i}-1-2c\right\}.\label{equ3.7}
\end{eqnarray}
If there are $h=o(p)$ bulks in the support of $F^{c,H}$, and let $\mu_{r_j}$ be the left bound of the    %$\Gamma_{F^{c,H}}$, 
from (\ref{equ2.5new}), we have
\begin{equation}\label{equ3.5.4}
\begin{split}
& \theta_{r_{j}}=O(1)~~\text{a.s.},~r_{j}\in\{K+1,\ldots,p-1\},j=1,\ldots,h-1,\\
& \theta_{r}\rightarrow1~~\text{a.s.},~ r\in\mathbb{L}\triangleq\{K+1,\ldots,p-1\}\setminus\{r_{1},\ldots,r_{h-1}\}.
\end{split}
\end{equation}
Combining it with (\ref{equ3.3}), for $i\in[k,K]$, it yields
\begin{equation*}
\begin{split}
1\leq\widetilde{\theta}_{p,i}=&\frac{1}{p-i-1}\left(\theta_{i+1}^2+\cdots+\theta_{K}^2+\theta_{K+1}^2+\cdots+\theta_{p-1}^2\right)\\
\leq&\frac{1}{p-i-1}\left((K-i+h-1)\max\limits_{j\in\{i+1,\ldots, K,r_{1},\ldots,r_{h-1}\}}\theta_{j}^2+\sum_{r\in\mathbb{L}}\theta_{r}^2\right)\rightarrow1,
\end{split}
\end{equation*}
as $p\rightarrow\infty$. Thus, (\ref{equ3.7}) is equivalent to
%±ÈÖµ·½·¨µÄµÈ¼ÛÊ½
%\begin{equation}\label{equ3.7.1}
%\sum_{i=k+1}^{K}\left\{\beta_{i}^{m}+\log\beta_{i}^m-1-2c\right\}.
%\end{equation}
\begin{equation}\label{equ3.7.1}
\sum_{i=k+1}^{K}\left\{\zeta_{i}^{2}+\log\zeta_{i}-1-2c\right\}.
\end{equation}
If the gap condition (\ref{equ3.6}) does not hold, (\ref{equ3.7.1}) can be negative, so that $\hat{K}_{\text{EDA}}$ is not consistent. Otherwise, for $k<K$ and sufficiently large $p$, we have $\frac{1}{n}\left(\text{EDA}_{k}-\text{EDA}_{K}\right)>0.$
%\begin{equation*}
%\frac{1}{n}\left(\text{EDA}_{k}-%\text{EDA}_{K}\right)>0.
%\end{equation*}
In other words,
\begin{equation}\label{equ3.8}
\hat{K}_{\text{EDA}}=\arg\min\limits_{k=1,\ldots,K}\frac{1}{n}\text{EDA}_{k}=K,~a.s..
\end{equation}
Next, consider the case that $K<k\leq w$. It follows that
\begin{eqnarray}
&&\frac{1}{n}\left(\text{EDA}_{k}-\text{EDA}_{K}\right)=\frac{1}{n}\sum_{i=K+1}^{k}\left(\text{EDA}_{i}-\text{EDA}_{i-1}\right)\notag\\
&=&\sum_{i=K+1}^{k}\left\{-\log\theta_{i}-(p-i)\log\left[1-\frac{1}{p-i}\left(1-\frac{\theta_{i}^{2}}{\widetilde{\theta}_{p,i}}\right)\right]-\log\widetilde{\theta}_{p,i}+2\frac{p}{n}\right\}\notag\\
&\sim&\sum_{i=K+1}^{k}\left\{1-\frac{\theta_{i}^{2}}{\widetilde{\theta}_{p,i}}-\log\theta_{i}-\log\widetilde{\theta}_{p,i}+2c\right\}.\label{equ3.9}
\end{eqnarray}
By (\ref{equ3.5.4}), for $i=K+1,\ldots,w$, we have
%±ÈÖµ·½·¨µÄ½¥½ü½á¹û
%\begin{equation*}
%\alpha_{i}\rightarrow1~~\text{and}~~\overline{\alpha}_{i}(m)=\frac{1}{p-i-1}\sum_{j=i+1}^{p-1}(\alpha_{j})^{m}\rightarrow1,~a.s..
%\end{equation*}
\begin{equation*}
\widetilde{\theta}_{p,i}=\frac{1}{p-i-1}\sum_{j=i+1}^{p-1}\theta_{j}^2\rightarrow1~~~\text{a.s.}.
\end{equation*}
Hence, (\ref{equ3.9}) is equivalent to $2(k-K)c>0$, which follows from $w=o(p)$. Then,
\begin{equation*}
\hat{K}_{\text{EDA}}=\arg\min\limits_{k=K,\ldots,w}\frac{1}{n}\text{EDA}_{k}=K~~~~\text{a.s.},
\end{equation*}
from which with (\ref{equ3.8}) conclusion (i) follows.

 If $\gamma_{K}\rightarrow\infty$, note that the proof for the case where $K<k\leq w$ proceeds in the same manner as before, which will not be repeated here. 
 
 For $k<K$, from Lemma \ref{lem2.4} and the second assertion in Theorem \ref{thm2},  it yields
\begin{eqnarray}
& &\frac{1}{n}\left(\text{EDA}_{k}-\text{EDA}_{K}\right)\notag\\
&=&\log(\theta_{k+1}\cdots\theta_{K})+(p-K-1)\log\frac{\widetilde{\theta}_{p,k}}{\widetilde{\theta}_{p,K}}+(K-k)\log\widetilde{\theta}_{p,k}-2(K-k)\frac{p}{n}\notag\\
&=&\hat{\lambda}_{k+1}-\lambda_{K+1}+(p-K-1)\log\left[1+\frac{1}{p-k-1}\left(\frac{\theta_{k+1}^2+\cdots+\theta_{K}^2}{\widetilde{\theta}_{p,K}}-(K-k)\right)\right]\notag\\
& &+(K-k)\log\widetilde{\theta}_{p,k}-2(K-k)\frac{p}{n}\notag\\
&\sim&\gamma_{k+1}-\mu_{K+1}+(p-K-1)\log\left[1+\frac{1}{p-k-1}\left(\frac{\theta_{k+1}^2+\cdots+\theta_{K}^2}{\widetilde{\theta}_{p,K}}-(K-k)\right)\right]\notag\\
& &+(K-k)\log\widetilde{\theta}_{p,k}-2(K-k)c\notag\\
&\geq&\gamma_{k+1}-\mu_{K+1}+(p-K-1)\log\left[1+\frac{1}{p-k-1}\left(\frac{\theta_{K}^2}{\widetilde{\theta}_{p,K}}-(K-k)\right)\right]\notag\\
& &+(K-k)\log\widetilde{\theta}_{p,k}-2(K-k)c.\label{equ3.14}
\end{eqnarray}
Since
\begin{equation*}
\frac{1}{p-k-1}\left(\frac{\theta_{K}^2}{\widetilde{\theta}_{p,K}}-(K-k)\right)\sim\frac{1}{p-k-1}\left(\exp\{2(\gamma_{K}-\mu_{K+1})\}-(K-k)\right)>0,
\end{equation*}
the second term of (\ref{equ3.14}) and then also (\ref{equ3.14}) tend to infinity as $p\rightarrow\infty$. Hence the second assertion holds.
\qed
\\

\noindent\textbf{Proof of Theorem \ref{thm3.2}.}
The proof of Theorem \ref{thm3.2} is identical to that of Theorem \ref{thm3.1} and hence omitted.	\qed\\

\noindent \textbf{Proof of Corollary \ref{cor1}.}
 We first verify \eqref{eivecconsistent}. From \eqref{l2normvandhatu} we find that for any fixed unit vector $\bu \in \mathbb{R}^p$,   \begin{equation}\label{buukdif}\inf_{t\in \mathbb{R}} \| t\bu-\hat{\bu}_k\|^2 = 1- \bu^* \hat{\bu}_k\hat{\bu}_k^* \bu = 1- \eta_k \frac{\bu^* \bA_n^* \bxi_k \bxi_k^* \bA_n \bu}{\gamma_k}+ O_P \left(\frac{1}{\sqrt{n}}\right),\end{equation}
 where the first step holds by taking $t =  \hat{\bu}_k^* \bu$. 
 Note that $\bA_n^*\bxi_k \bxi_k^* \bA_n$ is a rank one matrix and its eigenvector associated with the non-zero eigenvalue is $\tilde{\bu}:=\bA_n^* \bxi_k/\|\bA_n^*\bxi_k\|$.
 Then \eqref{eivecconsistent} follows by substituting $\bu = \tilde{\bu}$ into \eqref{buukdif} and using the fact that $(\bA_n\bA_n^*+\bSigma)\bxi_k = \gamma_k \bxi_k$.
 
The second statement \eqref{FnormhatU} can be concluded by finding that $\boldsymbol\Lambda =\bV_r^* \hat{\bU}_r$ minimizes $\|\bV_r \boldsymbol\Lambda - \hat{\bU}_r\|_F^2$ and its minimum value is obtained also by \eqref{l2normvandhatu}.
\qed
\\
 
\noindent\textbf{Proof of Theorem \ref{thm3.3}.}
We sketch the proofs here, which is quite similar to that of Theorem \ref{thm3.1}. For $k<K$, we have
\begin{eqnarray}
&&\frac{1}{n}\left(\acute{\text{EDA}}_{k}-\acute{\text{EDA}}_{K}\right)=\frac{1}{n}\sum_{i=k+1}^{K}\left(\acute{\text{EDA}}_{i-1}-\acute{\text{EDA}}_{i}\right)\notag\\
&=&\sum_{i=k+1}^{K}\left\{\log\theta_{i}+(n-i)\log\left[1-\frac{1}{n-i}\left(1-\frac{\theta_{i}^{2}}{\widetilde{\theta}_{n,i}}\right)\right]+\log\widetilde{\theta}_{n,i}-2\frac{p}{n}\right\},\label{equ3.12}\\
&&\frac{1}{n}\left(\acute{\text{EDB}}_{k}-\acute{\text{EDB}}_{K}\right)=\frac{1}{n}\sum_{i=k+1}^{K}\left(\acute{\text{EDB}}_{i-1}-\acute{\text{EDB}}_{i}\right)\notag\\
&=&\sum_{i=k+1}^{K}\left\{(\log p)(\log\theta_{i})+(n-i)\log\left[1-\frac{1}{n-i}\left(1-\frac{\theta_{i}^{2}}{\widetilde{\theta}_{n,i}}\right)\right]+\log\widetilde{\theta}_{n,i}-(\log n)\frac{p}{n}\right\}.\notag\\
\label{equ3.13}
\end{eqnarray}
According to the second assertion in Theorem \ref{thm2}, due to $h=o(p)$ bulks in $\Gamma_{F^{c,H}}$, we also have (\ref{equ3.5.4}). Thus,
\begin{equation*}
\widetilde{\theta}_{n,i}\sim1,~~i=2,\ldots,K.
\end{equation*}
When $\gamma_{1}<\infty$, for $i=2,\ldots,K$, we have $\theta_{i}\sim \zeta_{i}$  defined in (\ref{equ3.3}).
Hence, if the gap conditions (\ref{equ3.6}) and (\ref{equ3.6.5}) are satisfied, then
\begin{equation*}
\begin{split}
& \frac{1}{n}\left(\acute{\text{EDA}}_{k}-\acute{\text{EDA}}_{K}\right)\sim\sum_{i=k+1}^{K}\left\{\zeta_{i}^{2}+\log\zeta_{i}-1-2c\right\}\geq(K-k)\min\limits_{s=1,\ldots,K}a_{s}>0,\\
& \frac{1}{n}\left(\acute{\text{EDB}}_{k}-\acute{\text{EDB}}_{K}\right)\sim\sum_{i=k+1}^{K}\left\{\zeta_{i}^{2}+\log\zeta_{i}^{\log p}-1-(\log n)c\right\}\geq(K-k)\min\limits_{s=1,\ldots,K}b_{s}>0.
\end{split}
\end{equation*}
When $\gamma_{K}\rightarrow\infty$, by the similar discussion to that of Theorem \ref{thm3.1} with $n$ instead of $p$, without any gap conditions, we have
%%\begin{equation*}
%%\begin{split}
%%& \frac{1}{n}\left(\acute{\text{EDA}}_{k}-\acute{\text{EDA}}_{K}\right)\geq\left(\log n^{\delta_{K}}+n^{2\delta_{K}}-1-2(K-k)c\right)I_{\{\delta_{K}<\frac{1}{2}\}}\rightarrow\infty,\\
%%& \frac{1}{n}\left(\acute{\text{EDB}}_{k}-\acute{\text{EDB}}_{K}\right)\geq\left((\log p)(\log n^{\delta_{K}})+n^{2\delta_{K}}-1-(\log n)(K-k)c\right)I_{\{\delta_{K}<\frac{1}{2}\}}\rightarrow\infty,\\
%%\end{split}
%%\end{equation*}
%\begin{equation*}
%\begin{split}
%& \frac{1}{n}\left(\acute{\text{EDA}}_{k}-\acute{\text{EDA}}_{K}\right)\geq\lambda_{K}+(n-K)\log\left[1-\frac{1}{n-K}\left(1-\exp\{2\lambda_{K}\}\right)\right]-(K-k)(1+2c)\rightarrow\infty\\
%& \frac{1}{n}\left(\acute{\text{EDB}}_{k}-\acute{\text{EDB}}_{K}\right)\geq\lambda_{K}\log p+(n-K)\log\left[1-\frac{1}{n-K}\left(1-\exp\{2\lambda_{K}\}\right)\right]-(K-k)(1+c\log n)\rightarrow\infty.
%\end{split}
%\end{equation*}
\begin{equation*}
\begin{split}
\frac{1}{n}\left(\acute{\text{EDA}}_{k}-\acute{\text{EDA}}_{K}\right)
\geq  & \gamma_{k+1}-\mu_{K+1}+(n-K-1)\log\left[1+\frac{1}{n-k-1}\left(\frac{\theta_{K}^2}{\widetilde{\theta}_{n,K}}-(K-k)\right)\right]\\
& +(K-k)\log\widetilde{\theta}_{n,k}-2(K-k)c \rightarrow \infty,\\
\frac{1}{n}\left(\acute{\text{EDB}}_{k}-\acute{\text{EDB}}_{K}\right)
\geq& (\log p)(\gamma_{k+1}-\mu_{K+1})+(n-K-1)\log\left[1+\frac{1}{n-k-1}\left(\frac{\theta_{K}^2}{\widetilde{\theta}_{n,K}}-(K-k)\right)\right]\\
& +(K-k)\log\widetilde{\theta}_{n,k}-(\log n)(K-k)c \rightarrow \infty.
\end{split}
\end{equation*}

Next, consider $K<k\leq w=o(p)$. Analogously, it yields
\begin{equation*}
\begin{split}
&\frac{1}{n}\left(\acute{\text{EDA}_{k}}-
\acute{\text{EDA}_{K}}\right)\\
=&\sum_{i=K+1}^{k}\left\{-\log\theta_{i}-(n-i)\log\left[1-\frac{1}{n-i}\left(1-\frac{\theta_{i}^{2}}{\widetilde{\theta}_{n,i}}\right)\right]-\log\widetilde{\theta}_{n,i}+2\frac{p}{n}\right\}\notag\\
\sim&\sum_{i=K+1}^{k}\left\{1-\frac{\theta_{i}^{2}}{\widetilde{\theta}_{n,i}}-\log\theta_{i}-\log\widetilde{\theta}_{n,i}+2c\right\}\sim2(k-K)c>0,
\end{split}
\end{equation*}
and
\begin{equation*}
\begin{split}
&\frac{1}{n}\left(\acute{\text{EDB}_{k}}-\acute{\text{EDB}_{K}}\right)\\
=&\sum_{i=K+1}^{k}\left\{-(\log p)(\log\theta_{i})-(n-i)\log\left[1-\frac{1}{n-i}\left(1-\frac{\theta_{i}^{2}}{\widetilde{\theta}_{n,i}}\right)\right]-\log\widetilde{\theta}_{n,i}+(\log n)\frac{p}{n}\right\}\notag\\
\sim&\sum_{i=K+1}^{k}\left\{1-\frac{\theta_{i}^{2}}{\widetilde{\theta}_{n,i}}-(\log p)(\log\theta_{i})-\log\widetilde{\theta}_{n,i}+c\log n\right\}\sim(k-K)c\log n>0,
\end{split}
\end{equation*}
which completes the proof.
\qed

\subsection{Remaining proof}\label{secproofprop23}
This section includes the proofs of Propositions \ref{cenprop}, \ref{prop01}, Lemma \ref{sepl1}, and Theorem \ref{centhm3}.
Below are some lemmas required.
\begin{lemma}\label{m2}
For $n\times n$ invertible matrix $\bA$ and $n\times1$ vectors $\bq,\bv$ where $\bA$ and $\bA+\bv\bv^*$
are invertible, we have
\begin{equation*}
\bq^{*}\left(\bA+\bv \bv^{*}\right)^{-1}=\bq^{*} \bA^{-1}-\frac{\bq^{*} \bA^{-1} \bv}{1+\bv^{*} \bA^{-1} \bv} \bv^{*} \bA^{-1}.
\end{equation*}	
\end{lemma}

\begin{lemma}\label{m3}
Let $\bB=(b_{ij})\in\mR^{n\times n}$ with $\|\bB\|=O(1)$ and $\bx=(x_1,\ldots,x_n)^*$, where $x_i$ are i.i.d. satisfying $\Exp x_i=0$, ${\mathbb{E}}|x_i|^2=1$. Then, we have 
$$
\left|\bx^{*} \bB \bx-\trace \bB\right|^{q} \leq C_{q}\left(\left({\mathbb{E}}\left|x_{1}\right|^{4} \trace \bB \bB^{*}\right)^{q / 2}+{\mathbb{E}}\left|x_{1}\right|^{2 q} \trace\left(\bB \bB^{*}\right)^{q / 2}\right).
$$
\end{lemma}
%\begin{lemma}
%(Wely's inequality) Suppose that $A$ and  $B$ are two $n\times n$  Hermitian matrices, and there is %\begin{eqnarray*}
%	d_j+\lambda_k\le\gamma_i\le\lambda_r+d_s,
%\end{eqnarray*}
% where $j+k-n\ge i\ge r+s-1$	, $d_i$, $\lambda_i$ and $\gamma_i$ are the $i$-th eigenvalue of $A$, $B$ and $A+B$, respectively. 
%\end{lemma}

\begin{lemma}\label{burk}
	(Burkholder inequality) Let $\{X_k\}$ be a complex martingale difference sequence with respect to the filtration $\mathcal F_k$. For every $q\ge 1$, there exists $C_q>0$ such that: $${\mathbb{E}}\left|\sum^n_{k=1}X_k\right|^{2q}\le C_q\left({\mathbb{E}} \left(\sum^n_{k=1}{\mathbb{E}}(|X_k|^2|\mathcal F_{k-1})\right)^q+\sum^n_{k=1}{\mathbb{E}}|X_k|^{2q}\right).$$
\end{lemma}
For simplicity, we remove the subscripts of ``n''.
Let $\bX=[\bx_1,\ldots,\bx_n]$, $\bx_i=\ba_i+\bSigma^{1/2}\bw_i$, $\bX_k=\bX-\bx_k\be_k^*$, and hence define $$ \bQ_k(z)=(\bX_k\bX_k^*-z\bI)^{-1}.$$ 
Moreover, we also introduce some basic notations and formulas. For $n\times n$ invertible matrix $\bA,\bB$ and $n$-dimensional vector $\bq$, there are
\begin{equation}\label{eq03}
\bq^{*}\left(\bB+\bq \bq^{*}\right)^{-1}=\frac{1}{1+\bq^{*} \bB^{-1} \bq} \bq^{*} \bB^{-1},
\end{equation}
\begin{equation}\label{eq04}
\bA^{-1}-\bB^{-1}=\bB^{-1}(\bB-\bA)\bA^{-1}.
\end{equation}
Moreover, define \begin{equation}
\beta_k=\frac1{1+\bx_k^*\bQ_k(z)\bx_k},
\end{equation}
\begin{equation}\label{bk}
b_k=\frac1{1+\trace(\bSigma \bQ_k(z))/n+\ba_k^*\bQ_k(z)\ba_k}	.
\end{equation}
The following lemma is useful in calculating some moments bounds below:
\begin{lemma}\label{bd} For $z\in\mC_+$, there are
$|\beta_k|\le\frac{|z|}{|\Im z|}$, $|b_k|\le  \frac{|z|}{|\Im z|} $ and 
$\|\bQ_k(z)\bX_k\|\le (\frac1{|\Im z|}+\frac{|z|}{|\Im z|^2})^{1/2}$.	
\end{lemma}
\begin{proof}
We have $$|z^{-1}\beta_k| \le \frac{1}{\Im (z + z \bx_k^* \bQ_k(z) \bx_k)} \leq \frac{1}{\Im z},$$
where the second step uses the fact that $\Im (z \bx_k^* \bQ_k(z) \bx_k)>0 $. Therefore $|\beta_k|\le |z|/\Im z$. The bound for $|b_k|$ is checked similarly. 
For the last one, we have 
\begin{eqnarray*}
	\|\bQ_k(z)\bX_k\|&= &\|\bQ_k(z)\bX_k\bX_k^*\bQ_k(z)\|^{\frac12}\\
	&= &\|\bQ_k(z)(\bX_k\bX_k^*-z\bI+z\bI)\bQ_k(z)\|^{1/2}\\
	&\le &\|\bQ_k(z)+z\bQ_k(z)\bQ_k(z)\|^{1/2}\\
	&\le &(\frac1{|\Im z|}+\frac{|z|}{|\Im z|^2})^{1/2}.
\end{eqnarray*}	
\end{proof}

\noindent\textbf{Proof of Proposition \ref{prop01}.}
We first truncate, recentralize and renormalize the entries of $\bW$ following the steps in \cite{bai2007asymptotics}.
Select $\eta_n\to 0$ and satisfies 
$\eta_n^{-4} \int_{\left\{\left|n^{1/2}W_{11}\right| \geq \eta_n n^{1 / 4}\right\}}\left|n^{1/2}W_{11}\right|^4 \rightarrow 0$.
Let $\hat{W}_{i j}=W_{i j} I\left(\left|W_{i j}\right| \leq \eta_n n^{-1 / 4}\right)-{\mathbb{E}} W_{i j} I\left(\left|X_{i j}\right| \leq \eta_n n^{-1 / 4}\right), \tilde{\bW}_n=\bW_n-\hat{\bW}_n$, and $\hat{\bX}_n=\bA_n+\bSigma^{1/2}\hat{\bW}_n$, where $\hat{\bW}_n=\left(\hat{W}_{i j}\right)$. Let $\sigma_n^2={\mathbb{E}}\left|\hat{W}_{11}\right|^2$ and $\breve{\bX}_n=\bA+\sigma_n^{-1}\bSigma^{1/2} \hat{\bW}_n$. Write %$\hat{Q}^{-1}(z)=\left(\hat{\bW}_n\hat{\bW}_n^*-z \bI\right)^{-1}$ and 
$\breve{\bQ}(z)=\left(\breve{\bX}_n\breve{\bX}_n^*-z I\right)^{-1}$. Then following the arguments used in the proof Lemma 4 therein, we can show that ${\mathbb{E}} |\bv^*(\breve{\bQ}(z)-\bQ(z))\bv|^2 = o(n^{-1}).$

With this truncation and centralization, we have the following simple bound that can be checked using Lemma \ref{m3} and will be used frequently later: For any deterministic $p\times p$ matrix $\bA$ of bounded spectral norm,
\begin{equation}\label{quadbw}
{\mathbb{E}}\left|\bw_1^* \bA \bw_1-n^{-1} \trace \bA\right|^q \leq C\left(\eta_n^{2q-4} n^{-q / 2}+n^{-q / 2}\right) \le C n^{-q/2}.
\end{equation}
With this bound, and \begin{equation}\label{bdbilinear}{\mathbb{E}} |\bw_1^* \bv|^4 \le C n^{-2},
 \end{equation}
 for any deterministic unit norm vector $\bv$,
 we can obtain the following lemma without difficulty.
\begin{lemma}\label{delta}
  Let $\Delta_k=\bx_k^*\bQ_k(z)\bx_k-\frac{\trace \bSigma \bQ_k(z)}{n}-\ba_k^*\bQ_k(z)\ba_k$ and ${\mathbb{E}}_{\bw_k}$ be the conditional expectation with respect to the $\sigma$-field generated by  $\{\bw_l,l\neq k\}$. Under Assumption \ref{ass1} and Assumption \ref{ass3new}, for $1\le q\le 4$, we have 
 	\begin{equation*}
{\mathbb{E}}_{w_{k}}\left|\Delta_{k}\right|^{q}=O_P\left(\frac{1}{n^{q / 2} |\Im z|^{q}}\right).
\end{equation*}
A direct consequence is that 
\begin{equation*}
    {\mathbb{E}}_{w_k} |\beta_k - b_k|^q = O_P\left(\frac{|z|^q}{n^{q / 2} |\Im z|^{3q}}\right).
\end{equation*}
 \end{lemma}

\begin{proof}
%Note that under Assumption \ref{ass3}, it is easy to obtain that $\|\ba_i\|^2=O(1/n)$.
Using the Woodbury identity in Lemma \ref{wbm}, we write \begin{equation}\label{woodburrytQ}
\begin{aligned}
 \tilde \bQ(z)&= (-z\bI)^{-1}-(-z\bI)^{-1}\bX^*[\bI+\bX(-z)^{-1}\bX^*]^{-1}\bX (-z\bI)^{-1}
\\ 
 &=(-z\bI)^{-1}+z^{-1}\bX^*(\bX\bX^*-z\bI)^{-1}\bX.
 \end{aligned}
 \end{equation}
 To prove Proposition \ref{prop01}, it   suffices to prove \begin{eqnarray}
 	{\mathbb{E}}|\bu^*\bX^*(\bX\bX^*-z\bI)^{-1}\bX \bu-{\mathbb{E}} \bu^*\bX^*(\bX\bX^*-z\bI)^{-1}\bX \bu|^{2}\le Cn^{-1}
  \label{eq01},
  \end{eqnarray}
  and
 \begin{eqnarray} 
 |{\mathbb{E}} \bu^*\bX^*(\bX\bX^*-z\bI)^{-1}\bX \bu-{\mathbb{E}} \bu^*\bX_0^*(\bX_0\bX_0^*-z\bI)^{-1}\bX_0 \bu|\le C\frac{1}{\sqrt n},\label{eq02}
 \end{eqnarray}
where $\bu=(u_1,\ldots,u_n)^*$ is a fixed unit vector and $\bX_0$ represents the case of the $\bW$ being gaussian, denoted by $\bW_0$. Suppose, by singular value decomposition, $\bSigma=\bU \bD \bU^*$, we then have 
\begin{eqnarray*}
&&{\mathbb{E}} \bu^*\left[(\bA+ \bSigma^{1/2}\bW_0)^*(\bA+ \bSigma^{1/2}\bW_0)-z\bI\right]^{-1}\bu\\
&=&	{\mathbb{E}} \bu^*\left[(\bA+\bU\bD^{1/2}\bU^{*}\bW_0)^*(\bA+\bU\bD^{1/2}\bU^{*}\bW_0)-z\bI\right]^{-1}\bu\\
&=&	{\mathbb{E}} \bu^*\left[(\bA+\bU\bD^{1/2}\bW_0)^*(\bA+\bU\bD^{1/2}\bW_0)-z\bI\right]^{-1}\bu\\
&=&{\mathbb{E}} \bu^*\left[\big(\bU(\bU^*\bA+\bD^{1/2}\bW_0)\big)^*\big(\bU(\bU^*\bA+\bD^{1/2}\bW_0)\big)-z\bI\right]^{-1}\bu\\
&=&{\mathbb{E}} \bu^*\left[\big(\bU^*\bA+\bD^{1/2}\bW_
0\big)^*\big(\bU^*\bA+\bD^{1/2}\bW_0\big)-z\bI\right]^{-1}\bu.
\end{eqnarray*} 
 Letting $\bU^*\bA$ as $\bA$, it satisfies the model in \cite{hachem2013bilinear}. Hence we have \begin{equation}\label{Egausscase}
 	{\mathbb{E}} \left|\bu^*\left(\left[\big(\bU^*\bA+\bD^{1/2}\bW_
0\big)^*\big(\bU^*\bA+\bD^{1/2}\bW_0\big)-z\bI\right]^{-1}-\mathbf{T}'(z)\right)\bu\right|\le C\frac{1}{\sqrt n},
 \end{equation}
where $\mathbf{T}'(z)=\left(-z(1+\delta(z))\bI+\bA^*\bU(\bI+\tilde\delta(z)\bD)^{-1}\bU^*\bA\right)^{-1}=\tilde{\mathbf{T}}(z)$.
Moreover, combing Proposition 3.8 and Proposition 3.9 in \cite{hachem2013bilinear}, the conclusion follows.

{\textbf{Proof of \eqref{eq01}}: 
We will write the term in \eqref{eq01} as the sum of the martingale difference sequence first.
Using two basic matrix equality \eqref{eq03} and  \eqref{eq04}} and, we have
 \begin{eqnarray}\label{abc}
 	&&\bu^*\bX^*(\bX\bX^*-z\bI)^{-1}\bX \bu-\bu^*\bX_k^*(\bX_k\bX_k^*-z\bI)^{-1}\bX_k \bu\nonumber\\
 	&=&\bu^*(\bX^*-\bX_k^*)(\bX\bX^*-z\bI)^{-1}\bX \bu+\bu^*\bX_k^*\big(\bQ(z)-\bQ_k(z)\big)\bX \bu+\bu^*\bX_k^*\bQ_k(z)(\bX-\bX_k) \bu\nonumber\\
 	&=&\bu^*\be_k\bx_k^*\bQ_k(z)\bX \bu\beta_k-\bu^*\bX_k^*\bQ_k(z)\bx_k\bx_k^*\bQ_k(z)\bX \bu\beta_k+\bu^*\bX_k^*\bQ_k(z)\bx_k\be_k^*\bu\nonumber\\
 	&:=&A_k-B_k+C_k.
 \end{eqnarray}
 Denote by ${\mathbb{E}}_k$ the conditional expectation with respect to the $\sigma$-field generated by $\{w_i, i\le k\}$.
 With the above expansion, 
it is equivalent to obtaining a bound of 
$$ {\mathbb{E}} \left|\sum_{k=1}^n ({\mathbb{E}}_k - {\mathbb{E}}_{k-1})(A_k-B_k+C_k) \right|^{2}.$$

  We split $A_k$ as:
 \begin{eqnarray*}
 A_k&=&\bu^*\be_k	\bx_k^*\bQ_k(z)\bX \bu\beta_k\\
 &=& \bu^*\be_k(\ba_k^*+\bw_k^* \bSigma^{1/2})\bQ_k(z)(\bX_k+\bx_k\be_k^*)\bu\beta_k\\
 &=& \bu^*\be_k\ba_k^*\bQ_k(z)\bX_k \bu\beta_k+ \bu^*\be_k\ba_k^*\bQ_k(z)\bx_k\be_k^*\bu\beta_k\\
  &&+ \bu^*\be_k\bw_k^* \bSigma^{1/2}\bQ_k(z)\bX_k\bu\beta_k+ \bu^*\be_k\bw_k^* \bSigma^{1/2}\bQ_k(z)\bx_k\be_k^*\bu\beta_k\\
  &:=&A_{1k}+A_{2k}+A_{3k}+A_{4k}.
 \end{eqnarray*}
To obtain the bound for the term involving $A_k$, we consider the bounds of $A_{1k}$ to $A_{4k}$, respectively.

For $A_{1k}= \bu^*\be_k\ba_k^*\bQ_k(z)\bX_k\bu\beta_k$, we can decompose it as the sum of two components: 
\\$\bu^*\be_k\ba_k^*\bQ_k(z)\bX_k\bu b_k$ and $\bu^*\be_k\ba_k^*\bQ_k(z)\bX_k\bu(\beta_k-b_k)$. Since $({\mathbb{E}}_k - {\mathbb{E}}_{k-1})\bu^*\be_k\ba_k^*\bQ_k(z)\bX_k\bu b_k=0$, we have 
\begin{eqnarray*}
\sum_{k=1}^n{\mathbb{E}}_{k-1}\left|({\mathbb{E}}_k-{\mathbb{E}}_{k-1})A_{1k}\right|^{2}&\le &C\sum_{k=1}^n{\mathbb{E}}_{k-1}\left|\bu^*\be_k\ba_k^*\bQ_k(z)\bX_k\bu(\beta_k-b_k)\right|^2\\
&\le &C\sum_{k=1}^n {\mathbb{E}}_{k-1} \left\{|u_k|^2\|\bQ_k(z)\bX_k\|^2 {\mathbb{E}}_{w_k}|\beta_k-b_k|^2\right\}\\
&\le &Cn^{-1},
\end{eqnarray*}
where $u_k$ is the $k$-th coordinate of $\bu$, and the third lines uses Lemmas \ref{bd}, \ref{delta} and $\sum_{k=1} |u_k|^2 = 1$. 

Similarly,
\begin{equation}\label{A1k2q}
\begin{aligned}
\sum_{k=1}^n{\mathbb{E}}|({\mathbb{E}}_k-{\mathbb{E}}_{k-1})A_{1k}|^{2}&\le C\sum_{k=1}^n{\mathbb{E}} |\bu^*\be_k\ba_k^*\bQ_k(z)\bX_k\bu(\beta_k-b_k)|^{2}\\
&\le C\sum_{k=1}^n|u_k|^{2} {\mathbb{E}}\left(\|\bQ_k(z)\bX_k\| \cdot|\beta_k-b_k|\right)^{2}\\
&\le Cn^{-1}.
\end{aligned}
\end{equation}
Thus,  applying the Burkholder inequality in Lemma \ref{burk}, we obtain \begin{equation}\label{a1k}
 	{\mathbb{E}}\left|\sum^n_{k=1}({\mathbb{E}}_k-{\mathbb{E}}_{k-1})A_{1k}\right|^{2}\le Cn^{-1}.
 \end{equation}

\iffalse 
{\color{purple} For $A_{2k}=\bu^*\be_k\ba_k^*\bQ_k(z)\bx_k\be_k^*\bu\beta_k$ ({\color{red}?}), by expanding $\bx_k$ and using Lemma \ref{bd}, we have \begin{eqnarray*}
 	&&\sum_{k=1}^n\Exp_{k-1}\left|(\Exp_k-\Exp_{k-1})A_{2k}\right|^2\\&\le &C\sum_{k=1}^n\Exp_{k-1}|\bu^*\be_k\ba_k^*\bQ_k(z)\ba_k\be_k^*\bu\beta_k|^2+C\sum_{k=1}^n\Exp_{k-1}|\bu^*\be_k\ba_k^*\bQ_k(z) \bSigma^{1/2}\bw_k\be_k^*\bu\beta_k|^2\\
 	&\le &C\frac{\sum^n_{k=1}|u_k|^4}{n^2}
 \end{eqnarray*}
and \begin{eqnarray*}
\sum_{k=1}^n\Exp|(\Exp_k-\Exp_{k-1})A_{2k}|^{2q}&\le &C\sum_{k=1}^n\Exp |A_{2k}|^{2q}\\
&\le &C\Exp|\bu^*\be_k\ba_k^*\bQ_k(z)\ba_k\be_k^*\bu\beta_k|^{2q}+C\Exp|\bu^*\be_k\ba_k^*\bQ_k(z) \bSigma^{1/2}\bw_k\be_k^*\bu\beta_k|^{2q}\\
&\le &\frac{C\sum_{k=1}^n|u_k|^{4q}}{n^{2q}},
\end{eqnarray*}
where we use the Lemma 3.1 in \cite{hachem2013bilinear} in last inequality above. Combining with Lemma \ref{burk}, we also have \begin{equation}\label{a2k}
 	\Exp\left|\sum^n_{k=1}(\Exp_k-\Exp_{k-1})A_{2k}\right|^{2q}\le \frac{C}{n^{2q}}.
 \end{equation}
}
\fi

For $A_{2k}=\bu^*\be_k\ba_k^*\bQ_k(z)\bx_k\be_k^*\bu\beta_k$, by $\bx_k = \ba_k + \bSigma^{1/2}\bw_k$ and using Lemma \ref{bd}, we have
\begin{eqnarray}\label{A2kdecomp}
     ({\mathbb{E}}_k- {\mathbb{E}}_{k-1}) A_{2k} = ({\mathbb{E}}_k- {\mathbb{E}}_{k-1})\left[|u_k|^2 \ba_k^*\bQ_k(z)\ba_k(\beta_k-b_k)+|u_k|^2 \ba_k^*\bQ_k(z)\bSigma^{1/2}\bw_k\beta_k\right]
\end{eqnarray}
Then
\begin{eqnarray*}
 &&\sum_{k=1}^n{\mathbb{E}}_{k-1}\left|({\mathbb{E}}p_k-{\mathbb{E}}_{k-1})A_{2k}\right|^2\\
 &\le &C\sum_{k=1}^n |u_k|^4{\mathbb{E}}_{k-1}| \ba_k^*\bQ_k(z)\ba_k(\beta_k-b_k)|^2+C\sum_{k=1}^n  |u_k|^4 {\mathbb{E}}_{k-1}| \ba_k^*\bQ_k(z) \bSigma^{1/2}\bw_k\beta_k|^2\\
 & \le & C\sum_{k=1}^n |u_k|^4 {\mathbb{E}}_{k-1}  |\beta_k-b_k|^2+C\sum_{k=1}^n|u_k|^4 {\mathbb{E}}_{k-1}|\ba_k^*\bQ_k(z) \bSigma^{1/2}\bw_k|^2\\
 &\le & Cn^{-1}
 \end{eqnarray*}
where the second step uses $|\beta_k|=O(1)$ and $\|\bQ_k(z)\| =O(1)$, and the third step uses Lemma \ref{delta} and \eqref{bdbilinear}.
By \eqref{A2kdecomp} and an argument similar to \eqref{A1k2q}, it can also be checked that \begin{eqnarray*}
\sum_{k=1}^n{\mathbb{E}}|({\mathbb{E}}_k-{\mathbb{E}}_{k-1})A_{2k}|^{2}&\le &C\sum_{k=1}^n{\mathbb{E}} |A_{2k}|^{2}\le Cn^{-1}.
\end{eqnarray*}
Therefore an application of the Burkholder inequality yields \begin{equation*}\label{a2k}
 	{\mathbb{E}}\left|\sum^n_{k=1}({\mathbb{E}}_k-{\mathbb{E}}_{k-1})A_{2k}\right|^{2}\le Cn^{-1}.
 \end{equation*}

For $A_{3k}=\bu^*\be_k\bw_k^* \bSigma^{1/2}\bQ_k(z)\bX_k \bu\beta_k$, it can be handled following an argument similar to the one that leads to the bound for $A_{2k}$.
%we use $$(\Exp_k - \Exp_{k-1})A_{3k} = \bu^*\be_k\bw_k^* \bSigma^{1/2}\bQ_k(z)\bX_k \bu(\beta_k-b_k) $$

 For $A_{4k}=\bu^*\be_k\bw_k^* \bSigma^{1/2}\bQ_k(z)\bx_k\be_k^*\bu\beta_k$, we have \begin{eqnarray}
 A_{4k}&=&|u_k|^2 \bw_k^* \bSigma^{1/2}\bQ_k(z)( \bSigma^{1/2}\bw_k+\ba_k)(\beta_k-b_k+b_k)\nonumber\\
 	&=& |u_k|^2\bw_k^* \bSigma^{1/2}\bQ_k(z) \bSigma^{1/2}\bw_k (\beta_k-b_k)+|u_k|^2\bw_k^* \bSigma^{1/2}\bQ_k(z) \bSigma^{1/2}\bw_k b_k\nonumber\\
 	&&+ |u_k|^2\bw_k^* \bSigma^{1/2}\bQ_k(z)\ba_k \beta_k\nonumber\\
 	&:=&A_{5k}+A_{6k}+A_{7k}.
 \end{eqnarray}

 Now, we can continue to bound  $A_{5k}$ to $A_{7k}$. For $A_{5k}=|u_k|^2\bw_k^*\bSigma^{1/2}\bQ_k(z)\bSigma^{1/2}\bw_k (\beta_k-b_k)$, we find that \begin{eqnarray*}
	&&\sum_{k=1}^n{\mathbb{E}}_{k-1}\left|({\mathbb{E}}_k-{\mathbb{E}}_{k-1})A_{5k}\right|^2 \\
 &\le &C\sum_{k=1}^n  |u_k|^4{\mathbb{E}}_{k-1} \left\{\left({\mathbb{E}}_{\bw_k}|\bw_k^*\bSigma^{1/2}\bQ_k(z)\bSigma^{1/2}\bw_k|^4\right)^{1/2}\left({\mathbb{E}}_{\bw_k}|\beta_k-b_k|^4\right)^{1/2}\right\}\\
	&\le &C\frac{\sum_{k=1}^n|u_k|^4}{n},
\end{eqnarray*}
where we apply \eqref{quadbw} and Lemma \ref{delta}. Similarly, 
\begin{eqnarray*}
\sum_{k=1}^n{\mathbb{E}}|({\mathbb{E}}_k-{\mathbb{E}}_{k-1})A_{5k}|^{2}&\le &C\sum_{k=1}^n{\mathbb{E}} |A_{5k}|^{2}\\
&\le &C\sum_{k=1}^n |u_k|^{2}{\mathbb{E}}\left\{\left({\mathbb{E}}_{\bw_k}|\bw_k^*\bSigma^{1/2}\bQ_k(z)\bSigma^{1/2}\bw_k|^{4}\right)^{1/2}\left({\mathbb{E}}_{\bw_k}|\beta_k-b_k|^{4}\right)^{1/2}\right\}\\
&\le &Cn^{-1}.
\end{eqnarray*}
Thus, again using Lemma \ref{burk}, we have 
 \begin{equation*}\label{a5k}
 	{\mathbb{E}}\left|\sum^n_{k=1}({\mathbb{E}}_k-{\mathbb{E}}_{k-1})A_{5k}\right|^{2}\le Cn^{-1}.
 \end{equation*}
For $A_{6k}=|u_k|^2\bw_k^*\bSigma^{1/2}\bQ_k(z)\bSigma^{1/2}\bw_k b_k$, we have \begin{eqnarray*}
 	\sum_{k=1}^n{\mathbb{E}}_{k-1}\left|({\mathbb{E}}_k-{\mathbb{E}}_{k-1})A_{6k}\right|^2&=&\sum_{k=1}^n |u_k|^4{\mathbb{E}}_{k-1}\left\{\left({\mathbb{E}}_{\bw_k}|\bw_k^*\bSigma^{1/2}\bQ_k(z)\bSigma^{1/2}\bw_k-\frac1p\trace\bSigma \bQ(z)|^2\right)\right\}\\
 	&\le &C\frac{\sum_{k=1}^n|u_k|^4}{n}
 \end{eqnarray*}
 and 
 \begin{eqnarray*}
 	\sum_{k=1}^n{\mathbb{E}}|({\mathbb{E}}_k-{\mathbb{E}}_{k-1})A_{6k}|^{2}&\le &C\sum_{k=1}^n{\mathbb{E}}\left\{|u_k|^{2}\left|{\mathbb{E}}_{k}(\bw_k^*\bSigma^{1/2}\bQ_k(z)\bSigma^{1/2}\bw_k-\frac1p\trace\bSigma \bQ(z))\right|^{2}\right\}\\
 	&\le &Cn^{-1}.
 \end{eqnarray*} 
 Thus,  we find
\begin{equation*}
 	{\mathbb{E}}\left|\sum^n_{k=1}({\mathbb{E}}_k-{\mathbb{E}}_{k-1})A_{6k}\right|^{2}\le Cn^{-1}.
 \end{equation*}
 The term invovling  $A_{7k}=\bu^*\be_k\bw_k^*\bSigma^{1/2}\bQ_k(z)\ba_k\be_k^*\bu\beta_k$ can be bounded by an argument similar to that for $A_{2k}$. 
  %Then, there is \begin{equation}\label{a7k}\Exp\left|\sum^n_{k=1}(\Exp_k-\Exp_{k-1})A_{7k}\right|^{2q}\le \frac{C}{n^{q}}. \end{equation}
Combining the above discussions we obtain
\begin{equation}\label{bdak}
    {\mathbb{E}}\left|\sum^n_{k=1}({\mathbb{E}}_k-{\mathbb{E}}_{k-1})A_{k}\right|^{2}\le Cn^{-1}.
\end{equation}

 Now, we consider $B_k$ in \eqref{abc}. We split $B_k$ into several components: \begin{eqnarray*}
 B_k&=& \bu^*\bX_k^*\bQ_k(z)\bx_k\bx_k^*\bQ_k(z)\bX \bu\beta_k\\
 &=&	\bu^*\bX_k^*\bQ_k(z)\ba_k\ba_k^*\bQ_k(z)\bX_k \bu\beta_k+\bu^*\bX_k^*\bQ_k(z)\bSigma^{1/2}\bw_k\ba_k^*\bQ_k(z)\bX_k \bu\beta_k\\
 &&+\bu^*\bX_k^*\bQ_k(z)\ba_k\bw_k^*\bSigma^{1/2}\bQ_k(z)\bX_k \bu\beta_k+\bu^*\bX_k^*\bQ_k(z)\bSigma^{1/2}\bw_k\bw_k^* \bSigma^{1/2}\bQ_k(z)\bX_k \bu(\beta_k-b_k)\\
 &&+\bu^*\bX_k^*\bQ_k(z)\bSigma^{1/2}\bw_k\bw_k^*\bSigma^{1/2}\bQ_k(z)\bX_k \bu b_k+\bu^*\bX_k^*\bQ_k(z)\ba_k\bx_k^*\bQ_k(z)\bx_k\be_k^* \bu\beta_k\\
 &&+\bu^*\bX_k^*\bQ_k(z)\bSigma^{1/2}\bw_k\bx_k^* \bQ_k(z)\bx_k\be_k^* \bu\beta_k\\
 &:=&B_{1k}+B_{2k}+B_{3k}+B_{4k}+B_{5k}+B_{6k}+B_{7k}.
 \end{eqnarray*}

We obtain bounds for $B_{1k}$, $B_{2k}$ and $B_{3k}$ by arguments similar to those leading to the bounds for $A_{1k}$ and $A_{2k}$.  For the terms $B_{4k}$ and $B_{5k}$, we utilize \eqref{bdbilinear}. 
For $B_{6k}$ we use $\bx_k = \ba_k + \bSigma^{1/2}\bw_k$ to further decompose it into four components. For the component without $\bw_k$, it can be bounded following an argument similar to the one that leads to the bound for $A_{1k}$. For the component with one $\bw_k$, it can be handled similarly to $A_{2k}$. For the component involving the quadratic form $\bw_k^* \bSigma^{1/2}\bQ_k(z) \bSigma^{1/2}\bw_k$, we use arguments leading to the bound for $A_{5k}$ and $A_{6k}$. For $B_{7k}$, it is similar to $B_{6k}$, and owing to the presence of $\bu^*\bX_k^*\bQ_k(z)\bSigma^{1/2}\bw_k$, which has a $4$-th moment of $O(n^{-2})$, its analysis becomes even simpler. Therefore, we find 
\begin{equation}\label{bdbk}
  {\mathbb{E}}\left|\sum^n_{k=1}({\mathbb{E}}_k-{\mathbb{E}}_{k-1})B_{k}\right|^{2}\le Cn^{-1}.
\end{equation}

Recalling the definition of $C_k$ in \eqref{abc}, according to the analysis of $A_k$,
we readily obtain that 
 \begin{equation}\label{bdck}
	{\mathbb{E}}\left|\sum^n_{k=1}({\mathbb{E}}_k-{\mathbb{E}}_{k-1})C_{k}\right|^{2}\le Cn^{-1}.
 \end{equation}
 
 Combining the bounds in \eqref{bdak}, \eqref{bdbk},and \eqref{bdck}, we can obtain \eqref{eq01}. 
 %Taking $q=2$, we can prove that $$\bu^*\bX^*(\bX\bX^*-z\bI)^{-1}\bX \bu\rightarrow\Exp \bu^*\bX^*(\bX\bX^*-z\bI)^{-1}\bX \bu \text{ ~a.s. }$$ \qed

\textbf{Proof of \eqref{eq02}}: 
We first define \begin{eqnarray*}
 \bZ_k^1&=&\sum_{i=1}^k\bx_i\be_i^*+\sum_{i=k+1}^n\bx_i^0\be_i^*\\
 \bZ_k&=&\sum_{i=1}^{k-1}\bx_i\be_i^*+\sum_{i=k+1}^n\bx_i^0\be_i^*\\
 \bZ_k^0&=&\sum_{i=1}^{k-1}\bx_i\be_i^*+\sum_{i=k}^n\bx_i^0\be_i^*,\\
 \end{eqnarray*}
where $\bx^0_i=\ba_i+ \bSigma^{1/2}\bw_i^0$, and $\bw_i^0$ follows  normal distribution with mean $\mathbf0$ variance $1/n$. 
Define 
\begin{equation}
    \bG_k(z)=(\bZ_k\bZ_k^*-z\bI)^{-1}, ~~~~\beta_k^1=
 \frac1{1+\bx_k^*\bG_k(z)\bx_k},~~~~ \beta_k^0~=~\frac1{1+{\bx_k^0}^*\bG_k(z)\bx_k^0}.
\end{equation}
Write \begin{eqnarray*}
 &&{\mathbb{E}} \bu^*\bX^*(\bX\bX^*-z\bI)^{-1}\bX \bu-{\mathbb{E}} \bu^*\bX_0^*(\bX_0\bX_0^*-z\bI)^{-1}\bX_0 \bu\\
&=&\sum^n_{k=1}{\mathbb{E}}\left( \bu^* {\bZ_k^{1}}^*({\bZ_k^{1}}{\bZ_k^{1}}^*-z\bI)^{-1}{\bZ_k^{1}}\bu-\bu^* {\bZ_k}^*({\bZ_k}{\bZ_k}^*-z\bI)^{-1}{\bZ_k}\bu\right)\\
&&-\sum^n_{k=1}{\mathbb{E}}\left( \bu^* {\bZ_k^{0}}^*({\bZ_k^{0}}{\bZ_k^{0}}^*-z\bI)^{-1}{\bZ_k^{0}}\bu-\bu^* {\bZ_k}^*({\bZ_k}{\bZ_k}^*-z\bI)^{-1}{\bZ_k}\bu\right)\\
&:=&\sum_{k=1}^n\left[{\mathbb{E}}\left(A_k^1-B_k^1+C_k^1\right)-{\mathbb{E}}\left(A_k^0-B_k^0+C_k^0\right)\right],
 \end{eqnarray*}
 where
\begin{equation*}
\begin{aligned}
 A_k^1 &= \bu^*\be_k\bx_k \bG_k \bZ_k^1\bu\beta_k^1, &
 B_k^1 &= \bu^*\bZ_k^*\bG_k\bx_k\bx_k^*\bG_k \bZ_k^1\bu\beta_k^1, &
 C_k^1 &= \bu^*\bZ_k^*\bG_k\bx_k\be_k^*\bu, \\
 A_k^0 &= \bu^*\be_k\bx_k^0 \bG_K \bZ_k^0\bu\beta_k^0, &
 B_k^0 &= \bu^*\bZ_k^*\bG_k\bx_k^0{\bx_k^0}^*\bG_k \bZ_k^0\bu\beta_k^0, &
 C_k^0 &= \bu^*\bZ_k^* \bG_k\bx_k^0\be_k^*\bu.
\end{aligned}
\end{equation*}

 Similar to $A_k, B_k, C_k$ in \eqref{abc}, here, $A_k^1, B_k^1,C_k^1,A_k^0, B_k^0, C_k^0$ can be further decomposed as before, and we use the superscripts ``1'' and ``0'' to distinguish the general case and the gaussian case. Since the procedure is similar as before, for simplicity, we list two typical examples to illustrate the proof idea. 
 	 For example, consider $A_k^1$,
 \begin{eqnarray*}
 A_k^1&=&\bu^*\be_k	\bx_k^*\bG_k(z)\bZ_k^1 \bu\beta_k^1\\
 &=&\bu^*\be_k(\ba_k^*+\bw_k^* \bSigma^{1/2})\bG_k(z)(\bZ_k+\bx_k\be_k^*)\bu\beta_k^1\\
 &=&\bu^*\be_k\ba_k^*\bG_k(z)\bZ_k\bu\beta_k^1+\bu^*\be_k\ba_k^*\bG_k(z)\bx_k\be_k^*\bu\beta_k^1\\
  &&+\bu^*\be_k\bw_k^* \bSigma^{1/2}\bG_k(z)\bZ_k\bu\beta_k^1+\bu^*\be_k\bw_k^* \bSigma^{1/2}\bG_k(z)\bx_k\be_k^*\bu\beta_k^1\\
  &:=&A_{1k}^1+A_{2k}^1+A_{3k}^1+A_{4k}^1.
 \end{eqnarray*}
To handle the term involving $A_{1k}^1=\bu^*\be_k\ba_k^*\bG_k(z)\bZ_k\bu\beta_k^1$,  we find that \begin{eqnarray*}
 \begin{aligned}
 \left|\sum^n_{k=1}{\mathbb{E}} \left[\bu^*\be_k\ba_k^*\bG_k(z)\bZ_k\bu(\beta_k^1-b_k)\right]\right|&	\le \sum^n_{k=1}\left({\mathbb{E}}|\bu^*\be_k\ba_k^*\bG_k(z)\bZ_k\bu|^2{\mathbb{E}}|\beta_k^1-b_k|^2\right)^{1/2}\\
& \le \frac{C}{\sqrt{n}} \sum_{k=1}^n |u_k|\cdot \|\ba_k\| \le \frac{C}{2\sqrt{n}}\sum_{k=1}^n(|u_k|^2+ \|\ba_k\|^2) \le \frac C{\sqrt n},
 \end{aligned}
 \end{eqnarray*}
 where 	$b_k$ is defined in \eqref{bk}, and in the second step we use ${\mathbb{E}}|\beta_k^1-b_k|^2 = O(n^{-1})$. Thus, we have \begin{equation}\label{a1k1}
 \sum^n_{k=1}{\mathbb{E}} 	A_{1k}^1=\sum^n_{k=1}{\mathbb{E}} \bu^*\be_k\ba_k^*\bG_k(z)\bZ_k\bu b_k+O\left(\frac1{\sqrt n}\right).
 \end{equation}
Similarly, we have $\sum^n_{k=1}{\mathbb{E}} 	A_{1k}^0=\sum^n_{k=1}{\mathbb{E}} \bu^*\be_k\ba_k^*\bG_k(z)\bZ_k\bu b_k+O(\frac1{\sqrt n})$. 
For $A_{3k}^1$, we have \\${\mathbb{E}}\bu^*\be_k\bw_k^* \bSigma^{1/2}\bG_k(z)\bZ_k\bu b_k=0$. Then, by Cauchy–Schwarz inequality, write 
\begin{equation}
    \begin{aligned}
     \sum_{k=1}^n{\mathbb{E}} \left[u_k\bw_k^* \bSigma^{1/2}\bG_k(z)\bZ_k\bu(\beta_k^1-b_k)\right]\le \frac Cn\sum_{k=1}^n  |u_k|({\mathbb{E}}\bu^*\bZ_k^*\bG_k(z)\bSigma \bG_k(z)\bZ_k\bu)^{1/2}\le \frac{C}{\sqrt{n}}.
    \end{aligned}
\end{equation}
Similarly, it is easy to prove that 
$\left|\sum^n_{k=1}{\mathbb{E}} A_{jk}^1\right|=O(\frac1{\sqrt n})$ for $j=2,4$ by using Cauchy–Schwarz inequality.

For $B_k^1$ , we write
\begin{eqnarray*}
 B_k^1&=&\bu^*\bZ_k^*(\bZ_k\bZ_k^*-z\bI)^{-1}\bx_k\bx_k^*(\bZ_k\bZ_k^*-z\bI)^{-1}\bZ_k\bu\beta_k^1\\
 &=&	\bu^*\bZ_k^*\bG_k(z)\ba_k\ba_k^*\bG_k(z)\bZ_k \bu\beta_k^1+\bu^*\bZ_k^*\bG_k(z) \bSigma^{1/2}\bw_k\ba_k^*\bG_k(z)\bZ_k \bu\beta_k^1\\
 &&+\bu^*\bZ_k^*\bG_k(z)\ba_k\bw_k^* \bSigma^{1/2}\bG_k(z)\bZ_k \bu\beta_k^1+\bu^*\bZ_k^*\bG_k(z) \bSigma^{1/2}\bw_k\bw_k^* \bSigma^{1/2}\bG_k(z)\bZ_k \bu(\beta_k^1-b_k)\\
 &&+\bu^*\bZ_k^*\bG_k(z) \bSigma^{1/2}\bw_k\bw_k^* \bSigma^{1/2}\bG_k(z)\bZ_k \bu b_k+\bu^*\bZ_k^*\bG_k(z)\ba_k\bx_k^*\bG_k(z)\bx_k\be_k^* \bu\beta_k^1\\
 &&+\bu^*\bZ_k^*\bG_k(z) \bSigma^{1/2}\bw_k\bx_k\bG_k(z)\bx_k\be_k^* \bu\beta_k^1\\
 &:=&B_{1k}^1+B_{2k}^1+B_{3k}^1+B_{4k}^1+B_{5k}^1+B_{6k}^1+B_{7k}^1.
 \end{eqnarray*}
 Similar to \eqref{a1k1}, we have \begin{equation*} \sum^n_{k=1}{\mathbb{E}} B_{1k}^1=\sum^n_{k=1}{\mathbb{E}} \left[\bu^*\bZ_k^*\bG_k(z)\ba_k\ba_k^*\bG_k(z)\bZ_k \bu b_k\right]+O\left(\frac1{\sqrt n}\right),
 \end{equation*}
  and
 \begin{equation*} \sum^n_{k=1}{\mathbb{E}} 	B_{1k}^0=\sum^n_{k=1}{\mathbb{E}} \left[\bu^*\bZ_k^*\bG_k(z)\ba_k\ba_k^*\bG_k(z)\bZ_k \bu b_k\right]+O\left(\frac1{\sqrt n}\right)	.
 \end{equation*}
 For $B_{2k}^1=\bu^*\bZ_k^*\bG_k(z) \bSigma^{1/2}\bw_k\ba_k^*\bG_k(z)\bZ_k \bu\beta_k^1$, we have $$\sum_{k=1}^n {\mathbb{E}} B_{2k}^1= \sum_{k=1}^n {\mathbb{E}}\left[\bu^*\bZ_k^*\bG_k(z) \bSigma^{1/2}\bw_k\ba_k^*\bG_k(z)\bZ_k \bu(\beta_k^1-b_k)\right]=O\left(\frac{1}{\sqrt{n}}\right),
 $$ and by the same reason, we have such bound for $\sum_{k=1}^n{\mathbb{E}} B_{2k}^0$, $\sum_{k=1}^n{\mathbb{E}} B_{3k}^1$ and   $\sum_{k=1}^n{\mathbb{E}} B_{3k}^0$.

For $B_{4k}^1=\bu^*\bZ_k^*\bG_k(z)\bSigma^{1/2}\bw_k\bw_k^*\bSigma^{1/2}\bG_k(z)\bZ_k \bu(\beta_k^1-b_k)$, by lemma \ref{m3} we have $$
 \begin{aligned}
     |\sum_{k=1}^n {\mathbb{E}} B_{4k}^1| &\leq \sum_{k=1}^n\left[{\mathbb{E}}|\bw_k^*\bSigma^{1/2}\bG_k(z)\bZ_k \bu\bu^*\bZ_k^*\bG_k(z)\bSigma^{1/2}\bw_k|^2 {\mathbb{E}}|\beta_k^1-b_k|^2\right]^{1/2}\\
     &\le \frac{C}{\sqrt{n}}\sum_{k=1}^n\Big[{\mathbb{E}}|\bw_k^*\bSigma^{1/2}\bG_k(z)\bZ_k \bu\bu^*\bZ_k^*\bG_k(z)\bSigma^{1/2}\bw_k-n^{-1}\bu^*\bZ_k\bG_k(z)\bSigma \bG_k(z)\bZ_k\bu|^2\\
     &~~+{\mathbb{E}}(\bu^*\bZ_k\bG_k(z)\bSigma \bG_k(z)\bZ_k\bu/n)^2\Big]^{1/2}=O\left(\frac{1}{\sqrt{n}}\right).
 \end{aligned}
 $$
 and $\sum_{k=1}^p B_{4k}^0$ also has bound of order $O(\frac{1}{\sqrt{n}})$.
 
 For $B_{5k}^1$ and $B_{5k}^0$, we have 
 $$\sum_{k=1}^n {\mathbb{E}} B_{5k}^1=\sum_{k=1}^n {\mathbb{E}} B_{5k}^0=\sum_{k=1}^n \frac{\bu^*\bZ_k^*\bG_k(z)\bSigma \bG_k(z)\bZ_k \bu}{n}.$$
\noindent For $B_{6k}^1=\bu^*\bZ_k^*\bG_k(z)\ba_k\bx_k^*\bG_k(z)\bx_k\be_k^* \bu\beta_k^1$, it can be decomposed into  \begin{eqnarray*}
 	B_{6k}^1&=&\bu^*\bZ_k^*\bG_k(z)\ba_k\ba_k^*\bG_k(z)\ba_k\be_k^* \bu\beta_k^1+\bu^*\bZ_k^*\bG_k(z)\ba_k\bw_k^*\bSigma^{1/2}\bG_k(z)\ba_k\be_k^* \bu\beta_k^1\\
 	&&+\bu^*\bZ_k^*\bG_k(z)\ba_k\ba_k^*\bG_k(z)\bSigma^{1/2}\bw_k\be_k^* \bu\beta_k^1+\bu^*\bZ_k^*\bG_k(z)\ba_k\bw_k^*\bSigma^{1/2}\bG_k(z)\bSigma^{1/2}\bw_k\be_k^* \bu(\beta_k^1-b_k)\\
 &&+\bu^*\bZ_k^*\bG_k(z)\ba_k\bw_k^*\bSigma^{1/2}\bG_k(z)\bSigma^{1/2}\bw_k\be_k^* \bu b_k,
 \end{eqnarray*}
and it is readily verified that $$ \begin{aligned}
    \sum^n_{k=1}{\mathbb{E}} 	B_{6k}^1&=\sum^n_{k=1}{\mathbb{E}}\left[\frac{\trace [\bG_k(z)\bSigma]}{n}\bu^*\bZ_k^*\bG_k(z)\ba_ku_kb_k+\bu^*\bZ_k^*\bG_k(z)\ba_k\ba_k^*\bG_k(z)\ba_ku_kb_k\right]+O\left(\frac1{\sqrt n}\right)\\
    &=\sum^n_{k=1}{\mathbb{E}} 	B_{6k}^0+O\left(\frac1{\sqrt n}\right).
\end{aligned} $$
Similarly, by decomposing $B_{7k}^1=\bu^*\bZ_k^*\bG_k(z)\bSigma^{1/2}\bw_k\bx_k\bG_k(z)\bx_ku_k\beta_k^1$, one can prove that \begin{equation}\label{e3602}
    \sum^n_{k=1}{\mathbb{E}} B_{7k}^1=2\ba_k^*\bG_k(z)\bSigma \bG_k(z)\bZ_k^*\bu u_kb_k+O\left(\frac{1}{\sqrt{n}}\right)= \sum^n_{k=1}{\mathbb{E}} B_{7k}^0+O\left(\frac{1}{\sqrt{n}}\right).
    \end{equation}

Therefore, combining arguments above, \eqref{eq02} holds.
\end{proof}

\noindent\textbf{Proof of Lemma \ref{sepl1}.}
Theorem 1 in \cite{liu2022ieee} has established this conclusion with an additional assumption that $\bA$ contains a finite number of different columns. 
%We use two steps to extend to general low-rank $\bA$. The first step is to show that for general low rank $\bA$ the conclusion holds when $\bW$ is Gaussian. And the second step is to extend to general $\bW$ satisfying Assumption \ref{ass1}.
We use three steps to extend to $\bA$ satisfying Assumption \ref{ass3new}. First, we extend Proposition 1 to $\bA$ satisfying Assumptions \ref{ass1},\ref{ass2} and \ref{ass3new}. Second, we show that the exact separation of eigenvalues holds with an additional assumption that $\bW$ is Gaussian and holds almost surely.  Third, we extend the result to general $\bW$ satisfying Assumption \ref{ass1}.

We start from the first step. In \cite{liu2022ieee}, the fact that the spectral norm of $\ba_i$ is $O(n^{-1/2})$  for $i=1,\cdots,n$ has been used in the proof. This fact is a consequence of the assumption that $\bA$ contains a finite number of different columns and the spectral norm of $\bA$ is bounded. Such bound on $\|\ba_i\|$ leads to those terms with the form $\bw_k^* R_1 \ba_k$ or $\ba_k^* R_2\ba_k$ trivially negligible, where  $R_1$ and $R_2$ have bounded norm (or some power of $v_n^{-1}$ when dealing with some resolvents) and independent of $\bw_k$.

%some bounds that are trivially negligible especially for those terms with the form $\bw_k^* R_1 \ba_k$ or $\ba_k^* R_2\ba_k$, for some $R_1$ and $R_2$ that have bounded norm (or some power of $v_n^{-1}$ when dealing with some resolvents) and independent of $\bw_k$.  %The rank assumption in that work is a consequence of the norm assumption on the columns of $\bA$ as $K \left[\lambda_{K}(\bA)\right]^2\le \|\bA\|_{Fr}^2 = O(1)$. Another reason is to make sure the limiting spectral distribution is the same as that without the deterministic perturbation but assuming the rank is $o(n)$ suffices to guarantee this.
Through careful check,  assuming $\|\bA\|=O(1)$ and $\trace (\bA\bA^*) = O(n^{1/3})$ is sufficient to proceed the proof of Proposition 1 in \cite{liu2022ieee}. The key observation is that those bounds involve individual $\ba_i$ can be aggregated by summing from $i=1$ to $n$,  resulting in bounds that are still negligible. The summation from $i=1$ to $n$ appears naturally in the previous proof as we study the difference between $m_n$ and its deterministic counterpart.  
%Notice that $r \left[\lambda_{r}(\bA)\right]^2\le \|A\|_{Fr}^2 = \sum_{i=1}^n \|\ba_i\|^2 = O(n^{1/2})$, and the proof of Proposition 1 does not use the finite rank of $\bA$. Therefore 
We can finish our first step by validating the bounds in (23) and (25) therein with some modifications to steps used to bound terms involving $\bA$ or $\ba_k$. 
%This implies we can assume $K = O(n^{1/2})$ and then using the same strategy in Lemma 3 of the manuscript we conclude the exact separation results.
Below, we list several modifications under the new assumption $\|\bA\|=O(1)$ and $\trace (\bA\bA^*) = O(n^{1/3})$ (the equation labels and notations refer to those in \cite{liu2022ieee}):
\begin{enumerate}
\item For the expansion of $\hat{\omega}_n$ below equation (19), the 
     terms $\mathbb{E} (\beta_k \bw_k^* \bSigma^{1/2} T \bQ_k \ba_k)$
     is bounded by $O(n^{-1/2}v_n^{-3})$  using $\mathbb{E} |\bw_k^* \bv|^q = n^{-q/2}\|\bv\|^q$ for any fixed $\bv$. The second last term can be bounded via $p^{-1}\sum_{k=1}^n \mathbb{E}|\beta_k \ba_k^* T\bQ_k \ba_k| \le p^{-1}\sum_{k=1}^n \mathbb{E} \left(|\beta_k| \|\ba_k\|^2 \|T\bQ_k\|\right) = p^{-1}v_n^{-3} \trace (\bA\bA^*) = O(n^{-2/3}v_n^{-3})$. The last term can still be bounded by $O(n^{-1/2}v_n^{-3})$ using $\trace (\bA\bA^*) =O(n^{1/3})$.
     
    \item For the moment bound for $\hat{\gamma}_k=\bx_k^* \bQ_k \bx_k-n^{-1} \trace \bSigma \bQ_k-\ba_k^* \bQ_k \ba_k$ below equation (23), we can get the same bound using $\mathbb{E} |\bw_k^* \bv|^q = n^{-q/2}\|\bv\|^q$ for any fixed $\bv$. 
    \item To derive the constant upper bound for $\sup_{x\in[a,b]} |b_k|$ and $\sup_{x\in[a,b]} |\mathbb{E} \beta_k|$ (on page 306 of \cite{liu2022ieee}), we need different strategy since we cannot use $\mathbb{E}|\ba_k^* \bQ_k \ba_k|^2 = o(1)$. 
    It has been established that $b_k' = \sup_{x\in[a,b]}\left[1+n^{-1} \mathbb{E}\trace \bR \bQ_k'\right]^{-1}$ is bounded from above. 
    This implies that there exists a positive constant $c_0$ such that $|1+n^{-1} \mathbb{E}\trace \bR \bQ_k'|>c_0$.
    It can be checked through martingale expansion that $n^{-1}\mathbb{E}\trace [\bSigma (\bQ_k-\bQ_k')]|=o(1)$. Thus, together with the fact that $\trace(\bA\bA^*) = O(n^{1/3})$, and $\|\bQ\|,\|\bQ'\| \le v_n^{-1}$ we know that both
    $|1+n^{-1} \mathbb{E}\trace \bR \bQ_k|$ and $|1+n^{-1} \mathbb{E}\trace \bSigma \bQ_k|$ are also bounded from below, and still use $c_0$ to denote the common lower bound.
    
    Now to bound $|b_k|$ for $x\in[a,b]$, we consider two cases on the bound of $\ba_k^* \mathbb{E} \bQ_k \ba_k$.  If $|\ba_k^*\mathbb{E} \bQ_k \ba_k|<c_0/2$, then we have the trivial bound $|b_k|<2/c_0$. If $|\ba_k^*\mathbb{E} \bQ_k \ba_k|\ge c_0/2$, we use Lemma S.1 to write 
    \begin{equation}\label{akQ}
    \ba_k^* \bQ \ba_k = [\ba_k^*  \bQ_k \ba_k -\beta_k \ba_k^* \bQ_k \bx_k \bx_k^* \bQ_k \ba_k] =  \ba_k^* \bQ_k \ba_k(1-\beta_k \ba_k^* \bQ_k \ba_k)+o(1),
    \end{equation}
    where the equality holds with overwhelming probability (An event $\Omega$ holds with overwhelming probability if $\mathbb{P}(\Omega)>1-n^{-\ell}$ for any positive constant $\ell$). 
    The left term is bounded from above by a constant with overwhelming probability, since $\mathbb{E} |\bv^* \bQ \bv- \bv^* T \bv|^q = O(n^{-q/2} v_n^{-q})$ for any fixed vector $\bv$ of bounded norm and any constant $q>2$, and $\|T\|<\infty$ when $x\in[a,b]$.
    Note that $(1-\beta_k \ba_k^* \bQ_k \ba_k) = \beta_k [1+\bw_k^* \bSigma \bQ_k \bSigma^{1/2}\bw_k+o(1)]  = \beta_k [1+ n^{-1} \mathbb{E} \trace \bSigma \bQ_k +o(1)] $. Recall $|1+n^{-1} \mathbb{E}\trace \bSigma \bQ_k|>c_0$ and $|\ba_k^* \bQ_k \ba_k| \ge c_0/3$.
    According to \eqref{akQ} and these bounds, we know that for $x\in[a,b]$, $|\beta_k|$ is bounded from above with overwhelming probability. This together with the naive bound $|\beta_k|< |z|/v_n$ implies that $\sup_{x\in[a,b]}\mathbb{E}|\beta_k|^q < C_q$ for some constant $C_q$ depending only on $q$. By the relation $b_k = \beta_k - b_k \beta_k \gamma_k$, it can be readily checked that $|b_k|< C$ for some constant $C$.
    
    \item In the last step of (35), we can use \begin{equation*}\begin{aligned}|(1+n^{-1} \trace \bSigma \mathbb{E}\bQ)^{-1}+z\mathbb{E} \underline{m}_n | &= \frac{1}{n}\sum_{k=1}^n \mathbb{E} \beta_k \alpha_n (\bx_k^*\bQ_k \bx_k-\frac{1}{n}\trace \bSigma \mathbb{E}\bQ) \\ & \le n^{-1} \sum_{k=1}^n \|\ba_k\|^2 \|\bQ_k\| +O(n^{-1/2}) = O(n^{-1/2}) 
    \end{aligned}\end{equation*}
    where the second last step uses $\mathbb{E}|\beta_k|^q<C_q$, $\alpha_n < C$ and $\mathbb{E} |\bw_k^* \bSigma^{1/2}\bQ_k \bSigma^{1/2}\bw_k - n^{-1}\trace \bSigma \mathbb{E}\bQ |^q= O(n^{-q/2})$, and the last step uses $\trace(\bA\bA^*) = O(n^{1/3})$. 
    
  % \item The constant upper bound of $\sup_{x\in[a,b]} |E \beta_k|$ can be checked readily using the new bound on $\ba_k$.
\end{enumerate}
There are some other terms that can be bounded using very similar strategies, and we omit the details.

Then, we show the second step: the exact separation phenomenon holds under Gaussian $\bW$. This can be further divided into two steps based on the assumption on $\bA$: 1. $\bA$ contains $K$ different columns; 2. general $\bA$ satisfying Assumption \ref{ass3new}.
To conclude the first step here when $\bA$ contains $K$ different columns, we just need to modify step proving equation (46) in \cite{liu2022ieee}, where our previous analysis requires $\trace (\bA^*\bA) =O(1)$ but this does not hold if $\bA$ does not have fixed rank. 
To handle this, note that for the Gaussian case, $\bA\bW^* \bW \bA^*$ is a Wishart matrix with distribution $W_p(n^{-1}\bA\bA^*,p)$. Then it is easy to conclude that if the variance of the entries of Gaussian $\bW$ is normalized as $\mathbb{E}|W_{ij}|^2=n^{-1}$, $\|\bA\|<\infty$, and $p/n\to 0$, the spectral norm of $\bA\bW^* \bW \bA^*$ tends to zero almost surely, so as that of $\bW\bA^*$.   Other parts of the proof of Theorem 1 can be extended to the case that $\bA$ contains $K=O(n^{1/3})$ distinct column vectors straightforwardly. In this way, we extend their conclusion that Theorem 1 holds with probability one for all large $n$, instead of probability tending to one.

To handle the case for general $\bA$ satisfying Assumption \ref{ass3new} and Gaussian $\bW$, we take use of the invariance of distribution of $\bW$ under orthogonal matrices. Assume that $\bW$ is Gaussian, and $\bA$ has a singular value decomposition $\bU_1\boldsymbol\Lambda \bV_1^*$, where $\boldsymbol\Lambda$ is a $p\times n$ matrix with $K$ singular values on its first $K$ main diagonal positions for some $K=O(n^{1/3})$.  Let $\bV_2$ be an $n\times n$ orthogonal matrix where the first $K$ row is constructed such that the number of distinct columns is $K$. One can construct such $\bV_2$ by letting the first $\lfloor n/K \rfloor$ entries of the first row be $\sqrt{\lfloor n/K \rfloor ^{-1}}$, other positions being zero; the $\lfloor n/K \rfloor+1$ to $2 \lfloor n/K \rfloor$ entries of the second row be $\sqrt{\lfloor n/K \rfloor ^{-1}}$, other entries being zero;  and  the same way to define the $K$-th row.
Further define $\bO = \bV_1 \bV_2$. Then $\bU_1^* \bX \bO = \boldsymbol\Lambda \bV_2 + \bU^*\bSigma^{1/2} \bW \bO \stackrel{d}=\boldsymbol\Lambda \bV_2 + \bU^* \bSigma^{1/2} \bW  $ becomes Model 1 in \cite{liu2022ieee}, i.e., the columns of signal part contains $K$ different vectors. Therefore the conclusion of this lemma holds for $\bU_1^* \bX \bO \bO^* \bX^* \bU_1$. Thus the conclusion also holds for $\bX\bX^* $ since $\bU_1$ and $\bO$ are both orthogonal matrices.

To conclude the third step,  we introduce a continuous interpolation matrix defined as $\bW(t)=\sqrt{t}\bW_1 + \sqrt{1-t}\bW_0$ for $t\in [0,1]$, where $\bW_0$ is Gaussian, $\bW_1$ is general, and both satisfy the moment conditions in Assumption \ref{ass1}. Note that $\bW(t)$ satisfies Assumption \ref{ass1} for any $t\in [0,1]$. Define $\bX(t)$ and $\bS(t)$  by replacing $\bW$ with $\bW(t)$, respectively. Denote the $i$-th largest singular value of a matrix $M$ by $\sigma_i(M)$. For any $t_1, t_2 \in[0,1]$, we have
\begin{equation}\label{continuity}\begin{aligned}
   | \lambda_i (\bS(t_1)) - \lambda_i(\bS(t_2))|&\leq C' |\sigma_i(\bX(t_1))- \sigma_i(\bX(t_2)) | \\& \leq C'' \sigma_1 (\bW(t_1)-\bW(t_2)) \leq C''' \sqrt{|t_1-t_2|},
\end{aligned}\end{equation}
where $C', C'', C'''$ are some positive constants independent of $n, t_1, t_2$. In the first and third step we use the fact that $\sigma_1(\bS(t))$ and  $\sigma_1(\bW(t))$ are bounded, and the second step uses Wely's inequality.
   Now we can conclude the exact separation by the continuity of eigenvalues together with the general version of Proposition 1 in \cite{liu2022ieee}, established in step 1. More specifically, let $t_j = j/n$, we know that $\lambda_k(\bS(0))>b_k$, and the extended Proposition 1 in \cite{liu2022ieee} implies that there are no eigenvalues of $\bS(t_j)$, $j=1,\cdots, n$ in $[a_k, b_k]$. Therefore \eqref{continuity} implies $\lambda_k(\bS(1))>b_k$ almost surely for sufficiently large $n$.
\qed

\vspace{1em}
\noindent\textbf{Proof of Proposition \ref{cenprop}.}
We first consider the bound involving $\tilde{\mathcal{\bQ}}_n(z)$.
The strategy is to consider the difference between the quadratic form involving $\tilde{\mathcal{\bQ}}_n(z)$ and the one involving $\tilde{\bQ}_n$, as studied in Proposition \ref{prop1}, and to derive the limits of such difference terms. The approach also applies to deriving the second bound for $\mathcal{\bQ}_n(z)$. 

Similarly to \eqref{woodburrytQ}, by Lemma \ref{wbm}, we have
\begin{equation*}
    (\boldsymbol\Phi \bX^*\bX \boldsymbol\Phi-z\bI)^{-1} = (-z\bI)^{-1}+z^{-1}\boldsymbol\Phi \bX^*(\bX\boldsymbol\Phi \bX^*-z\bI)^{-1}\bX \boldsymbol\Phi.
\end{equation*}
%By \eqref{woodburrytQ}
%\begin{equation*}   \bu^*\boldsymbol\Phi \bX^*(\bX \bX^*- z \bI)^{-1} \bX \boldsymbol\Phi \bu =  \bu^*\boldsymbol\Phi (\bX^* \bX-z\bI)^{-1}\boldsymbol\Phi \bu +\bu^*\boldsymbol\Phi \bu.\end{equation*}
%According to ..., 
%\begin{equation}
%   \bu^*\boldsymbol\Phi (\bX^* \bX-z\bI)^{-1}\boldsymbol\Phi \bu -\left[  \tilde{r} \bu^* \boldsymbol\Phi \bu -\tilde{r}^2 \bu^*\boldsymbol\Phi \bA_n^*(\bI+\tilde{r}\bR)^{-1}\bA_n\boldsymbol\Phi \bu\right]
%\end{equation}
By calculations, we find 
\begin{equation}\label{diftQQ1}
\begin{aligned}
   \tilde{D}_c&:= \bu^*\boldsymbol\Phi \bX^*(\bX\boldsymbol\Phi \bX^*- z \bI)^{-1} \bX \boldsymbol\Phi \bu - \bu^*\boldsymbol\Phi \bX^*(\bX \bX^*- z \bI)^{-1} \bX \boldsymbol\Phi \bu \\
    & = \bu^*\boldsymbol\Phi \bX^*(\bX \bX^*-z\bI)^{-1} n^{-1}\bX \mathbf{1}\mathbf{1}^*\bX^*(\bX\boldsymbol\Phi \bX^*-z\bI)^{-1}\bX \boldsymbol\Phi \bu\\
    &= -\frac{n^{-1}\bu^* \boldsymbol\Phi \bX^* (\bX \bX^*- z\bI)^{-1} \bX \mathbf{1}\mathbf{1}^* \bX^*(\bX \bX^*-z\bI)^{-1}\bX\boldsymbol\Phi \bu}{1-n^{-1}\mathbf{1}^* \bX^*(\bX \bX^*-z\bI)^{-1}\bX \mathbf{1}}\\
    & = \frac{n^{-1}z^2\bu^* \boldsymbol\Phi (\bX^*\bX-z\bI)^{-1}\mathbf{1}\mathbf{1}^*(\bX^*\bX-z\bI)^{-1} \boldsymbol\Phi \bu}{zn^{-1}\mathbf{1}^*(\bX^*\bX-z\bI)^{-1}\mathbf{1}}
\end{aligned}
\end{equation}
where the first step uses \eqref{eq04}, the second step uses \eqref{eq03}, and the third step uses \eqref{woodburrytQ} and $\bu^*\boldsymbol\Phi \mathbf{1}=0$.
Following similar steps, we obtain
\begin{equation}\label{diftQQlim}
\begin{aligned}
\tilde{L}_c&:=\tilde{r}^2 \bu^*\boldsymbol\Phi \bA_n^* (\bI+\tilde{r}\bar{\bR}_n)^{-1}\bA_n\boldsymbol\Phi \bu- \tilde{r}^2\bu^*\boldsymbol\Phi \bA_n^* (\bI+\tilde{r}{\bR}_n)^{-1}\bA_n\boldsymbol\Phi \bu \\
& =\frac{\tilde{r}^2\bu^*\boldsymbol\Phi \bA_n^*(\bI+\tilde{r}\bR)^{-1}n^{-1}\bA_n \tilde{r}\mathbf{1}\mathbf{1}^* \bA_n^*(\bI+\tilde{r}\bR)^{-1}}{1-\tilde{r}n^{-1}\mathbf{1}^* \bA_n^* (\bI+\tilde{r}\bR)^{-1}\bA_n\mathbf{1}}.
\end{aligned}
\end{equation}

Next we verify that $z^{-1}\tilde{D}_c-\tilde{L}_c$ is $O_P(n^{-1/2})$.
By polarization, 
\eqref{key1} still holds with different sequences of deterministic vectors on both sides of $\tilde{\bQ}_n(z)-\tilde{R}_n(z)$. This together with $\bu^*\boldsymbol\Phi \mathbf{1}=0$ yields  $$n^{-1/2}\bu^* \boldsymbol\Phi (\bX^*\bX-z\bI)^{-1}\mathbf{1} + n^{-1/2}\tilde{r}^2 \bu^*\boldsymbol\Phi \bA_n^*(\bI+\tilde{r}\bR_n)^{-1}\bA_n \mathbf{1} = O_P(n^{-1/2}).$$
For the term in the denominator, we have
\begin{equation*}
    zn^{-1}\mathbf{1}^*(\bX^*\bX-z\bI)^{-1}\mathbf{1} - \left[z\tilde{r}- zn^{-1}\tilde{r}^2 \mathbf{1}^* \bA_n^*(\bI+\tilde{r}\bR)^{-1}\bA_n \mathbf{1} \right] = O_P(n^{-1/2}).
\end{equation*}
With these two bounds, and by calculating the difference of those two terms in the last step of \eqref{diftQQ1} and \eqref{diftQQlim}, respectively,
%by dividing $z$ on both sides of \eqref{diftQQ1} and then considering the difference with both sides of \eqref{diftQQlim},
we find that 
\begin{equation}\label{difDcLc}
\begin{aligned}
& z^{-1}\tilde{D}_c-\tilde{L}_c= z^{-1}\left[\bu^*\boldsymbol\Phi \bX^*(\bX\boldsymbol\Phi \bX^*- z \bI)^{-1} \bX \boldsymbol\Phi \bu - \bu^*\boldsymbol\Phi \bX^*(\bX \bX^*- z \bI)^{-1} \bX \boldsymbol\Phi \bu \right]\\ & \quad - \left[\tilde{r}^2 \bu^*\boldsymbol\Phi \bA_n^* (\bI+\tilde{r}\bar{\bR}_n)^{-1}\bA_n\boldsymbol\Phi \bu- \tilde{r}^2\bu^*\boldsymbol\Phi \bA_n^* (\bI+\tilde{r}{\bR}_n)^{-1}\bA_n\boldsymbol\Phi \bu\right] = O_P(n^{-1/2}).
\end{aligned}
\end{equation}

Since $(-z\bI)^{-1}\bu^*\boldsymbol\Phi \bu+z^{-1} \bu^*\boldsymbol\Phi \bX^*(\bX \bX^*- z \bI)^{-1} \bX \boldsymbol\Phi \bu - \left[\tilde{r}\bu^*\boldsymbol\Phi \bu-\tilde{r}^2 \bu^* \boldsymbol\Phi \bA_n^* (\bI+\tilde{r}{\bR}_n)^{-1}\bA_n\boldsymbol\Phi \bu \right]=O_P(n^{-1/2})$ by \eqref{key1}, we conclude from this and \eqref{difDcLc} that  
\begin{equation*}\begin{aligned}
&(-z\bI)^{-1}\bu^*\boldsymbol\Phi \bu+z^{-1}\bu^*\boldsymbol\Phi \bX^*(\bX\boldsymbol\Phi \bX^*- z \bI)^{-1} \bX \boldsymbol\Phi \bu\\
& \quad \quad - \left[\tilde{r}\bu^*\boldsymbol\Phi \bu- \tilde{r}^2 \bu^*\boldsymbol\Phi \bA_n^* (\bI+\tilde{r}\bar{\bR}_n)^{-1}\bA_n\boldsymbol\Phi \bu \right]=O_P(n^{-1/2}).\end{aligned}\end{equation*}

Therefore, 
\begin{equation*}
    \bu^*(\boldsymbol\Phi \bX^*\bX\boldsymbol\Phi-z\bI)^{-1} - \left[\tilde{r}\bu^*\boldsymbol\Phi \bu - \tilde{r}^2 \bu^*\boldsymbol\Phi \bA_n^* (\bI+\tilde{r}\bar{\bR}_n)^{-1}\bA_n\boldsymbol\Phi \bu -z^{-1}\bu^* n^{-1}\mathbf{1}\mathbf{1}^* \bu\right] = O_P(n^{-1/2}).
\end{equation*}
This concludes the bound for $\tilde{\mathcal{\bQ}}_n(z)$.
%Let $\tilde{r}_c$ be the solution to \eqref{32ra3} with $\bR_n$ replaced by $\bar{\bR}_n$. Following an argument similar to the one that lead to \eqref{difdeltatr}, we obtain that $|\tilde{r}_c -\tilde{r}| = O(n^{-1})$ if $\Im z$ is bounded away from zero. 
%Therefore we conclude the first bound on $\tilde{\bQ}_n(z)$.

Then we prove the second bound on $\mathcal{\bQ}_n(z)$.
We have
\begin{equation}\label{Dc}
\begin{aligned}
   D_c &= \bv^*\left[(\bX \boldsymbol\Phi \bX^*- z\bI)^{-1} - (\bX \bX^*-z\bI)^{-1}\right]\bv\\
    & =\frac{\bv^*\bQ_n n^{-1}\bX\mathbf{1}\mathbf{1}^*\bX^* \bQ_n\bv}{1-n^{-1}\mathbf{1}^*\bX^* \bQ \bX \mathbf{1}} = \frac{\bv^*\bQ_n n^{-1}\bA\mathbf{1}\mathbf{1}^*\bA^* \bQ_n\bv}{1-n^{-1}\mathbf{1}^*\bX^* \bQ \bX \mathbf{1}} + O_P(n^{-1/2}),
\end{aligned}
\end{equation}
where in the last step we use $n^{-1/2}\bv^*\bQ_n \bW \mathbf{1} = O_P(n^{-1/2})$, which can be checked by following arguments similar to (3.7)-(3.12) of \cite{pan2014comparison}.
We also find
\begin{equation*}
\begin{aligned}
   L_c  &:= (-z - z \tilde{r}\bar{\bR}_n)^{-1} -(-z - z \tilde{r}{\bR}_n)^{-1} \\
   &= \frac{-(-z - z\tilde{r} \bR)^{-1}z\tilde{r}n^{-1}\bA \mathbf{1}\mathbf{1}^* \bA^* (-z - z\tilde{r} \bR)^{-1}}{1-z n^{-1}\mathbf{1}^*\bA^* (-z - z\tilde{r} \bR)^{-1}\bA \mathbf{1}}.
\end{aligned}
\end{equation*}
By \eqref{Qnlim}, we have $n^{-1/2}\bv^*\bQ_n \bA\mathbf{1}-n^{-1/2}\bv^*(-z - z\tilde{r} \bR)^{-1}\bA\mathbf{1}=O_P(n^{-1/2}).$ Since $1-n^{-1}\mathbf{1}^*\bX^* \bQ \bX \mathbf{1} =zn^{-1}\mathbf{1}^* \tilde{\bQ}_n \mathbf{1},$
according to \eqref{key1}, we find $1-n^{-1}\mathbf{1}^*\bX^* \bQ \bX \mathbf{1} +zn^{-1}\mathbf{1}^*\left[\tilde{r}\bI+\tilde{r}^2 \bA^*(\bI+\tilde{r}(z)\bR_n)^{-1} \bA\right] \mathbf{1}=O_P(n^{-1/2})$. Therefore, we obtain that $D_c-L_c = O_P(n^{-1/2})$. This combining with the fact that $\bv^*(\bX\bX^*-z\bI)^{-1}\bv -\bv^*(-z - z \tilde{r}{\bR}_n)^{-1}\bv=O_P(n^{-1/2})$ 
concludes the proof.
\qed

\vspace{1em}
\noindent\textbf{Proof of Theorem \ref{centhm3}.}
We can follow those arguments used in the proof of Theorems \ref{thmeve1} and \ref{thm2} to conclude the result. 
Note that the term $z^{-1} n^{-1}\mathbf{1}\mathbf{1}^*$ appeared in $\boldsymbol{\tilde{\mathcal{D}}}_n$ does not contribute to the integral as in \eqref{defF1} and \eqref{defF2} since the contour of the integration does not enclose the origin. 
\qed

{\section*{Acknowledgements}
The authors thank the Editor, Associate Editor, and two referees for their valuable comments that improved the paper. The authors are listed alphabetically and have contributed equally to this work.
}

\section*{Funding}
Xiaoyu Liu was partially supported by the Basic Research Program of Guangzhou Municipal Science and Technology Bureau(Grant No, SL2022A04J01237). Yiming Liu was partially supported by grants 23JNQN19 of Fundamental Research Funds
for the Central University and 12301340 of NSFC, 2024A1515012761 of Natural Science Foundation of Guangdong Province, China, SL2023A04J00604 of Basic and Applied Basic Research Foundation, 2024.
G.M. Pan was partially supported by MOE Tier 1 grant RG76/21 and MOE-T2EP20123-0007 at NTU.
Lingyue Zhang was supported by the National Natural Science Foundation of China under Grant 12201435.
Zhixiang Zhang was partially supported by the SRG2023-00053-FST from University of Macau.

\bibliographystyle{apalike}
\bibliography{clu-ref2}

\end{document}